\documentclass[a4paper,11pt]{amsart}
\usepackage[cp1250]{inputenc}
\usepackage{amstext,amsbsy,amsmath,amscd,amssymb, graphics}
\usepackage[all]{xy}
\input{xy}

\let\ssize\scriptstyle
\newif\ifFIRST\newdimen\MAXright\MAXright0pt
\def\sdynkin{\bgroup\eightpoint\dynkin}
\def\endsdynkin{\enddynkin\egroup}
\def\dynkin{\bgroup\FIRSTtrue\hskip.5em\setbox1\hbox{$\diagup$}%
\setbox2\hbox{$\diagdown$}%
\setbox0\hbox to2\wd1{\hrulefill}%
\setbox3\hbox{$\bullet$}%
\setbox4\hbox{$\times$}%
\setbox7\hbox{$\circ$}
\def\whiteroot##1{\ifFIRST\setbox5\hbox{$##1$}\ifdim\wd5>1.3em
\hskip-.5em\hskip.5\wd5\fi\fi\FIRSTfalse
\hskip-.25em\raise1.5\wd3\hbox to0pt{\hss\hskip.45em$
\ssize##1$\hss}\copy7\hskip-.25em\setbox6\hbox{$##1$}
\MAXright\wd6}
\def\root##1{\ifFIRST\setbox5\hbox{$##1$}\ifdim\wd5>1.3em%
\hskip-.5em\hskip.5\wd5\fi\fi\FIRSTfalse%
\hskip-.25em\raise1.5\wd3\hbox to0pt{\hss\hskip.45em$%
\ssize##1$\hss}\copy3\hskip-.25em\setbox6\hbox{$##1$}%
\MAXright\wd6}%
\def\whitedroot##1{\ifFIRST\setbox5\hbox{$##1$}\ifdim\wd5>1.3em
\hskip-.5em\hskip.5\wd5\fi\fi\FIRSTfalse
\hskip-.25em\lower1.8\wd3\hbox to0pt{\hss\hskip.45em$
\ssize##1$\hss}\copy7\hskip-.25em\setbox6\hbox{$##1$}
\MAXright\wd6}%
\def\whiterroot##1{\hskip-.25em\copy7\hbox to0pt{\hskip.3em$\ssize##1$\hss}%
\hskip-.25em\setbox6\hbox{\hskip.6em$##1##1$}%
\MAXright\wd6}%
\def\droot##1{\ifFIRST\setbox5\hbox{$##1$}\ifdim\wd5>1.3em%
\hskip-.5em\hskip.5\wd5\fi\fi\FIRSTfalse%
\hskip-.25em\lower1.8\wd3\hbox to0pt{\hss\hskip.45em$%
\ssize##1$\hss}\copy3\hskip-.25em\setbox6\hbox{$##1$}%
\MAXright\wd6}%
\def\rroot##1{\hskip-.25em\copy3\hbox to0pt{\hskip.3em$\ssize##1$\hss}%
\hskip-.25em\setbox6\hbox{\hskip.6em$##1##1$}%
\MAXright\wd6}%
\def\norroot##1{\hskip-.36em\copy4\hbox to0pt{\hskip.3em$\ssize##1$\hss}%
\hskip-.48em\setbox6\hbox{\hskip.6em$##1##1$}%
\MAXright\wd6}%
\def\noroot##1{\ifFIRST\setbox5\hbox{$##1$}\ifdim\wd5>1.3em%
\hskip-.5em\hskip.5\wd5\fi\fi\FIRSTfalse%
\hskip-.36em\raise1.5\wd3\hbox to0pt{\hss\hskip.6em$%
\ssize##1$\hss}\copy4\hskip-.38em\setbox6\hbox{$##1$}%
\MAXright\wd6}%
\def\nodroot##1{\ifFIRST\setbox5\hbox{$##1$}\ifdim\wd5>1.3em%
\hskip-.5em\hskip.5\wd5\fi\fi\FIRSTfalse%
\hskip-.36em\lower1.8\wd3\hbox to0pt{\hss\hskip.6em$%
\ssize##1$\hss}\copy4\hskip-.38em\setbox6\hbox{$##1$}%
\MAXright\wd6}%
\def\nolink{\hskip\wd0}
\def\link{\raise.22em\copy0}%
\def\llink##1{\raise.32em\copy0\hskip-\wd0%
\raise.12em\copy0\hskip-.5\wd0\hbox to0pt{\hss$##1$\hss}\hskip.5\wd0}%
\def\lllink##1{\raise.22em\copy0\hskip-\wd0\raise.32em\copy0\hskip-\wd0%
\raise.12em\copy0\hskip-.5\wd0\hbox to0pt{\hss$##1$\hss}\hskip.5\wd0}%
\def\rootupright##1{\hbox to0pt{\raise.45em\copy1\hskip-.25em\raise1.3\ht1%
\hbox{\copy3\hskip.3em$\ssize##1$}\hss}%
\setbox6\hbox{\hskip.6em\copy1\copy1$##1##1$}%
\ifdim\MAXright<\wd6\MAXright\wd6\fi}%
\def\whiterootupright##1{\hbox to0pt{\raise.45em\copy1\hskip-.25em\raise1.3\ht1
\hbox{\copy7\hskip.3em$\ssize##1$}\hss}
\setbox6\hbox{\hskip.6em\copy1\copy1$##1##1$}
\ifdim\MAXright<\wd6\MAXright\wd6\fi}
\def\norootupright##1{\hbox to0pt{\raise.45em\copy1\hskip-.36em\raise1.3\ht1%
\hbox{\copy4\hskip.3em$\ssize##1$}\hss}%
\setbox6\hbox{\hskip.6em\copy1\copy1$##1##1$}%
\ifdim\MAXright<\wd6\MAXright\wd6\fi}%
\def\rootdownright##1{\hbox to0pt{\raise-.5em\copy2\hskip-.25em\raise-1.35\ht1%
\hbox{\copy3\hskip.3em$\ssize##1$}\hss}\setbox6%
\hbox{\hskip.6em\copy2\copy2$##1##1$}%
\ifdim\MAXright<\wd6\MAXright\wd6\fi}%
\def\whiterootdownright##1{\hbox to0pt{\raise-.5em\copy2\hskip-.25em\raise-1.35\ht1
\hbox{\copy7\hskip.3em$\ssize##1$}\hss}\setbox6
\hbox{\hskip.6em\copy2\copy2$##1##1$}
\ifdim\MAXright<\wd6\MAXright\wd6\fi}
\def\rootdown##1{\hbox to0pt{\hskip-.05em\vrule height.25em depth.65em%
\hskip-.25em\raise-.95em\hbox{\copy3\hskip.3em$\ssize##1$}\hss}%
\setbox6\hbox{$##1$}%
\ifdim\MAXright<\wd6\MAXright\wd6\fi}%
\def\whiterootdown##1{\hbox to0pt{\hskip-.05em\vrule height.25em depth.65em
\hskip-.25em\raise-.95em\hbox{\copy7\hskip.3em$\ssize##1$}\hss}
\setbox6\hbox{$##1$}
\ifdim\MAXright<\wd6\MAXright\wd6\fi}
\def\dots{\hskip.5em\cdots\hskip.5em}}%
\def\enddynkin{\ifdim\MAXright>1em\hskip.5\MAXright\else\hskip.5em\fi\egroup}%
 
This macro is used for creation of Dynkin diagrams in AMSTeX.
The usage:
In math-mode, between \dynkin and \enddynkin,
\root#1 makes a node with label #1 over this node,
\rroot#1 and \droot do the same with labels on right and down,
\noroot and \norroot do the same with cross instead of bullet,
\link links two adjacent \root's, \llink#1 links two adjacent nodes
with double-line with #1 (probably '>' or '<') in the middle, \lllink is does
same with three lines, \dots does \cdots between the two adjacent nodes,
\rootdown#1 makes a linked node down labeled by #1, \rootupright#1 and
\rootdownright#1 are labeled linked nodes in the indicated directions 
which do not count for the jorizontal dimensions.

\sdynkin ... \endsdynkin works exactly as above but in smaller
version (eightpoint)

Examples:
$\dynkin \root{a_1}\link\root{a_2}\dots\root{a_{n-1}}\link\root{a_n}
\enddynkin$
or
$\dynkin \root{}\lllink>\root{}\enddynkin$
or
$\dynkin \root{a}\link\root{b}\rootupright{c}\rootdownright{d}\enddynkin$
or
$\dynkin \root{}\link\root{}\rootupright{}\link\root{}\enddynkin$.
or (white roots!)
{$\dynkin \whiteroot{0}\link\noroot{0}\dots\root{0}\link\whiteroot{0}\rootupright{0}\whiterootdownright{0}\enddynkin$}. 

\advance\topmargin by -0.6cm

\numberwithin{equation}{section}

\newtheorem{Theo}{Theorem}[subsection]
\newtheorem{Cor}[Theo]{Corollary}
\newtheorem{Lem}[Theo]{Lemma}
\newtheorem{Exa}[Theo]{Example}

\newtheorem{Def}[Theo]{Definition}
\newtheorem{Rem}[Theo]{Remark}
\newtheorem{corollary}[Theo]{Corollary}
\newtheorem{hypothesis}[Theo]{Conjecture}

\let\cal=\mathcal

\def\Hom{\mathop {\rm Hom} \nolimits}

\def\R{\mathbb R}
\def\C{\mathbb C}
\def\N{\mathbb N}

\def\Z{\mathbb Z}
\def\V{\mathbb V}
\def\W{\mathbb W}

\def\lieg{\mathfrak{g}}
\def\lieh{\mathfrak{h}}
\def\liep{\mathfrak{p}}
\def\lieb{\mathfrak{b}}
\def\lien{\mathfrak{n}}
\def\liez{\mathfrak{z}}
\def\univ{{\cal U}}
\def\so{{\mathfrak{so}}}
\def\SO{{\mathrm{SO}}}
\def\Euc{{\mathrm{Euc}}}
\def\Spin{{\mathrm{Spin}}}
\def\sl{{\mathfrak{sl}}}
\def\gl{{\mathfrak{gl}}}

\def\Sl{{\mathrm{SL}}}
\def\So{{\mathrm{SO}}}

\def\root{\mathrm{root}}
\def\and{\mathrm{and}}

\def\SL{\mathrm{SL}}

\def\Tr{\mathrm{Tr}}
\def\diag{\mathrm{diag}}
\def\Cliff{\mathrm{Cliff}}
\def\Ad{{\mathrm{Ad}}}
\def\ad{{\mathrm{ad}}}

\def\half{{\frac{1}{2}}}
\def\gvmhom{$M_\liep(\mu)\to M_\liep(\lambda)$ }

\def\cV{{\cal V}}

\def\cD{{\cal D}}

\def\br#1{\langle #1 \rangle}

\def\wr{\whiteroot{}}

\def\mmmm{\hspace{-15pt}}

\begin{document}
{\thispagestyle{empty}
\pagestyle{empty}
\begin{center} 
\vskip 3cm
Doctoral dissertational thesis
\vskip 2cm

{\Huge
Several Dirac operators in Parabolic Geometry
}
\vskip 1.4cm
{Peter Franek}
\vskip 1cm
supervisor: prof. RNDr. Vladim\'\i r Sou\v cek, DrSc.
\end{center}
\newpage
{\bf Acknowledgement:} I thank to my supervisor Vladim\'\i r Sou\v cek for directing 
my study and research. Further, I thank to Svatopluk Kr\'ysl for his permanent help.

\newpage
\begin{center}
{\Large Table of Contents}
\end{center}
\smallskip
\contentsline {section}{\tocsection {}{1}{\bf Preface}}{1}
\contentsline {section}{\tocsection {}{2}{\bf Introduction}}{7}
\contentsline {subsection}{\quad\tocsubsection {}{2.1}{Semisimple Lie algebras}\dotfill}{7}
\contentsline {subsection}{\quad\tocsubsection {}{2.2}{Parabolic subalgebras and grading}\dotfill}{9}
\contentsline {subsection}{\quad\tocsubsection {}{2.3}{True and generalized Verma modules}\dotfill}{12}
\contentsline {subsection}{\quad\tocsubsection {}{2.4}{Invariant differential operators}\dotfill}{13}
\contentsline {subsection}{\quad\tocsubsection {}{2.5}{Dirac operator and Clifford algebras}\dotfill}{18}
\contentsline {section}{\tocsection {}{3}{\bf Homomorphisms of generalized Verma modules}}{20}
\contentsline {subsection}{\quad\tocsubsection {}{3.1}{True Verma module homomorphism}\dotfill}{20}
\contentsline {subsection}{\quad\tocsubsection {}{3.2}{Parabolic Hasse graph}\dotfill}{22}
\contentsline {subsection}{\quad\tocsubsection {}{3.3}{Standard and nonstandard homomorphism}\dotfill}{23}
\contentsline {subsection}{\quad\tocsubsection {}{3.4}{BGG graph}\dotfill}{26}
\contentsline {section}{\tocsection {}{4}{\bf Dirac operator in orthogonal parabolic geometry}}{31}
\contentsline {subsection}{\quad\tocsubsection {}{4.1}{Singular orbit of a particular weight}\dotfill}{31}
\contentsline {subsection}{\quad\tocsubsection {}{4.2}{The real version}\dotfill}{32}
\contentsline {subsection}{\quad\tocsubsection {}{4.3}{Description of the differential operator}\dotfill}{33}
\contentsline {subsection}{\quad\tocsubsection {}{4.4}{More Dirac operators}\dotfill}{37}
\contentsline {subsection}{\quad\tocsubsection {}{4.5}{Description of the operator}\dotfill}{38}
\contentsline {section}{\tocsection {}{5}{\bf Singular orbit corresponding to $k$ Dirac operators}}{41}
\contentsline {subsection}{\quad\tocsubsection {}{5.1}{Even dimension, $k\leq n$}\dotfill}{41}
\contentsline {subsection}{\quad\tocsubsection {}{5.2}{Even dimension, $k>n$}\dotfill}{48}
\contentsline {subsection}{\quad\tocsubsection {}{5.3}{Odd dimension}\dotfill}{52}
\contentsline {section}{\tocsection {}{6}{\bf Calculus of extremal vectors}}{59}
\contentsline {subsection}{\quad\tocsubsection {}{6.1}{Computation of extremal vectors}\dotfill}{59}
\contentsline {subsection}{\quad\tocsubsection {}{6.2}{Translation of the extremal vector to an operator}\dotfill}{75}
\contentsline {section}{\tocsection {}{}{\bf References}}{81}

\newpage
}

\setcounter{page}{1}
\section{Preface}

There are two basic generalizations of the space of holomorphic functions to higher dimensions. 
One of them is the notion of holomorphic functions in several variables, 
$f: \R^{2k}\simeq \C^k\to \C$, $\bar{\partial_j}f=0$
for $j=1,\ldots,k$. 


The second possible generalization deals with functions defined on $\R^n$ with values in the
{\it Clifford algebra} (a particular generalization of complex numbers). The functions in question are
solutions of a first order elliptic system of partial differential equations called {\it Dirac equation}, 
which is a higher dimensional analogue of Cauchy-Riemann equations. It is easy to describe
the system in dimension 4 using quaternionic notation (it was done by Fueter in 30's).
A hyperbolic analogue of the system in dimension 4 was discovered
in theoretical physics by P. Dirac. The (elliptic version of) Dirac equation has been used
extensively in mathematics during the last 40 years (e.g. \cite{Atiyah-Bott-Patodi}). 
Solutions of the Dirac are often called {\it monogenic functions} (or {\it harmonic spinors}). 
The function theory for them is now known under the
name {\it Clifford analysis} (\cite{BDS, Soucek-book}).
In explicit terms, monogenic functions $f$ are defined on real Euclidean space $\R^n$ with values
in the Clifford algebra $\Cliff(n,\R)$ (or the {\it space of spinors} $S$) such that $Df=0$, 
where $D=\sum_j e_j \cdot \partial_j$ is
given by multiplication in the Clifford algebra (or the Clifford multiplication $\R^n\otimes S\to S$). 
Monogenic functions have similar nice properties as holomorphic functions 
(Cauchy integral formula, theory of residues, analyticity, maximum principle, unique continuation property, etc.) 
and they coincide with usual holomorphic functions on $\R^2\simeq \C$ for $n=2$ (\cite{GM}). 

Eigenfunctions of the (hyperbolic version of) the {\it Dirac operator} $D$
describe spin $1/2$ particles with mass in relativistic quantum mechanics.
As all basic equations of relativistic physics, it is {\it invariant} with respect to
the Poincar\'e group. Similar nice invariance properties are true also for the elliptic version
of the Dirac equation. In particular,  its
solutions are invariant with respect to the group $\Spin(n)$. It means that if a 
function $f$ is monogenic, the same is true for the function $g\cdot f$ defined by
$(g\cdot f)(x)=g(f(g^{-1}\cdot x))$  
(here $\R^n$ and $S$ are considered as the fundamental defining and fundamental spinor 
representation of $\Spin(n)$). Moreover, it turns out that the symmetry group is much larger
then $\Spin(n)$. It is a group $G$ which is a double-cover of the group of all M\"obius transformations
of $\R^n$. It contains $\Spin(n)$ similarly as the group of M\"obius transformations in the plane contains 
rotations. 
This is an analogue of the fact that holomorphic functions
are preserved by conformal transformations.

It is a natural idea to consider monogenic functions of several Clifford variables,
which form a common
generalization of the space of holomorphic functions in several variables and 
of the space of monogenic functions in one variable. 
A {\it monogenic function of several Clifford variables} is a function $f: (\R^n)^k\to S$,
where $S$ is the spinor module over $\Cliff(n,\R)$ resp. $\Cliff(n,\C)$ such that
$D_i f=0$, where $D_i=\sum_j e_j\cdot \partial_{ij}$ for $i=1,\ldots,k$ 
($x_{uv}$ are variables on $(\R^n)^k$, $u=1,\ldots, k, v=1,\ldots,n$).
The whole system can be written as $Df=0$ where $D=(D_1,\ldots, D_k),$
is called {\it Dirac operator in several Clifford variables.} 

One important theorem in the theory of holomorphic functions of several complex variables 
is the Hartog theorem. It says that not all open sets are natural
domains of definition for holomorphic functions of several complex variables. There are domains $\Omega$
with the property that any function holomorphic on $\Omega$ can be extended holomorphically to a larger
domain. This was a completely new phenomenon, which is not true for one complex variable. 
This property is a consequence of the fact that the Cauchy-Riemann operator $D$ in higher dimensions defines an overdetermined system of PDE's. The Hartog type theorems can be systematically
studied using a resolution of $D$, i.e. the (locally exact) complex of PDE's starting with the operator $D.$
This is the Dolbeault sequence, which is nowadays a standard basic part of the theory of functions
of several complex variables.

It can be expected that the Dirac operator in several Clifford variables will also define
an overdetermined system of PDE's and that Hartog type theorem will hold for monogenic
functions in several variables as well. As in the complex case, an adequate tool for a study of
such properties of monogenic functions in several variables would be an analogue of the Dolbeault
sequence starting with the Dirac operator in several Clifford variables. 
It is not an easy task to find such a resolution
and its general form is still not known. However, many special cases of the problem are already understood. 

For some reason, the dimension $n=4$ is special 
and an analogy with complex analysis is stronger than in higher dimensions.
The elliptic version of the Dirac equation in dimension 4 (the Fueter equation)
was studied by Fueter already in 40's (\cite{Fueter}) and the resolution for it
is already well understood for any number of variables  
(\cite{Salamon,Baston,Struppa1, Struppa2, Struppa3, Somberg, Sabadini_Fabrizio_Bures,  Struppa_Soucek, Bures_Soucek}). 

In higher dimensions, the situation is more complicated. 
In \cite{bluebook}, the authors used 
Fourier transform and translated the problem into the language of commutative algebra. Instead of
sequences of differential operators, they were computing the resolution of a module over the
ring of polynomials by using the Hilbert syzygy theorem and Gr\"obner bases. 

The Dirac operator in $k$ variables is invariant with respect to the group $\Sl(k)\times \Spin(n)$ ($\C$
is considered to be the trivial and $\C^k$ the defining representation of $\Sl(k)$), 
similarly as the Dirac operator in one variable is $\Spin$-invariant.
We already mentioned that the usual Dirac operator in one variable is invariant
with respect to a group that is a double-cover of the group of M\"obius transformations.
It is shown in this work that a similar fact is true for the Dirac operator in several variables.
It is  invariant with respect to a larger group that contains $\Sl(k)\times \Spin(n)$ as its semisimple
subgroup (in a similar way  as the group of M\"obius transformations contains $\So(n)$ as it semisimple
subgroup). 

While the notion of holomorphic function can be extended to functions on complex manifolds,
the usual Dirac operator can be defined on {\it spin-manifolds} which are
manifolds with a given {\it spin structure} (see \cite{Friedrich}).
The Dirac operator
acts not on functions but rather on sections of the associated 
spinor bundle over such spin-manifolds. 
An example of a spin structure on the sphere is the projection
$\Spin(n+1)\to \Spin(n+1)/\Spin(n)$.  The curved version (in the Cartan sense)
of this are oriented Riemannian manifolds with a chosen spin structure and the Dirac operator acts
between sections of the vector bundle associated to the spinor representation of $\Spin(n)$.

It can be shown that the spin structure on the sphere described above is a reduction of 
the bundle $G\to G/P$ where $G=\Spin(n+1,1)$ and $P$ is the {\it parabolic subgroup}
fixing a line in the null-cone of the Minkowski metric. 
If  $\V$ and $\W$ are two spinor representations of $P$
(it means that we extend the spinor action of $\Spin(n)\subset P$ on these modules
by a suitable action of the center of $P$ and let the unipotent part of $P$ act trivially), 
the Dirac operator acts
between sections of $G\times_P \V$ and $G\times_P \W$ and is $G$-invariant ($G$ is the group of
invariance of the Dirac operator mentioned above).
The curved analogues of the homogenous bundle $G\to G/P$ are principal $P$-bundles
${\cal G}\to M$ over a manifold $M$ called {\it conformal spin structures} on $M$. 

This geometric structure on $M$ (together with a {\it Cartan connection} $\omega$ on ${\cal G}$
which is an analogue of the {\it Maurer-Cartan form} on $G$) is an example
of the so called {\it parabolic geometry}. These are geometries
modeled on a homogeneous space $(G,P)$, where $G$ is a semisimple
Lie group and $P$ a {\it parabolic subgroup}.
The choice of $P$ is equivalent to the choice of
a gradation $\lieg=\oplus_{i=-k}^k \lieg_i$ of $\lieg$. It is well known that conformal
geometries and projective geometries can be described as parabolic geometries with one-graded
Lie algebra $\lieg=\lieg_{-1}\oplus\lieg_0\oplus\lieg_1$, see \cite{Kobayashi72}.
The topic was
studied in details in e.g. \cite{Slovak_DrSc, Sharpe, Kobayashi-Nagano1, Kobayashi-Nagano2, Tanaka, FefG, Cap-Schichl},
and many properties of invariant operators on such manifolds are known 
(\cite{Css1, CD}).

 
The facts indicated above suggest that we may find
a suitable parabolic geometry, which would correspond well to the symmetry
of the Dirac operator in several variables. 
We show in this thesis that it is indeed possible. The corresponding couple
is the Lie group $\Spin(n+k,k)$ and its parabolic subgroup $P$ having 
$\Sl(k)\times \Spin(n)$ as its Levi factor. 
If we consider the $\SL(k)\times \Spin(n)$-spinors representations 
$\V\simeq \C\otimes S,\,\,\, \W\simeq \C^k\otimes S$ as irreducible $P$-modules (choosing a particular action of the center of $P$),  
there is a unique $G$-invariant differential operator $\Gamma(G\times_P \V)\to \Gamma(G\times_P \W)$. 
Further, we show that after suitable local identifications between sections of the vector bundles and functions on
a vector space, this operator really reduces to the Dirac operator in $k$ variables. 
Again, it is possible to define an appropriate curved version of the
operator $D$ and to study properties of solutions of $D$ on manifolds with a given parabolic
structure of type $(G,P).$
In this work, we shall not study these questions in a curved situation and restrict our
attention to the homogeneous model.

Methods used in the thesis are completely algebraic.
It is well-known that there is a duality 
between invariant differential operators on the flat model of a parabolic geometry and 
homomorphisms of so called {\it  generalized Verma modules}, which are
$\lieg$-modules dual to the space of infinite jets of sections of some associated vector bundle:
$$M_\liep(\V)\simeq (J^\infty_{eP}(\Gamma(G\times_P \V^*)))^*$$
On such modules, the action of $\lieg$, the Lie algebra of $G$, is defined naturally as the
derivative of the action of $G$ on sections. 
More precisely, due to the transitive action of $G$, it is possible to reduce the classification of
$G$-invariant differential operators between sections of homogeneous bundles
associated to $P$-modules $\V$ and $\W$ to a classification of $P$-invariant maps
between infinite jets of sections of the corresponding associated vector bundle at the origin.
The dual $P$-homomorphism can be always extended to a $(\lieg, P)$-homomorphism of the generalized Verma modules
$M_\liep(\W^*)\to M_\liep(\V^*)$ (and vice versa).

The theory of generalized Verma modules (further denoted simply by GVM)
was created mostly by Lepowski (\cite{Lepowski}) who generalized the results about true Verma modules of
Bernstein-Gelfand-Gelfand and Verma (\cite{bgg1, bgg2, Verma}).

For an irreducible finite dimensional $P$-module $\V$, there exist
only a finite sequence of homomorphisms of GVM's starting with $M_\liep(\V)$
and all the GVM's in this sequence are  induced by highest weights that are linked 
by the {\it affine action} of the Weyl group $W$ (associated to $\lieg$).
The homomorphisms of GVM's that have regular infinitesimal character
are quite well understood (\cite{Boe}). But the homomorphism dual to the Dirac operator in several
variables described above acts between GVM's that have singular character.
The general theory gives only few tools to deal with this case, but despite this, we prove the
existence of a nonzero homomorphism of GVM's so that the dual operator is the Dirac operator in several variables. 

All GVM's with highest weights on a chosen affine orbit 
and homomorphisms between them form the {\it BGG graph}.
In the thesis, the structure of this graph 
is described (for the specific singular orbit that contains the homomorphism dual to the Dirac operator)
in case the dimension of the variables is
odd. A conjecture is formulated that, if the dimension $n$ of the variables is even and $n/2 \geq k$, where $k$ 
is the number of Clifford variables,
the {\it BGG} graph has the same form as in the odd case. This implies that the Dirac operator in $k$ variables
is only the first operator in a sequence of $G$-invariant operators. For $k=2$ (Dirac operator
in $2$ variables), there are only $2$ further operators in this sequence and their form was found explicitly.
It turned out that these operators
were identical with the operators forming a resolution of the Dirac operator in $2$ variables computed
in \cite{bluebook}. So, these operators not only form a resolution but each of them
is also $G$-invariant. After the identifications between sections of spinor bundles and functions on the flat 
space, action of $G$ includes translation ($f(x,y)\mapsto f(x+u, y+v)$), $\Sl(k)\times \Spin(n)$-action, and much more. 
The invariance of these operators is one of the main results of this work.

Now we give an overview of the chapters and their content.

The first chapter is the preface.

In the second chapter  we give an overview of the ingredients used later.
This includes the highest weight theory, parabolic subalgebras and associated grading, definition
of GVM and the duality between homomorphisms of GVM's
and invariant differential operators. The only new result is Theorem \ref{grading2degree} that
gives a simple  tool how to determine the order of an operator dual to a homomorphism of GVM for
first and second order operators.
Finally, the classical Dirac operator on the flat space is introduced.

The third chapter is devoted purely to properties of GVM's. Results of
Bernstein-Gelfand-Gelfand and Lepowski are presented, partially with comprehensive proofs
(lemma \ref{kernel} and Theorem \ref{zeromap})
and  with a slight extension of
Lepowski's theorem for non-integral weights (Theorem \ref{genlepowski}). Further, we show that
this theorem cannot be generalized to weights of singular characters.
We give a precise definition of {\it singular Hasse graph}
and {\it BGG graph} and give a conjecture that, in most cases, these graphs coincide.

In the fourth chapter, we choose the pair $(G=\Spin(n+1,1),P)$ and particular weights $\mu, \lambda$
and show the existence of a nonzero homomorphism between GVM's induced by $P$-representations
with these highest weights. We show that this is dual to a differential operator
that is, locally, the usual Dirac operator. Further, we choose
the pair $(G=\Spin(n+k,k), P)$ and particular weights $\mu, \lambda$ and we again show the
existence of the homomorphism of the GVM's. We claim that this is dual to the
Dirac operator in $k$ variables.
To give a meaning to this, we need to assign to each section of the spinor bundle
a function $f$ on $\lieg_-$ and then to restrict the operator to such sections, so that
the corresponding functions do not depend on $\lieg_{-2}$ and can be considered as
functions of $\lieg_{-1}$ only.

In the fifth chapter, we compute the affine orbit of the weights $\mu, \lambda$ introduced
in chapter four. We show that in case $n$ is even, the structure of the singular Hasse graph
does not
depend on $n$ for $n/2\geq k$ but becomes larger for $n/2<k$ (its form is described in the case
$k=n/2+1$). This does not happen in odd dimension. 
If $n$ is odd, the BGG graph was computed and its structure is
independent of $n$. We show that all homomorphisms of GVM's in the odd case
are standard whereas in the even case the second order homomorphisms are nonstandard,
if they exist. The existence of these second order homomorphisms, however, is not proved in general
and it remains as a conjecture.

The sixth chapter is devoted to a study of the sequence of GVM's that corresponds to $k=2$
(Dirac operator in $2$ variables). In this case, the sequence consists of only $3$ operators,
whereas the second one is nonstandard and of second order. The existence of the
nonstandard operator
is proved by computing the {\it extremal vector}. This is a vector in the GVM
that is the image of the highest weight vector in another GVM by the homomorphism in question.
This vector determines the homomorphism uniquely and the only condition on it
is its weight and the fact that it is annihilated by the action of all positive root
spaces in $\lieg$. Finally, we used the explicit form of this extremal vector and
``translated'' the GVM homomorphism into operators. We found the form of the three
operators to be (assigning functions on $\lieg_{-}$ to sections and restricting to
functions constant in $\lieg_{-2}$)
\begin{eqnarray*}
f\mapsto &&\mmmm {D_1 f\choose D_2 f}\\
{g_1\choose g_2}\mapsto &&\mmmm {D_1 D_1 g_2 - D_2 D_1 g_1\choose D_1 D_2 g_2-D_2 D_2 g_1}\\
{h_1\choose h_2}\mapsto &&\mmmm D_1 h_2 - D_2 h_1
\end{eqnarray*}
where $D_i=\sum_j e_j \cdot\partial_{ij}$ for $i=1,2$.
This coincides with the resolution computed in \cite[pp. 238]{bluebook}. The equation given by the third
operator $D_1 h_2 - D_2 h_1=i$ has solution for any smooth function $i$ so the resolution ends here.

\newpage
%
\newpage
\section{Introduction}

\subsection{Semisimple Lie algebras}
\label{intro}
Let $\lieg$ be a (real or complex) {\it semisimple Lie algebra}.
We fix a maximal commutative subalgebra $\lieh$ of $\lieg$ called {\it Cartan subalgebra}
and a set of {\it positive roots} $\Phi^+$ for
($\lieg, \lieh$). Let $\Phi=\Phi^+\cup -\Phi^+$ be the set of all {\it roots.} 
The {\it root space} $\lieg_\phi$ corresponding to the root $\phi$ is one-dimensional and 
consists of elements $e_\phi$ such that $[h,e_\phi]=\phi(h)e_\phi$ for each $h\in\lieh$.

It is well known that, for $\lieg$ semisimple, $\lieg=\lieh\oplus_{\phi\in\Phi} \lieg_\phi$.
For a fixed $(\lieg,\lieh,\Phi^+)$, we define $\Delta=\{\alpha_1,\ldots,\alpha_n\}$ 
to be the set of {\it simple roots} (basis of $\lieh^*$ so that each positive root is 
an integral combination of $\alpha_i$'s with nonnegative coefficients).

The {\it Killing form} $(a,b)\mapsto \Tr(\ad(a)\ad(b))$ defines a duality between $\lieh$ and $\lieh^*$.
For each root $\phi$ we define the {\it coroot} $H_\phi:=\frac{2}{(\phi,\phi)}\phi\in\lieh$ where
we identified $\phi$ with an element of $\lieh$ via the Killing form.

We will call elements of $\lieh^*$ {\it weights}. 
The {\it fundamental weights} $\varpi_1,\ldots,\varpi_n$ are elements of $\lieh^*$ dual
to the simple coroots $H_{\alpha_1},\ldots,H_{\alpha_n}$. 
Fundamental weights form a basis of $\lieh^*$ 
and we will sometimes denote weights by its coefficients in this basis. We define
an ordering on the weights by $\mu\leq\lambda$ if and only if $\lambda-\mu$ is a sum
of positive roots with nonnegative integral coefficients.

A weight $\mu\in\lieh$ is said to be {\it dominant}, if 
$H_\alpha(\mu)\geq 0$ for all $\alpha\in\Delta$ and {\it strictly dominant}, if 
$H_\alpha(\mu)>0$ for all $\alpha\in\Delta$. The set of dominant weights is 
sometimes called {\it fundamental Weyl chamber}.
A weight $\mu\in\lieh$ is called to be {\it integral},
if $H_{\alpha}(\mu)\in\Z$ for all $\alpha\in\Delta$. We denote by $P$ the set of integral weights
(it is called {\it weight lattice}) and by $P^{++}$ the set of dominant integral weights.
We see that
a weight is integral, if it is an integral combination of fundamental weights and dominant,
if it is a nonnegative combination of fundamental weights. 

A {\it representation} of $\lieg$ is a vector space $\V$ that is a  $\lieg$-module, i.e.
there is a homomorphisms of Lie algebras $\varphi: \lieg\to \gl(\V)$. For a representation
$\V$ and $\mu\in\lieh$, we define the {\it weight space} $\V^\mu$ of weight $\mu$ by 
$$\V^\mu:=\{v\in\V, \,\,\forall h\in\lieh \,\,\, \varphi(h)(v)=\mu(h)v\ \}.$$ 
A representation $\V$ is called a {\it highest weight module}, if
it is generated by a weight vector $v_\mu$ that is annihilated by all the positive root spaces in $\lieg$.

The following statements can be found in, e.g. \cite{GW, FH, Humphries}:

For $\lieg$ semisimple, each finite-dimensional representation $\V$ 
of $\lieg$ is a direct sum of irreducible representations. $\V$ is a direct sum of its weight spaces.
There is a one to one correspondence between isomorphism classes of 
finite dimensional irreducible representations and $P^{++}$. This correspondence assigns to
each $\mu\in P^{++}$ an irreducible finite dimensional highest weight module $\V_\mu$ 
with highest weight $\mu$ (and all such modules are isomorphic). 

For $(\lieg, \lieh, \Phi^+)$ we define the {\it lowest form} $\delta\in\lieh^*$ to be 
$\delta:=\half\sum_{\phi\in\Phi^+} \phi$. It is easy to show that $\delta=\sum_j \varpi_j$.
For each root $\phi$ we define
the {\it root reflection} $s_\phi:\lieh^*\to\lieh^*$, $\mu\mapsto \mu-H_{\phi}(\mu)\phi$ 
(it is a reflection in $\lieh^*$ with respect to the hyperplane orthogonal to $\phi$, and the Killing metric).
The {\it Weyl group} $W$ is the (finite) group generated by all simple root reflections $s_\alpha$. For any
element $w\in W$, we define the {\it length} $l(w)$ to be the minimal number $k$ so that 
$w=s_{\alpha_1}\ldots s_{\alpha_k}$ for some $\alpha_j\in\Delta$.

We will be interested mainly in the case when $\lieg$ is the complex orthogonal Lie algebra $\so(l,\beta)$
consisting of endomorphisms $A$ of $\C^l$ such that $\beta(x,Ay)+\beta(Ax,y)=0$ for $x,y\in\C^l$, $\beta$
being some non-degenerate symmetric bilinear form on $\C^l$. This is, clearly, the Lie algebra of the orthogonal
Lie group $\so(l,\beta)$. All such algebras, for different choices of $\beta$, are isomorphic, but for
convenient computing, we will choose the form $\beta(x,y)=\sum_{j=1}^l x_j y_{l+1-j}$. In this case, we
can describe elements of $\so(l,\beta)$ explicitly.

Let us suppose that $l=2n$ is even. The Lie algebra has rank $n$ and is sometimes denoted by $D_n$.
In the formalisms of \cite{GW}, it consists of matrices of size $2n\times 2n$, 
antisymmetric with respect to the anti-diagonal and 
the Cartan subalgebra $\lieh$ can be chosen to be the subalgebra of diagonal matrices. It has dimension $n$.
In this formalism, we introduce the orthogonal (with respect to the Killing form) basis $\{\epsilon_1,\ldots,\epsilon_n\}$ 
of $\lieh^*$ defined
by $\epsilon_i(\diag(a_1,\ldots,a_n,-a_n,\ldots,-a_1)):=a_i$. The positive roots can be chosen to be 
$\epsilon_i\pm \epsilon_j,\,\,i<j$.
For $n=3$, the matrices have the following form:
\begin{equation*}
\left(
\begin{tabular}{ccc|ccc}
$h_{1}$ & $x_{12}$ & $x_{13}$ & $x_{14}$ & $x_{15}$ & 0 \\
$y_{21}$ & $h_{2}$ & $x_{23}$ & $x_{24}$ & 0 & $-x_{15}$\\ 
$y_{31}$ & $y_{32}$ & $h_{3}$ & $0$ & $-x_{24}$ & $-x_{34}$\\ 
\hline
$y_{41}$ & $y_{42}$ & $0$ & $-h_{3}$ & $-x_{23}$ & $-x_{13}$\\ 
$y_{51}$ & $0$ & $-y_{42}$ & $-y_{32}$ & $-h_{2}$ & $-x_{12}$\\ 
$0$ & $-y_{51}$ & $-y_{41}$ & $-y_{31}$ & $-y_{21}$ & $-h_1$
\end{tabular}\right)
\label{matrix_even}
\end{equation*}
The matrices with $x_{ij}=1$ for some $i,j$ and having $0$ on other positions are generators of
positive root spaces. Similarly, matrices with $y_{ij}=1$ and $0$ elsewhere are generators of negative root spaces
and the matrices with $h_j=1$ generate the Cartan subalgebra $\lieh$. 
Further, we will usually denote by $x_{ij}$ resp. $y_{ij}$ the corresponding root space generator and by $h_j$
the elements of $\lieh$.
In the upper left square, the generator of a positive root space $x_{ij}$ corresponds to 
the root $\epsilon_i-\epsilon_j$ and in the upper right square,  $x_{ij}$ is a generator of the root space
corresponding to the root $\epsilon_{i}+\epsilon_{2n+1-j}$ (in our case, $n=3$).
The generator of negative root space for the root $-\phi$ is just the transposed generator of the $\phi$-root space.
In the $\epsilon_j$-basis, fundamental weights are $\varpi_j=[1,\ldots,1,0,\ldots,0]$ where the last
$1$ is on the $j$-th position for $j\leq n-2$ and $\varpi_{n-1}=\half[1,\ldots,1]$, $\varpi_n=\half[1,\ldots,1,-1]$.
A weight $\mu=[a_1,\ldots,a_n]$ is dominant exactly if
$a_1\geq a_2\geq\ldots\geq a_{n-1}\geq |a_n|$ and strictly dominant if strict inequalities hold. 
The weight $\mu$ is integral, if the $a_i$'s are all integers or all half-integers. The lowest form is
$\delta=[n-1,n-2,\ldots,1,0]$.

In case $l=2n+1$, the algebra is denoted by $B_n$ and consists of matrices $(2n+1)\times (2n+1)$ 
antisymmetric with respect to the anti-diagonal.
The Cartan subalgebra $\lieh$ can be again chosen to be the subalgebra of diagonal matrices. It has dimension $n$.
We define the orthogonal basis $\{\epsilon_1,\ldots,\epsilon_n\}$ of $\lieh^*$ by
$\epsilon_i(\diag(a_1,\ldots,a_n,0,-a_n,\ldots,-a_1)):=a_i$. The positive roots are then 
$\epsilon_i\pm \epsilon_j,\,\,i<j$ and $\epsilon_j$, $1\leq j\leq n$.
For $n=3$, the matrices have the following form:
\begin{equation*}
\left(
\begin{tabular}{ccc|c|ccc}
$h_{1}$ & $x_{12}$ & $x_{13}$ & $x_{14}$ & $x_{15}$ & $x_{16}$ & 0 \\
$y_{21}$ & $h_{2}$ & $x_{23}$ & $x_{24}$ & $x_{25}$ & 0 & $-x_{15}$\\ 
$y_{31}$ & $y_{32}$ & $h_{3}$ & $x_{34}$ & $0$ & $-x_{24}$ & $-x_{34}$\\ 
\hline
$y_{41}$ & $y_{42}$ & $y_{43}$ & $0$ & $-x_{34}$ & $-x_{24}$ & $-x_{14}$\\
\hline
$y_{51}$ & $y_{52}$ & $0$ & $-y_{43}$ & $-h_{3}$ & $-x_{23}$ & $-x_{13}$\\ 
$y_{61}$ & $0$ & $-y_{52}$ & $-y_{42}$ & $-y_{32}$ & $-h_{2}$ & $-x_{12}$\\ 
$0$ & $-y_{61}$ & $-y_{51}$ & $-y_{41}$ & $-y_{31}$ & $-y_{21}$ & $-h_1$
\end{tabular}\right)
\label{matrix_odd}
\end{equation*}
Similarly as before, the matrices with $x_{ij}=1$ for some $i,j$ and $0$ on other positions are generators of
positive root spaces and matrices with $y_{ij}=1$ for some $i,j$ and $0$ elsewhere are generators of negative root spaces.
In the upper left square, $x_{ij}$ corresponds to roots $\epsilon_i-\epsilon_j$. In the upper right square,  
$x_{ij}$ is a generator of the root spaces corresponding to $\epsilon_{i}+\epsilon_{2n+2-j}$ and in the middle-column,
$x_{i,n+1}$ generate the root space for $\epsilon_i$.
Generators of negative root spaces are again just transposed generators of positive root spaces.
In the $\epsilon_j$-basis, fundamental weights are $\varpi_j=[1,\ldots,1,0,\ldots,0]$ where the last
$1$ is on the $j$-th position for $j\leq n-1$ and $\varpi_{n}=\half[1,\ldots,1]$.
A weight $\mu=[a_1,\ldots,a_n]$ is dominant exactly if
$a_1\geq a_2\geq\ldots\geq a_{n-1}\geq a_n\geq 0$ and strictly dominant if strict inequalities holds. 
The weight $\mu$ is integral, if the $a_i$'s are all integers or all half-integers. The lowest form is 
$\delta=\half[2n-1,\ldots,3,1]$.

\subsection{Parabolic subalgebras and grading}
\label{parapara}
For the triple $(\lieg,\lieh,\Phi^+)$ where $\lieg$ is semisimple, $\lieh$ a Cartan subalgebra and $\Phi^+$
the set of positive roots, we define the Lie algebra $\lien:=\oplus_{\phi\in\Phi^+} \lieg_\phi$.
In the orthogonal case, choosing $\lieh$ and $\Phi^+$ as above, $\lien$ consists of strictly
upper triangular matrices in $\lieg$.  Further, we define the {\it Borel subalgebra} 
$\lieb:=\lieh\oplus\lien$.

We call any subalgebra $\liep$ of $\lieg$ a {\it parabolic subalgebra}, if it contains $\lieb$ 
(associated to some choice of ($\lieh,\Phi^+$). In that case, $\lieh$ is also a Cartan subalgebra
of $\liep$ and $\Phi^+$ is a set of positive roots for $(\liep,\lieh)$ as well (however, $\liep$
usually does not contain all negative root spaces of $\lieg$). It was shown in \cite{Cap} that
there is a $1-1$ correspondence between parabolic subalgebras of $\lieg$ (by fixed $\lieh,\Phi^+$)
and subsets of $\Delta$: to any $\Sigma\subset\Delta$ we assign the parabolic Lie algebra
$
\liep_\Sigma:=\sum_{\phi\in A} \lieg_{-\phi}\oplus\lieb
$
where $A\subset\Phi^+$ consists of those positive roots that can be expressed as a sum of simple 
roots that are not in $\Sigma$.
Each parabolic subalgebra $\liep$ is of this type. Further, there is a $1-1$ correspondence between
parabolic subalgebras of $\lieg$ and gradings $\lieg=\oplus_{i=-k}^k \lieg_i$ of $\lieg$. 
Given $\Sigma\subset\Delta$, 
the set $\lieg_i$ ($i\neq 0)$ is defined to be 
$\oplus_{\phi\in A_i} \lieg_\phi$, where $A_i$ contains elements 
$\phi=\sum_{\alpha_j\in\Delta} c_j \alpha_j$ 
such that $\sum_{\alpha_j\in\Sigma} c_j=i$, and $\lieg_0=\lieh\oplus_{\phi\in A_0} \lieg_\phi$. 
Given a grading $\oplus_j \lieg_j$, the parabolic subalgebra is then $\liep=\oplus_{j\geq 0}\lieg_j$.

The set of simple roots $\Delta$ can be described by the {\it Dynkin diagram}: it is a diagram
where nodes are simple roots and edges indicate angles between them, see, e.g. \cite{Humphries}.

A parabolic subalgebra can be given by crossing the nodes representing simple roots in $\Sigma$ in the Dynkin diagram. 
For example, the Dynkin diagram for $\so(2(k+n),\C)$ with $\alpha_k$ crossed (i.e. $\Sigma=\{\alpha_k\}$)
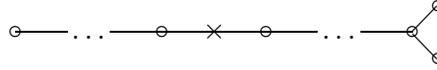
$$
\dynkin\whiteroot{}\link\ldots\link\whiteroot{}\link\noroot{}\link\whiteroot{}\link\ldots\link\whiteroot{}\whiterootdownright{}\whiterootupright{} \enddynkin
$$
induces a gradation 
\begin{equation}
\label{2grading}
\lieg=
\left(
\begin{tabular}{c|c|c}
$\lieg_0$ & $\lieg_{1}$ & $\lieg_{2}$ \\
\hline
$\lieg_{-1}$ & $\lieg_0$ & $\lieg_{1}$ \\
\hline
$\lieg_{-2}$ & $\lieg_{-1}$ & $\lieg_0$
\end{tabular}
\right)
\end{equation}
where $\lieg_0=\gl(k,\C)\oplus\so(2n,\C)$, $\lieg_1$ consists of a block of size $k\times n$ 
and its negative transpose with respect to the anti-diagonal, and $\lieg_2$ is a block of size $k\times k$ 
(antisymmetric with respect to the anti-diagonal).
The associated parabolic subalgebra is then 
$\liep=\lieg_0\oplus\lieg_1\oplus\lieg_2$. 
If $k=1$, then $\lieg_{2}=0$ and $\lieg$ is only $1$-graded.

Similarly, in the odd orthogonal case, $\Sigma=\{\alpha_1\}$ induces a $1$-grading of $\lieg$ 
and $\Sigma=\{\alpha_k\}$ induces a $2$-grading for $k\geq 2$.

For each grading $\oplus_j \lieg_j$ of $\lieg$, there is a unique element $E\in \lieg_0$ called
{\it grading element} defined by the property $[E,g_j]=j g_j$ for $g_j\in\lieg_j$.

\begin{lemma}
\label{gradingelement}
Let $\lieg=\so(2(n+k),\C)$, $\liep$ be the parabolic subalgebra corresponding to $\Sigma=\{\alpha_k\}$ and $n\geq 1$. 
Then the grading element associated to $(\lieg, \liep)$ is 
\begin{equation*}
E=\diag(1,\ldots,1,0,\ldots,0,-1,\ldots,-1)=
\left(
\begin{tabular}{c|c|c}
$1$ & $0$ & $0$ \\
\hline
$0$ & $0$ & $0$ \\
\hline
$0$ & $0$ & $-1$
\end{tabular}
\right)
\end{equation*}
the diagonal blocks are of size $k\times k, n\times n$ and $k\times k$.

For $\lieg=\so(2(k+n)+1,\C)$, $\Sigma=\{\alpha_k\}$, $n\geq 1$
the grading element is $$E=\diag(1,\ldots,1,0,\ldots,0,-1,\ldots,-1)$$ (there are $k$ one's and $2n+1$ zeroes 
in the expression).

If the parabolic subalgebra is the Borel subalgebra $\lieb$, then the grading element corresponding to $(\lieg,\lieb)$
is $$E_\lieb=\diag(\ldots,5/2,3/2,1/2,-1/2,-3/2,-5/2,\ldots)$$ in the even orthogonal case and 
$$E_\lieb=\diag(\ldots,2,1,0,-1,-2,\ldots)$$
in the odd orthogonal case.
\end{lemma}

\begin{proof}
In the even $2$-graded case, consider the structure of the grading given by (\ref{2grading}) and we can verify the condition
$[E,g_j]=j g_j$ by commuting the (antisymmetric with respect to the anti-diagonal) matrices. The odd case is similar.
In the even Borel case, the root space $E_{i,j}-E_{2(n+k)+1-j,2(n+k)+1-i}\in\lieg_{j-i}$  for $i+j\leq 2(n+k), i\neq j$ and
the equality  $[E,g_j]=j g_j$ can be again easily verified by commuting ($E_{ij}$ is a matrix with a $1$ on position $(i,j)$
and zeros on other positions). 
\end{proof}

It can be shown that $\lieg_0=\lieg_0^{ss}\oplus \liez$ where $\lieg_0^{ss}$ is semisimple and $\liez$
is the center of $\lieg_0$. Dimension of $\liez$ is equal to the cardinality of $\Sigma$, so in case
$\Sigma=\{\alpha_k\}$, $\liez$ is generated by the grading element. Clearly, $\lieh$ is a Cartan subalgebra for
$\lieg_0$, $\Delta-\Sigma$ is a set of simple roots and $\{\varpi_j;\,\,\alpha_j\notin\Sigma\}$ is a set
of fundamental weights for it.  Any irreducible representation $\V$ of $\lieg_0^{ss}$ can be extended to a 
representation of $\lieg_0$ letting $\liez$ act by $z\cdot v=\nu(z)v$
where $\nu\in\liez^*$ is arbitrary. Let $\mu$ be any weight of $\V$. The number $\mu(E)$ is called
{\it generalized conformal weight} (because $\V$ is generated, as a $\lieg_0$-module, by a highest weight vector and $E$ is in 
the center of $\lieg_0$, this is independent of the choice of $\mu$). 
Further, any irreducible $\lieg_0$-module $\V$ can be extended to an irreducible representation of the whole $\liep$, letting 
$\liep^+=\sum_{j>0} \lieg_j$ act trivially. On the other hand,
if $\V$ is an irreducible representation of $\liep$, the action of $\liep^+$ must
be trivial on it (it follows from Engel's theorem about nilpotent Lie algebras), so,
finite-dimensional irreducible representations of $\liep$ are completely 
described by the highest weight of $\V$ as a $\lieg_0$-module.

Let us denote by $P_\liep$ the set of weights $\mu$ such that $H_\alpha(\mu)\in\Z$ for
all $\alpha\in\Delta-\Sigma$ and call it $\liep$-integral weights (or $\lieg_0^{ss}$-integral).
We say that a weight $\mu$ is $\liep$-dominant (or $\lieg_0^{ss}$-dominant) if 
$H_\alpha(\mu)\geq 0$ for all $\alpha\in\Delta-\Sigma$. 
Similarly, we define a strictly $\liep$-dominant weight.
We denote by $P_\liep^{++}$ the set of $\liep$-integral and $\liep$-dominant weights.
Note that $P_\liep$ is not a lattice, but it consists of 
$\dim(\liez)$-dimensional planes in $\lieh^*$. A weight is in $P_\liep^{++}$ exactly if
it is expressed as $\sum_j c_j \varpi_j$ so that $c_j$ is a nonnegative integer
for each $j$ such that $\alpha_j\notin\Sigma$. We see that there is a $1-1$ correspondence
between $P_\liep^{++}$ and the set of isomorphism classes of irreducible finite dimensional $\liep$-modules. Clearly, 
$P^{++}\subset P_\liep^{++}$.

In case $\so(l,\C)$, we usually express the weights in the $\epsilon_j$-basis defined in \ref{intro}. 
Assume that $\lieg=\so(2(k+n),\C)$, $\Sigma=\{\alpha_k\}$ and $n\geq 1$. 
Then a weight $\mu=[a_1,\ldots,a_k|b_1,\ldots,b_n]$ is $\liep$-dominant, if $a_1\geq a_2\geq\ldots\geq a_k$
and $b_1\geq b_2\geq\ldots\geq b_{n-1}\geq |b_n|$ and $\liep$-integral if $a_i-a_j$ are all integers and 
$b_j$ are either all integers, or all half-integers.
Similarly, in the odd case $\lieg=\so(2(k+n)+1,\C)$, $\Sigma=\{\alpha_k\}$ and $n\geq 1$, a weight 
$\mu=[a_1,\ldots,a_k|b_1,\ldots,b_n]$ is $\liep$-dominant if $a_1\geq\ldots\geq a_k$ and $b_1\geq\ldots\geq b_n\geq 0$
and $\liep$-integral, if $a_i-a_j\in\Z$ and the $b_i$'s are all integers or all half-integers.

From lemma \ref{gradingelement} it follows that the grading element evaluation is
$$[a_1,\ldots,a_k|b_1,\ldots,b_n](E)=a_1+\ldots+a_k.$$
in even as well as in the odd case.

\subsection{True and generalized Verma modules}
\label{vermamodules-general}
The representations of any Lie algebra $\lieg$ are in a natural correspondence
with representations of the associative algebra
$$\univ(\lieg):=T(\lieg)/I$$
where $I$ is the ideal in the tensor algebra generated by all $[x,y]-x\otimes y-y\otimes x$.
It is called {\it universal enveloping algebra}. The filtration 
$T_k(\lieg)=\oplus_{i=0}^k T^i(\lieg)$ projects to a filtration $\univ_k(\lieg)$
of $\univ(\lieg)$.

We see from the definition that $\univ(\lieg)$ is a $\lieg$-module, if the
action of $\lieg$ on $\univ(\lieg)$ is the left multiplication, i.e.
$g\cdot (g_1\otimes\ldots\otimes g_k \mathrm{mod}\, I)=
g\otimes g_1\otimes\ldots\otimes g_k \mathrm{mod}\, I$
and, extending this action, a left $\univ(\lieg)$-module. 
Considering multiplication from the right, $\univ(\lieg)$
is a right $\univ(\lieg)$-module and for any subalgebra $\lieg_1$
of $\lieg$, $\univ(\lieg)$ is a right $\univ(\lieg_1)$-module. The right and left actions
obviously commute.

The PBW-theorem (see, e.g. \cite{Humphries}) states that for any basis 
$v_1,\ldots,v_n$ of $\lieg$, the set 
$$\{v_{i_1}\otimes v_{i_2} \otimes\ldots\otimes v_{i_j} \,
\mathrm{mod}\, I;\,\, i_1\leq i_2\leq\ldots\leq i_j, j\in\N\}$$
is a basis of $\univ(\lieg)$ (as a vector space).

Assume that $\lieg$ is semisimple and $(\lieh,\Phi^+)$ is chosen.
For any $\mu\in\lieh^*$, $\C$ can be given the structure of a $\lieb$-module, 
where the action of $\lieh$ is 
$h\cdot v=\mu(h) v$ for $h\in\lieh, v\in\C$ and the nilpotent algebra
$\lien$ acts trivially. We will
denote this representation by $\C_\mu$. We see that $\C_\mu$ is 
a left $\univ(\lieb)$-module. Further, we know that 
$\univ(\lieg)$ is a right $\univ(\lieb)$ module
and left $\univ(\lieg)$-module and that both actions commute. So, we
can define the {\it Verma module} to be the $\lieg$-module
$$M(\mu):=\univ(\lieg)\otimes_{\univ(\lieb)}\C_\mu$$
If we choose $\{y_1,\ldots,y_k,h_1,\ldots,h_n,x_1,\ldots,x_k\}$ to be 
an ordered basis of $\lieg$ so that $y_j$ are generators of the
negative root spaces in $\lieg$, $h_j$ generates $\lieh$ and $x_j$ are generators
of the positive root spaces in $\lieg$, it follows easily from the PBW theorem
that, as a vector space (and as $\lieg_-$-module), $M(\mu)\simeq\univ(\lieg_-)\otimes \C_\mu$
where $\lieg_-:=\sum_{\phi\in\Phi^+} \lieg_{-\phi}$. The tensor product is
over the $\C$, resp. $\R$, if $\lieg$ is a complex, resp. real, Lie algebra.
$M(\mu)$ is a highest weight module with highest weight $\mu$ and
$\mu$-weight vector $1\otimes 1$. Simple commutation relations imply that 
the weight of $y_1\ldots y_k\otimes 1$ is $\mu-\sum_j \root(y_j)$, where
$\root(y_j)$ is the root to which $y_j$ is the root space generator.

The Verma modules are maximal in the sense that each highest weight module
with highest weight $\mu$ is isomorphic to a factor of the Verma module $M(\mu)$, see
e.g. \cite{Humphries, Dixmier}.

Let $\liep$ be a parabolic subalgebra of $\lieg$ and $\mu\in P_\liep^{++}$.
Then there is an irreducible representation $\V_\mu$ of $\liep$ with highest 
weight $\mu$ and we define the {\it generalized Verma module} 
$$
M_\liep(\mu):=\univ(\lieg)\otimes_{\univ(\liep)} \V_\mu\simeq \univ(\lieg_-)\otimes \V_\mu,
$$
where $\lieg_-=\lieg_{-k}\oplus \ldots \oplus \lieg_{-1}$ is the negative part
of $\lieg$ in the grading corresponding to $\liep$, see \ref{parapara}. The second tensor
product, again, is over $\C$ resp. $\R$ for a complex resp. real Lie algebra $\lieg$.
Sometimes we will write $M_\liep(\V)$ for the generalized Verma module induced
by a $\liep$-module $\V$.

$M_\liep(\mu)$ is naturally a $\lieg$-module by the left action of $\lieg$ on $\univ(\lieg)$. We call
$M_\liep(\mu)\to M_\liep(\lambda)$ a {\it homomorphism of generalized Verma modules}, if 
it is a $\lieg$-homomorphism.
If $\V_\mu$ is a $P$-module, where $P$ is a Lie group with Lie 
algebra $\liep$, $M_\liep(\mu)$ is also a 
$(\lieg, P)$-module, the action of $p\in P$ being 
$$y_1\ldots y_k\otimes v\mapsto \Ad(p)(y_1)\ldots \Ad(p)(y_k)\otimes (p\cdot v).$$
$M_\liep(\mu)$ is a highest weight module with highest weight $\mu$ and 
$\mu$-weight 
vector $1\otimes v_\mu$, where $v_\mu$ is the highest weight vector in 
$\V_\mu$. Let us denote the generators of negative root spaces in $\lieg_-$
by $y_j$ and the generators of negative root spaces in $\lieg_0$ by $Y_k$.
It follows from the PBW theorem and the fact that $\V_\mu$ is generated by $v_\mu$
that the vectors
\begin{equation}
\label{yYv}
y_{i_1}\ldots y_{i_l}\otimes Y_{j_1}\ldots Y_{j_s} v_\mu
\end{equation}
generates $M_\liep(\mu)$ and, if $Y_{j_1}\ldots Y_{j_s} v_\mu$ 
is nonzero, then $(\ref{yYv})$ is a weight vector of the weight 
$\mu-\sum_u \root(Y_{j_u})-\sum_v \root(y_{i_v})$. However,
$Y_{j_1}\ldots Y_{j_s}v_\mu$ may be zero in $\V_\mu$, so the generalized 
Verma modules are in general ``smaller'' then (true) Verma modules.

\subsection{Invariant differential operators}
\label{invdifop}

Let $G$ be a semisimple Lie group and $H$ a 
(closed) Lie subgroup. The {\it homogeneous space} $G/H$ is a smooth manifold
and the projection $\pi:G\to G/H$ induces a {\it principal fiber bundle}
on $G/H$ with structure group $H$, see, e.g. \cite{Sharpe, Cap, Slovak1}.
To each $H$-module $\V$ we can associate a vector bundle 
$V(G/H):=G\times_H \V$.
The equivalence class of $[g,v]\in G\times \V$ in $G\times_H \V$ will be denoted by
$[g,v]_H$ and the projection $G\times_H \V\to G/H$ is given by $[g,v]_H\mapsto gH$. 
There is an isomorphism between sections $\Gamma(V(G/H))$ and
$H$-equivariant functions $C^\infty(G,\V)^H$, i.e. functions
$f: G\to \V$ such that $f(h^{-1}v)=h\cdot f(v)$ for each $h\in H$. This isomorphism
assigns to a section $s\in\Gamma(V(G/H))$ the function $f$ so that $f(g)$
is the unique $v\in\V$ for that $s(gH)=[g,v]_H$. On the other hand, 
given an equivariant function $f$, we can define the section $s$ by
$s(gH)=[g,f(g)]_H$ and this definition is independent of the choice of $g$.

The group $G$ has a natural left action on the following spaces:
\begin{itemize}
\item{on $G/H$: $g_1\cdot gH=g_1 gH$}
\item{on $V(G/H)$: $g_1\cdot [g,v]_H=[g_1g,v]_H$}
\item{on $\Gamma(V(G/H))$: $(g_1\cdot s)(gH)=g_1\cdot (s(g_1^{-1}gH))$}
\item{on $C^\infty(G,\V)^H$: $(g_1\cdot f)(g)=f(g_1^{-1}g)$}
\end{itemize}
It is easy to verify that the last two actions are compatible with
the isomorphism $\Gamma(V(G/H))\simeq C^\infty(G,\V)^H$ described 
above.

Let $\V$ and $\W$ be representations of $H$ and $V(G/H)$ and $W(G/H)$ the 
associated vector bundles. A map $D:\Gamma(V(G/H))\to\Gamma(W(G/H))$ is called 
an {\it operator}. The operator $D$ is called {\it differential of order $k$},
if the value $Ds(x)$ depends only on $s(x)$ and derivations of $s$ in $x$
up to order $k$ (this makes sense, because in a neighborhood $U$ of $x$,
the bundle is trivialized to $U\times \V$ by some bundle map and the 
transitions between bundle maps are $C^\infty$). An operator $D$  is 
called {\it invariant} if it commutes with the action of $G$ on sections.
We will deal only with linear invariant differential operators in this work.

It is easy to see that an invariant operator $D$ 
is completely determined by the values $Ds(0)$ on sections, 
where $0=eH\in G/H$ is the image of the identity element of $G$ in $G/H$. 
If the operator is of order $k$, the values $Ds(0)$ depend only on
the $k$-jet $J^k_0 s$ of sections in $0$ and the operator is determined by a 
$H$-homomorphism $\varphi_D: J^k_0(V(G/H))\to W$, where the action of $H$ 
on $J^k_0(V(G/H))$ is the action on representatives. The operator
acts by $Ds(0)=\varphi_D J^k_0 (s)$ and $\varphi_D$ must be a 
$H$-homomorphism because $D$ is $G$-invariant and the action of
$H$ takes $k$-jets of sections in $0$ to themselves.

The most trivial example is the Euclidean space $\R^n$
considered as a homogeneous space $\Euc(n)/\SO(n)$. Let
$\V=\R$ be the trivial representation of $\SO(n)$, then the
sections $\Gamma(V(G/H))$ are just functions on $\R^n$. Its $1$-jets
in $0$ can be identified with 
$\R\oplus (\R^n)^*\otimes \R\simeq \R\oplus (\R^n)^*$ 
by $f\mapsto (f(0),df(0))$.
The action of $H$ on $1$-jets in $0$ implies that $(\R^n)^*$ is the
dual defining representation of $\SO(n)$ and we see that the
only (nontrivial) $H$-homomorphisms from $J^1_0$ to $\W$ exists
for $\W=(\R^n)^*$ and $\varphi: \R\oplus(\R^n)^*\to (\R^n)^*$ the projection
onto the second component. The vector bundle associated to $(\R^n)^*$
is the cotangent bundle and $D$ is the differential.

Let $\lieg$ be the (semisimple) Lie algebra of $G$.
We will call $P$ a {\it parabolic subgroup} of $G$ if it is a Lie subgroup
and its Lie algebra $\liep$ is a parabolic subalgebra of $\lieg$. 
For a parabolic subalgebra $\liep$ of $\lieg$, $P$ can be defined 
explicitly as the set $\{p\in G, \Ad(p)(\liep)\subset\liep\}$.
We know that each irreducible finite dimensional representation
of $P$ is a representation of $\liep$ and can be characterized by
its highest weight $\mu\in P_\liep^{++}$. 
We shall assume that $G$ and $P$ are 
real Lie groups and consider smooth sections,
not holomorphic. So, the algebras ($\lieg, \liep$) are real as well.

The space $J_{0}^k(V(G/P))$ of sections can be identified with
$J_e^k(C^\infty(G,\V))^P$, where $e\in G$ is the identity element. 
It can be shown that this is dual, as a $P$-module, to 
$\univ_k(\lieg)\otimes_{\univ(\liep)} \V_\lambda^*$
($\univ_k(\lieg)$ is the $k$-th
filtration of $\univ(\lieg)$) and the duality is given by 
\begin{equation}
\label{duality}
\br{Y_1\ldots Y_l\otimes_{\univ(\liep)} A, j_e^k f}=A((L_{Y_1}\ldots L_{Y_l} f)(e))
\end{equation}
for $l\leq k$, $A\in \V_{\lambda}^*$, $j_e^k f$ the $k$-jet of $f$ in $e$, 
$Y_j\in\lieg$ and $L_{Y_j}$ the derivation with respect to the
left invariant vector fields induced by $Y_j$. 

An invariant differential operator of order $k$ is determined by some
$P$-homomorphism $\varphi_D: J_e^k(C^\infty(G,\V_\lambda))^P\to\V_\mu$
and, applying $(\ref{duality})$, the dual map 
$\varphi_D^*: \V_\mu^*\to \univ_k(\lieg)\otimes_{\univ(\liep)} \V_\lambda^*$.
The right side is a $P$-submodule of $M_\liep(\V_\lambda^*)$. 
Further, each $P$-homomorphism $\phi: \V_\mu^*\to M_\liep(\V_\lambda^*)$ can be 
extended to a $(\lieg,P)$-homomorphism $M_\liep(\V_\mu^*)\to M_\liep(\V_\lambda^*)$
of generalized Verma modules by 
$y_1\ldots y_l\otimes_\R v\mapsto y_1\ldots y_l\otimes_\R \phi(v)$ for $y_j\in\lieg_-$, 
the action of $\liep$ on $\V_\lambda^*$ being the infinitesimal action of $P$ 
(we identified $M_\liep(\nu)\simeq \univ(\lieg_-)\otimes_\R \V_\nu)$.

It follows that there is a duality between invariant linear differential operators 
$D: \Gamma(G\times_P \V_\mu)\to \Gamma(G\times_P \V_\lambda)$ of any finite order and $(\lieg,P)$-homomorphisms
of generalized Verma modules
$M_\liep(\V_\mu^*)\to M_\liep(\V_\lambda^*)$ (see \cite{Cap} for details).

Let $\mu,\lambda\in P_\liep^{++}$ and assume that $\V_\mu$ and $\V_\lambda$ are $P$-modules, so that
the action of $\liep$ is the infinitesimal action (for example, if $P$ is simply connected, this is true for all
$\mu, \lambda\in P_\liep^{++}$). 
Then each $\lieg$-homomorphism $M_\liep(\mu)\to M_\liep(\lambda)$ lifts to a $(\lieg,P)$-homomorphism,
defining the action of $P$ in a natural way by 
$p\cdot (y_1\ldots y_l\otimes_\R v):=\Ad(y_1)\ldots \Ad(y_l)\otimes_\R p\cdot v$. 

The following theorem is an important tool for determining the order of an operator, dual to
a homomorphism of generalized Verma modules.

\begin{theorem}
\label{grading2degree}
Let $\mu, \lambda$ be highest weights of some irreducible finite-dimensional 
$P$-modules $\V_\mu, \, \V_\lambda$ and $\phi: M_\liep(\mu)\to M_\liep(\lambda)$ be a nonzero 
homomorphism of generalized Verma modules. Let $E$ be the grading element
for $(\lieg,\liep)$ and assume that $(\lambda-\mu)(E)\in \{1,2\}$.
Then the corresponding (dual) invariant differential operator 
$\Gamma(G\times_P \V_\lambda^*)\to \Gamma(G\times_P \V_\mu^*)$ has degree $(\lambda-\mu)(E)$.
\end{theorem}

\begin{proof}
Let $v_\mu$ be the highest weight vector of $\V_\mu$. 
Then $\phi$ is completely determined by $\phi(1\otimes v_\mu)$,
the image of the highest weight vector in $M_\liep(\mu)$, because generalized
Verma modules are highest weight modules. Define $y_j\in\univ(\lieg_-)$ and
$v_j\in\V_\lambda$ so that $\phi(1\otimes v_\mu)=\sum_j y_j\otimes v_j$ (this is
possible since $M_\liep(\lambda)\simeq \univ(\lieg_-)\otimes \V_\lambda$).

Let $k$ be the maximal integer such that $y_i\in\univ_k (\lieg_-)$ for some $y_i$ and let $0\neq g_0\in\univ(\lieg_0)$.
Then
$\phi$ maps $1\otimes g_0\cdot v_\mu=g_0\otimes_{\univ(\liep)}v_\mu$ to 
$$\sum_j g_0 y_j\otimes_{\univ(\liep)} v_j=\sum_j (y_j g_0 +[g_0,y_j])\otimes_{\univ(\liep)}v_j=
\sum_j (y_j\otimes_\R g_0\cdot v_j + [g_0,y]\otimes_\R v_j)$$
because $[g_0,y_j]\in\univ(\lieg_-)$ for $g_0\in \univ(\lieg_0), y_j\in\univ(\lieg_-)$ ($[a,b]=ab-ba$ is the
commutator in the associative algebra $\univ(\lieg)$). Simple commutation relations show that, if
$y_j\in\univ_{l}(\lieg_-)-\univ_{l-1}(\lieg_-)$, then  $[g_0, y_j]\in\univ_{l}(\lieg_-)-\univ_{l-1}(\lieg_-)$ as well. 
This implies that, if $k$ is the smallest integer such that
$\phi$ maps $1\otimes v_\mu$ into $\univ_{k}(\lieg_-)\otimes \V_\mu$, then $\phi$ maps
$1\otimes g_0\cdot v_\mu$ into $\univ_{k}(\lieg_-)\otimes \V_\mu$ as well (but not to $\univ_{k-1}(\lieg_-)\otimes \V_\lambda$). 
$\V_\mu$ is generated
by $v_\mu$, so  we proved that $\phi$ maps $1\otimes \V_\mu$ into $\univ_{k}(\lieg_-)\otimes \V_\mu$.

In particular, for any $v\in\V_\mu$, 
$\phi(1\otimes_\R v)=\sum_j \tilde{y}_j\otimes_\R \tilde{v}_j$
for some $\tilde{v}_j\in\V_\lambda$, $\tilde{y}_j\in\univ_{k}(\lieg_-)$ and 
$\tilde{y}_i\notin \univ_{k-1}(\lieg_-)$ for some $i$. 
Without loss of generality, we can assume that $\tilde{y}_j=y_1^{(j)}\ldots y_{l(j)}^{(j)}$
for some $y_u^{(j)}\in\lieg_-$, $l(j)\leq k$ and $l(i)=k$.

Applying the duality (\ref{duality}),
the differential operator $D$ satisfies
$$v((Df)(0))=\sum_j \tilde{v}_j(L_{y_1^{(j)}} \ldots L_{y_{l(j)}^{(j)}} (f)(0)),$$ where $L_{y_u^{(j)}}$ are the left
invariant vector fields generated by $y_u^{(j)}\in\lieg_-$. So, the operator $D$
dual to the homomorphism is of order $k$.

Let us suppose that the operator has order $k$, i.e. $\phi$ maps $1\otimes v_\mu$ into 
$\univ_k(\lieg_-)\otimes \V_\lambda$ but not into $\univ_{k-1}(\lieg_-)\otimes \V_\lambda$.
Let $\{y_1,\ldots y_n\}$ be an ordered basis of $\lieg_-$ that consists of generators of 
negative root spaces in $\lieg_-$.

Let $\phi(1\otimes v_\mu)=\sum_j \tilde{y}_j\otimes v_j$ and assume that all the $v_j$'s are 
weight vectors in $\V_\lambda$ and $\tilde{y}_j$ is a product of the $y_j$'s
(it follows from the PBW theorem that such expression is always possible). Then all $\tilde{y}_j\otimes v_j$ are weight vectors
and, because their sum is a weight vector of weight $\mu$, each $\tilde{y}_j\otimes v_j$ is
a weight vector of weight $\mu$ as well.

Because $\phi(1\otimes v_\mu)\notin \univ_{k-1}(\lieg_-)\otimes \V_\lambda$, there exists $i$ such that 
$\tilde{y}_i=y_{i_1}\ldots y_{i_k}$ is a product of $k$ elements. Let $u_j\in\N$ be defined by 
$y_{i_j}\in\lieg_{-u_j}$.
The action of the grading element on $y_{i_1}\ldots y_{i_k}\otimes v_i$ is
\begin{eqnarray*}
&&\mmmm E\cdot (y_{i_1}\ldots y_{i_k}\otimes_\R v_i)=
Ey_{i_1}\ldots y_{i_k}\otimes_{\univ(\liep)}v_i=\\
&&\mmmm=(y_{i_1}E+[E,y_{i_1}])y_{i_2}\ldots y_{i_k}\otimes_{\univ(\liep)}v_i=\ldots=\\
&&\mmmm=y_{i_1}\ldots y_{i_k} (\lambda(E)-u_1-\ldots -u_k)\otimes_\R v_i
\end{eqnarray*}
But $y_{i_1}\ldots y_{i_k}\otimes v_i$ is a weight vector of weight $\mu$, so the 
left hand side equals $\mu(E)(y_{i_1}\ldots y_{i_k}\otimes_\R v_i)$. It follows
\begin{equation}
\label{lambda-mu}
(\lambda-\mu)(E)=\sum_j u_j\geq k
\end{equation}
because $u_j\geq 1$ for all $j$.
So, we see that $(\lambda-\mu)(E)$ is always an integer larger or 
equal to the order of the operator.

It follows immediately that $(\lambda-\mu)(E)=1$ implies that the operator is of first
order. To finish the proof, it remains to show that for a first order operator,
$(\lambda-\mu)(E)$ is $1$ (and not $2$).

Assume that $D$ is an operator of first order. This means that 
$\phi(1\otimes v_\mu)=\sum_j y_j\otimes v_j$ for $y_j\in\univ_1(\lieg_-)$ and again,
assume that $y_j$ are either constants or generators of negative root spaces and $v_i$
are weight vectors. All the terms $y_j\otimes v_j$ are of weight $\mu$, and therefore, 
$$\mu(E)(y_j\otimes v_j)=E (y_j\otimes v_j)=(\lambda(E)-[E,y_j])(y_j\otimes v_j)$$
so $[E,y_j]=(\mu-\lambda)(E)$ for all $j$ and it follows that all the $y_j$'s are from the
same graded components of $\lieg$. If $y_j\in\lieg_{-1}$, so $(\lambda-\mu)(E)=1$ and we are done.
Assume, for contradiction, that $y_j\in\lieg_{-k}$ for $k>1$.

Because $\sum_j y_j\otimes v_j\in\lieg_{-k}\otimes \V_\lambda$, choosing a basis $\{\tilde{v}_1,\ldots, \tilde{v}_m\}$ of
$\V_\lambda$, $\sum_j y_j\otimes v_j$ can be uniquely expressed as 
$\sum_{j=1}^k \tilde{y}_j\otimes \tilde{v}_j$ for some $\tilde{y}_j\in\lieg_{-k}$.
Because it is 
a homomorphic image of a highest weight vector in $M_\liep(\mu)$, 
it must be annihilated by all positive root spaces in $\lieg$, in particular, by any generator $x$ of a root space
in $\lieg_1$:
\begin{eqnarray*}
&&\mmmm x\cdot (\sum_j \tilde{y}_j\otimes \tilde{v}_j)=\sum_j x \tilde{y}_j\otimes_{\univ(\liep)}\tilde{v}_j=(\tilde{y}_jx+[x,\tilde{y}_j])
\otimes_{\univ(\liep)}\tilde{v}_j=\\
&&\mmmm =\sum_j \tilde{y}_j\otimes_{\univ(\liep)}x\cdot \tilde{v}_j+[x,\tilde{y}_j]\otimes_{\univ(\liep)}\tilde{v}_j=
\sum_j [x,\tilde{y}_j]\otimes \tilde{v}_j=0
\end{eqnarray*}
because $[x,\tilde{y}_j]\in \lieg_{-k+1}\subset\lieg_-$ and $x\cdot \tilde{v}_\lambda=0$.
Because $\tilde{v}_j$ forms a basis of $\V_\mu$, it follows that for each $j$, $[x,\tilde{y}_j]=0$
for all $x\in\lieg_1$. The grading must fulfill that $\lieg_{-1}$
generates $\lieg_-$ and it follows that $\lieg_{1}$ generates 
$\liep^+=\sum_{i=1}^k \lieg_{i}$. The Jacobi identity implies that
if $\tilde{y}_j$ commutes with $\lieg_1$, it commutes with all the $\liep^+$
as well. Let $\tilde{y}_j=\sum_i a_i y_{-\phi_i}$ where $y_{-\phi_i}$ is a generator
of the $-\phi_i$-root space. Define $x:=\sum_i a_i x_{\phi_i}$, where 
$x_{\phi_i}$ is a generator of the $\phi$-root space. We see that 
$x\in\lieg_{k}$ and $[x, \tilde{y}_j]=\sum_i a_i^2 [x_\phi, y_{-\phi}]\neq 0$ because
$[x_\phi,y_{-\psi}]\neq 0$ if and only if $\phi=\psi$. So, we obtained
a contradiction and therefore $y_j\in\lieg_{-1}$. 
\end{proof}

We see that the degree of an operator of first or second order can be recognized
immediately by the weights $\mu, \lambda$.

The operators of first order can be found easily. As we know, they are determined
by $P$-homomorphism $J_e^1(V(G/H))\to \W$, where $\V$ and $\W$ are the $P$-modules
inducing the associated vector bundles. As vector spaces,
$$J_e^1(V(G/H))\simeq \V\oplus (\lieg_-^*\otimes \V)$$
if we assign to a any $1$-jet in $e$ of $f\in C^\infty(G,\V)^P$ 
the pair $(v,\phi)\in \V\oplus (\lieg_-^*\otimes \V)$ defined by
$v=f(e)$ and $\phi(X)=(\omega^{-1}{X}f)(0)$ where $\omega$ is the Maurer-Cartan 
form trivializing the tangent bundle $TG$.

The action of $P$ on $\V\oplus (\lieg_-^*\otimes \V)$ was computed explicitly 
in \cite{Slovak_DrSc} and it follows that the action of $G_0$ (as a subgroup of
P) on $J_e^1\simeq \V\oplus(\lieg_-^*\otimes \V)$ 
is the usual action on the sum and product of $G_0$-modules. The tensor product
$\lieg_-^*\otimes \V$ decomposes, as a $G_0$-module, into a sum of irreducible 
$G_0$-modules $F_1\oplus\ldots\oplus F_k$ so that all irreducible modules in 
the decomposition have multiplicity one, see \cite{ss}. Therefore, any
homomorphism $\lieg_-^*\otimes \V$ to an irreducible module $F$ is a projection
onto one of these irreducible components. So, if none of the $F_j's$ is isomorphic
to $\V$, then there is a direct $G_0$-module decomposition 
$$J_e^1(V(G/P))\simeq \V\oplus F_1\oplus\ldots \oplus F_k$$ 
and each first order operator from $\Gamma(V(G/H))$ to somewhere
is determined by the projection to some $F_j$ (but the opposite statement is not true,
the projections need not to be $P$-homomorphisms in general).

We will usually consider complex representations of the (real) Lie group $P$.
If $\C$ is the trivial representation, $\Gamma(\C(G/P))$ is the space of smooth
functions on $G/P$, so $J_e^1\simeq \C\oplus \lieg_-^*\otimes_\R \C\simeq 
\C\oplus (\liep^+)^c$, where $(\liep^+)^c$ is the complexification
of the positive part $\liep^+=\oplus_{k\geq 1}\lieg_k$. In this case,
$J_e^1\simeq \C\oplus (\liep^+)^c$ also as $P$-module (if the action of $P$ is the adjoint action), so
there exists a differential operator of first order determined by 
$\C\oplus(\liep^+)^c \to (\liep^+)^c$ and it is the unique $G$-invariant first order operator
acting on functions on $G/P$. The vector bundle associated to the complexified adjoint representation $(\liep^+)^c$ of $P$
is the complex cotangent space $T^*(G/P)^c$ and the operator
is (up to multiple) the De Rham differential $d$, because we know that
$g\cdot df=d(g\cdot f)$ for all diffeomorphisms $g$, especially for all
$g\in G$ (it is easy to see that the action of $g$ on a section is
the pullback of the differential form or function in this case). 
In the language of generalized Verma modules, the De Rham differential is dual 
to a homomorphism $M_\liep((\lieg/\liep)^c)\to M_\liep(\C)$, where 
$\lieg/\liep$ is the representation of $\liep$ dual to $\liep^+$.

\subsection{Dirac operator and Clifford algebras}
\label{diracintro}
Fix a (positive definite) metric on $\R^n$, resp. $\C^n$. 
The real, resp. complex, Clifford algebra $\Cliff(n,\R)$, resp. $\Cliff(n,\C)$, is defined to be the 
associative algebra generated by $1$ and $\R^n$, resp. $\C^n$, with the relations 
$e_i\cdot e_j=-e_j\cdot e_i$ for $i\neq j$ and $e_i\cdot e_i=-1$, $e_j$ being an orthonormal basis of
$\R^n$, resp. $\C^n$. 
The spin group $\Spin(n)$, resp. $\Spin(n,\C)$ is the multiplicative subgroup 
generated by all $v_i\cdot v_j$, where 
$v_i, v_j$ are unit vectors. It has a natural action on $\R^n$ (resp. $\C^n$) defined by
$$(v_1\cdot\ldots\cdot v_{2n})\cdot v:=v_1\cdot\ldots\cdot v_{2n} \cdot v\cdot v_{2n}\cdot\ldots\cdot v_1\in\Cliff$$
and it is easy to check that the result is again a vector. This action is a rotation for any $x\in\Spin(n)$ 
and the vector space $\R^n$, resp. $\C^n$, with this action is the fundamental 
defining representation of $\Spin(n)$, resp. $\Spin(n,\C)$. 
It can be shown that $\Spin(n)$, resp. $\Spin(n,\C)$, is the unique connected and simply connected 
Lie group with Lie algebra $\so(n)$, resp. $\so(n,\C)$ 
(however, if we represent elements in $\so(n,\C)$ as matrices like
in section \ref{intro}, one has to take the defining condition of 
$\Cliff$ $v_1\cdot v_2+v_2\cdot v_1=2\beta(v_1,v_2)=
\sum v_1^i v_2^{n+1-i}$).

If $n$ is odd, the fundamental spinor representation $S$ of $\Spin(n,\C)$ can be realized as a minimal
left ideal in $\Cliff(n,\C)$.  If $n$ is even, there exists 
a unique minimal left ideal $S$ in $\Cliff(n,\C)$ that decomposes,
%
as a $\Spin(n,\C)$-module, into two irreducible modules $S^+$ and $S^-$ that 
are the fundamental spinor $\so(n,\C)$-modules with highest weights $\varpi_{n/2-1}$ and $\varpi_{n/2}$. 
These two representations of $\so(n,\C)$ are either self-dual (if $n/2$ is even) or dual to each other 
(if $n/2$ is odd). 
In both cases (odd and even), $S$ is called the {\it space of spinors}. By restriction, $S$ (resp. $S^+$, $S^-$)
are complex representations of the real Lie group $\Spin(n)$ as well.
The action of $\R^n\subset\Cliff(n,\R)\subset\Cliff(n,\C)$ on $S\subset\Cliff(n,\C)$ given by left multiplication 
in $\Cliff(n,\C)$ is called Clifford multiplication. 

The Dirac operator is defined via Clifford multiplication 
to be $$D:C^\infty(\R^n,S)\to C^\infty(\R^n,S),\,\,\,\, f\mapsto \sum_j e_j\cdot\partial_j f$$
and is independent of the choice of the orthonormal basis $e_j$.

\newpage
\section{Homomorphisms of generalized Verma modules}\label{gvm}
\subsection{True Verma module homomorphism}
Let $W$ be the Weyl group of $\lieg$ assuming a fixed choice of $\lieh$.
We define the {\it affine action} of $W$ on $\lieh$ to be
$w\cdot \mu:=w(\mu+\delta)-\delta$. Usually, we will denote by dot ($w\cdot \lambda$) 
the affine action and the usual action of $W$ without dot ($w\lambda$). 
The theorem of Harish-Chandra says that a homomorphism of Verma modules
$M(\mu)\to M(\lambda)$ can be nonzero only if $\mu$ and $\lambda$ are on the same affine orbit of $W$, i.e.
$\mu=w\cdot\lambda$ for some $w\in W$. It is well known (see, e.g. \cite[pp. 251]{Dixmier}) that
$\dim(\Hom(M(\mu),M(\lambda)))\leq 1$ and that each homomorphism $M(\mu)\to M(\lambda)$ is 
injective. So, there exists a nonzero homomorphism $M(\mu)\to M(\lambda)$ if and only if $M(\mu)$
is isomorphic to a unique submodule of $M(\lambda)$ isomorphic to $M(\mu)$. We will write
simply $M(\mu)\subset M(\lambda)$ in that case.

We can define the Bruhat ordering on the Weyl group by the relation $w\leq w'$, if there
exists a sequence $w=w_0, w_1, \ldots, w_k=w'$ so that $w_j=s_{\gamma_j} w_{j-1}$ for some
positive roots $\gamma_j$ and the length 
$l(w_j)>l(w_{j-1})$. 
We call {\it Hasse graph} the graph whose vertices are elements of $W$
and there is an arrow $w\to w'$ if and only if $w'=s_\gamma w$ for some $\gamma$ and 
$l(w')=l(w)+1$. It can be shown (\cite[pp. 265]{Dixmier}) that for any $w\in W$ and $\gamma\in\Phi$ either $w\leq s_\gamma w$
or $s_\gamma w\leq w$, i.e. the Hasse graph is the minimal partial ordered set of the Bruhat ordering.

\begin{lemma}
\label{bruhat2weight}
Let $\tilde\lambda\in P^{++}$, 
then $w\leq w'$ if and only if $w'(\tilde\lambda+\delta)\leq w(\tilde\lambda+\delta)$
\end{lemma}

\begin{proof}
See, e.g. \cite[pp. 117]{Roggenkamp}, or \cite[pp. 79]{Dlab}.
\end{proof}

\begin{lemma}
\label{W2lieh}
Let $w'=s_\gamma w$ for some positive root $\gamma$ and  $\lambda$ be a strictly dominant weight. Then 
$l(w')>l(w)$ if and only if $(w\lambda)(H_\gamma)>0$ 
($H_\gamma$ is the $\gamma$-coroot).
\end{lemma}

\begin{proof}
Because $\lambda$ is strictly dominant, it is in the interior of the fundamental Weyl chamber and 
$w\lambda$ is in the interior of another Weyl chamber. 
Therefore, $s_\gamma (w\lambda)=w\lambda-(w\lambda)(H_\gamma)\gamma\neq w\lambda$ and $(w\lambda)(H_\gamma)\neq 0$.
For $\lambda\in P^{++}+\delta$, it follows from the previous lemma that $w'\geq w$ if and only if
$w'\lambda=w\lambda-\lambda(H_\gamma) \gamma\leq w\lambda$, so $w\lambda(H_\gamma)$ is a nonnegative integer. 
The function $f(\mu)=\mu(H_\gamma)$ defined on the interior of the Weyl chamber containing $w P^{++}$ 
is a continuous nonzero function with positive values on integral weights, so it has positive values for all
$\mu=w\lambda$, where $\lambda$ is strictly dominant.
\end{proof}

The following theorem, proved by Bernstein-Gelfand-Gelfand and by Verma (\cite{bgg2, Verma}), 
describes the homomorphisms of Verma modules in terms of their highest weights:

\begin{theorem}
\label{truevermamodulesmap}
Let $\mu, \lambda\in\lieh^*$ (not necessarily integral).
There exists a nonzero homomorphism of Verma modules $M(\mu)\to M(\lambda)$ if and
only if there exists weights $\lambda=\lambda_0, \lambda_1,\ldots,\lambda_k=\mu$
so that $\lambda_{i+1}=s_{\beta_i}\cdot \lambda_{i}$ for some positive roots $\beta_i$
and $(\lambda_i+\delta)(H_{\beta_i})$ are  nonnegative integers for all $i$.
\end{theorem}

\begin{corollary}
Let $\tilde\lambda\in P^{++}$. Then $M(w'\cdot\tilde\lambda)\subset M(w\cdot\tilde\lambda)$ if and only if $w\leq w'$.
\end{corollary}

\begin{corollary}
\label{maxsubmod}
For $\tilde\lambda\in P^{++}$, the maximal (proper) submodule of $M(\tilde\lambda)$ is 
$\sum_{\alpha\in\Delta} M(s_\alpha\cdot\tilde\lambda)\subset M(\tilde\lambda)$.
\end{corollary}

The following lemma is a useful tool for comparing length of two elements in $W$:
\begin{lemma}
\label{grading2lengthbor}
Let  $M(\mu)\subset M(\lambda)$, $\mu=s_\gamma\cdot\lambda$ for some positive root $\gamma$ and suppose that
$\lambda\neq\mu$.
Let $E_\lieb$ be the grading element
associated to $(\lieg,\lieb)$.
Then $(\lambda-\mu)(E_\lieb)\in\N$.
\end{lemma}

\begin{proof}
Let $\oplus_j \lieg_j$ be the grading associated to $(\lieg, \lieb)$ and let 
$X_\gamma$ be the $\gamma$-root space generator. In this grading, all the positive root spaces are in
$\lieb^+$, so $X_\gamma\in\lieg_i$ for some $i>0$.
The defining equation for $E_\lieb$ implies $[E_\lieb,X_\gamma]=i X_\gamma=\gamma(E_\lieb) X_\gamma$, so $\gamma(E_\lieb)=i$.
We obtain 
\begin{eqnarray*}
&& \mmmm \mu(E_\lieb)=(s_\gamma (\lambda+\delta)-\delta)(E_\lieb)=
 (\lambda+\delta-(\lambda+\delta)(H_\gamma)\gamma-\delta)(E_\lieb)=\\
&& \mmmm =\lambda(E_\lieb)- i (\lambda+\delta)(H_\gamma)
\end{eqnarray*}
We know from Theorem \ref{truevermamodulesmap} that $(\lambda+\delta)(H_\gamma)$ is a nonnegative integer.
If $(\lambda+\delta)(H_\gamma)=0$, then
$\lambda=\mu$ contradicts the assumption. Therefore, both $i$ and $(\lambda+\delta)(H_\gamma)$ are positive integers
and $(\lambda-\mu)(E_\lieb)\in\N$.
\end{proof}

So, to construct the Hasse graph, we can find the Weyl group orbit of $\delta$ and for any two elements 
$\mu, \lambda$ on it such that $w'\delta=s_\gamma w \delta$, we can determine whether $w'\delta\geq w\delta$
(i.e. $w\leq w'$) or $w\delta\geq w'\delta$ (i.e. $w'\leq w$). In such a way, 
we construct the Bruhat ordering.

The following theorem describes the Jordan-H\"older series for Verma modules:

\begin{theorem}
\label{JH}
Let $0\subset A\subset B\subset M(\lambda)$ be $\lieg$-submodules and let $B/A=L(\mu)$ be an irreducible 
module with highest weight $\mu$. Then $M(\mu)\subset M(\lambda)$.
\end{theorem}

\begin{proof}
See, \cite[pp. 262]{Dixmier}.
\end{proof}

\begin{corollary}
\label{submodules2verma}
Let $\V_\mu$, resp. $\V_\lambda$, be highest weight $\lieg$-modules with highest weights $\mu$, resp. $\lambda$, 
and let $\V_\mu$ be a submodule of $\V_\lambda$.  Then $M(\mu)\subset M(\lambda)$.
\end{corollary}

\begin{proof}
Let $I$ be the ideal in $M(\lambda)$ so that $M(\lambda)/I\simeq \V_\lambda$ and let 
$\pi: M(\lambda)\to M(\lambda)/I$ be the
projection. Identifying $\V_\mu$ with a submodule of $M(\lambda)/I$, define $B:=\pi^{-1}(\V_\mu)$. Then 
$0\subset I\subset B\subset M(\lambda)$ is a sequence of $\lieg$-submodules with $B/I\simeq\V_\mu$. Because 
$\V_\mu$ is a highest weight module, take $J$ to be the maximal non-trivial submodule of $B/I$ so that
$(B/I)/J\simeq L(\mu)$ is irreducible with highest weight $\mu$. Take
$A:=\pi^{-1}(J)$. Then $L(\mu)\simeq (B/I)/(A/I)\simeq (B/A)$. So, $0\subset A\subset B\subset M(\lambda)$
is a sequence of $\lieg$-submodules and we can apply Theorem \ref{JH}.
\end{proof}

%

\subsection{Parabolic Hasse graph}
Let $\liep$ be a parabolic subgroup of $\lieg$, we denote by $W_p$ the
Weyl group of $\liep$. It is the subgroup of $W$ generated by the 
``uncrossed'' simple roots, i.e. simple roots with root spaces contained in $\lieg_0$
in the associated grading, see section \ref{parapara}. Denote the set of these simple roots by
$S$, i.e. $S=\Delta-\Sigma$.

Let us denote by $W^p$ the subset of $W$ consisting of those $w\in W$ so that 
$w\tilde\lambda$ is $\liep$-dominant for each $\lieg$-dominant weight $\tilde\lambda$.
It is shown in \cite[pp. 40]{Bjorner} that any $w\in W$ can be uniquely decomposed $w=w_p w^p$
where $w_p\in W_p$ and $w^p\in W^p$ and the length $l(w)=l(w_p)+l(w^p)$. A special case of this is
the following lemma that we will use later:

\begin{lemma}
\label{s_alpha w}
Let $w\in W^p$, $\alpha\in S$. Then $s_\alpha w\notin W^p$ and $l(s_\alpha w)=l(w)+1$.
\end{lemma}

%
The following statement shows that the Bruhat ordering is reasonable defined on $W^p$:
\begin{lemma}
\label{parabolicpath}
Let $w,w'\in W^p$, and $w\leq w'$ in the Bruhat ordering. Then there is a path 
$w\to w_1\to \ldots\to w_k=w'$ in the Hasse diagram so that all the elements $w_j\in W^p$.
\end{lemma}
\begin{proof}
\cite[pp. 45]{Bjorner}
\end{proof}

We define the {\it parabolic Hasse graph} to be the set $W^p$ of vertices with arrows $w\to w'$
if and only $w\to w'$ in $W$.  The lemma \ref{parabolicpath} says that the Bruhat ordering on $W$ induces 
the Bruhat ordering on $W^p$.

The following lemma is a tool for comparing length of two elements in $W^p$:

\begin{lemma}
\label{grading2length}
Let $E$ be the grading element for the pair $(\lieg, \liep)$ and let $w,w'\in W^p$, $w'=s_\gamma w$ and
$l(w')>l(w)$.
Then $w\delta(E)$ and $w'\delta(E)$ are integers or half-integers and 
$w\delta(E)-w'\delta(E)\in\N$. 
\end{lemma}

\begin{proof}
Let $\oplus_j \lieg_j$ be the grading associated to $(\lieg, \liep)$. It follows from lemma \ref{s_alpha w} 
that $\gamma\notin S$, so the $\gamma$-root space generator $X_\gamma\in\lieg_i$ for some $i>0$ and this implies 
$\gamma(E)=i\in\N$.
We obtain $w'\delta(E)=(s_\gamma w\delta)(E)=(w\delta-w\delta(H_\gamma) \gamma)(E)=w\delta(E)-i w\delta(H_\gamma)$.
Using lemma \ref{W2lieh}, we see that  $w\delta(H_\gamma)>0$. 
Recall that $\delta=\sum_j \varpi_j$ is 
integral and so is $w\delta$ for any $w\in W$.
Therefore, $w \delta(H_\gamma)\in\N$ and $w\delta(E)-w'\delta(E)=iw\delta(H_\gamma)\in\N$.

To show the half-integrality, recall that $\delta=1/2 \sum_{\beta\in\Phi^+} \beta$ and 
$w\delta(E)=1/2 \sum (w\beta)(E)$. For any root $\beta$, $w\beta$ is a root as well and 
$w\beta(E)$ is integral.
\end{proof}

So, one way how to construct the parabolic Hasse graph is to find all $\liep$-dominant weights on the
Weyl orbit of $\delta$ and for any $w'\delta=s_\gamma w\delta$, one may use the last lemma 
to determine whether $w'\geq w$ or $w\geq w'$.

\subsection{Standard and nonstandard homomorphism}
Let $\lieg$ be a semisimple Lie algebra and $\liep$ a parabolic subalgebra. Assume, as before,
that $S=\Delta-\Sigma$ is the set of uncrossed nodes, i.e. $S\subset\Delta$, $\alpha\in S \Leftrightarrow
Y_\alpha\in\liep$ where $Y_\alpha$ is the generator of the ($-\alpha$)-root space. 
We know from \ref{vermamodules-general} that for $\mu\in P_\liep^{++}$,
$M_\liep(\mu)$ is a highest weight module with highest weight $\mu$ and therefore
$$M_\liep(\mu)=M(\mu)/K_\mu$$
for some submodule $K_\mu$ of $M(\mu)$.

It follows from \ref{truevermamodulesmap} that for $\mu\in P_\liep^{++}$ and $\alpha\in S$, 
$M(s_\alpha \cdot \mu)\subset M(\mu)$ and the following lemma makes sense:

\begin{lemma}
\label{kernel}
The kernel of the projection $M(\mu)\to M_\liep(\mu)$ is 
$K_\mu=\sum_{\alpha\in S} M(s_\alpha \cdot \mu)$ (the sum of vector subspaces 
in $M(\mu)$).
\end{lemma}

\begin{proof}
Let $\{y_i\}$ be a basis of $\lieg_-$ and $\{Y_j\}$ be the set of generators of negative root spaces in $\lieg_0$.
Because $M(\mu)\to M(\mu)/K_\mu\simeq M_\liep(\mu)$ is a $\lieg$-homomorphism that maps $1\otimes 1$ to
$1\otimes v$, where $v\in \V_\mu$ is the highest weight vector in $\V_\mu$, it follows that the projection is
given by 
$$y_1\ldots y_k Y_1\ldots Y_l \otimes 1\mapsto y_1\ldots y_k\otimes (Y_1\ldots Y_l v)$$
The right side is zero if and only if $Y_1\ldots Y_l v$ is zero in $\V_\mu$.


The $\liep$-module $\V_\mu$ is an irreducible $\lieg_0^{ss}$-module, so it is a factor of the Verma module:
$\V_\mu\simeq M_{\lieg_0^{ss}}(\mu)/K'$ where $M_{\lieg_0^{ss}}(\mu)$ is the (true) Verma module for 
the Lie algebra $\lieg_0^{ss}$ and $K'$ is a maximal submodule of $M_{\lieg_0^{ss}}(\mu)$. $S$ is the set of simple
roots for $\lieg_0^{ss}$, so it follows from 
corollary \ref{maxsubmod} that $K'=\sum_{\alpha\in S} M_{\lieg_0^{ss}}(s_\alpha\cdot\mu)$. 
Using elementary representation theory, we obtain that $M_{\lieg_0^{ss}}(s_\alpha\cdot \mu)$
is a submodule of $M_{\lieg_0^{ss}}(\mu)$ generated by $Y_\alpha^{\mu(H_\alpha)+1}\otimes 1$ where $Y_\alpha$
is the generator of the root space of root $-\alpha$. The projection $M_{\lieg_0^{ss}}(\mu)\to \V_\mu$ is given
by $Y_1\ldots Y_l\otimes 1\mapsto Y_1\ldots Y_l v$. Therefore,
a vector $Y_1\ldots Y_l v$ is zero in $\V_\mu$ if and only if $Y_1\ldots Y_l\in\univ(\lieg_0)$ 
can be written down as a sum of vectors of type $Y_1'\ldots Y_m' Y_\alpha^{\mu(H_\alpha)+1}$,
summing over $\alpha\in S$. So, as a $\lieg_-$-module, $K_\mu$ is generated by 
the vectors $Y_\alpha^{\mu(H_\alpha)+1}\otimes 1$ 
what is the highest weight vector of $M(s_\alpha\cdot\mu)\subset M(\mu)$.
\end{proof}

Let $\mu,\lambda$ be $\liep$-dominant and $M(\mu)\subset M(\lambda)$. It follows from  
lemma \ref{s_alpha w} that  for $\alpha\in S$ $M(s_\alpha\cdot\mu)\subset M(\mu)$ and
$M(s_\alpha\cdot\lambda)\subset M(\lambda)$. This implies (\cite[pp. 252]{Dixmier})
that $M(s_\alpha\cdot\mu)\subset M(s_\alpha\cdot\lambda)$ and, consequently, 
$K_\mu\subset K_\lambda$ (representing both as submodules of $M(\lambda)$). Therefore,
for any $\mu,\lambda\in P_\liep^{++}$ and a homomorphism $i: M(\mu)\to M(\lambda)$ there is a well
defined factor homomorphism of generalized Verma modules $\tilde{i}: M_\liep(\mu)\to M_\liep(\lambda)$.

\begin{definition}
Generalized Verma module homomorphisms that are factors of true Verma module homomorphisms
are called standard.
\end{definition}

Standard homomorphisms are in general not injective, but a standard homomorphism \gvmhom
is unique up to multiple.

\begin{theorem}
\label{zeromap}
Let $\mu, \lambda\in P_\liep^{++}$, $i:M(\mu)\to M(\lambda)$ be a homomorphism of Verma modules.
Then the standard homomorphism $M_\liep(\mu)\to M_\liep(\lambda)$ is zero if and only if there exists 
$\alpha\in S$ so that $i(M(\mu))\subset M(s_\alpha \cdot\lambda)$.
\end{theorem}

\begin{proof}
If $M(\mu)\subset M(s_\alpha\cdot\lambda)$ for some $\alpha\in S$, it follows from lemma \ref{kernel} that the standard 
homomorphism is zero. On the other hand, suppose that the standard homomorphism $M_\liep(\mu)\to M_\liep(\lambda)$ is zero.
Lemma \ref{kernel} says that $M(\mu)\subset \sum_{j=1}^l M(s_{\alpha_j}\cdot\lambda)$ where $\alpha_1,\ldots,\alpha_l$
are the simple roots in $S$. Choose $i$ to be so that $M(\mu)\subset \sum_{j=1}^i M(s_{\lambda_j}\cdot \lambda)=:B$ but
$M(\mu)\nsubseteqq \sum_{j=1}^{i-1} M(s_{\alpha_j}\cdot\lambda)=:A$. The module $B/A$ is
a highest weight module with highest weight $s_{\alpha_i}\cdot\lambda$, because $M(s_{\alpha_i}\cdot\lambda)$ is such.
Let $v_\mu$ be the generator of $M(\mu)$ in $M(\lambda)$ (it is a weight vector of weight $\mu$) and
let $\pi: B\to B/A$ be the projection. 
Because $M(\mu)\subset B$, $M(\mu)\nsubseteqq A$, $\pi(v_\mu)\neq 0$ in $B/A$ and it generates some 
highest weight module $\V_\mu\subset B/A$ with highest weight $\mu$. Applying corollary \ref{submodules2verma} we 
obtain $M(\mu)\subset M(s_{\alpha_i}\cdot\lambda)$.
\end{proof}

\begin{theorem}
\label{genlepowski}
Let $\tilde{\lambda}+\delta$ be strictly dominant, $w,w'\in W^p$, $w\to w'$ in the parabolic Hasse graph 
for $(\lieg, \liep)$ 
and assume that $w\cdot \tilde\lambda,\, w'\cdot \tilde\lambda 
\in P_\liep^{++}$. Further, let $M(w'\cdot\tilde\lambda)\subset M(w\cdot\tilde\lambda)$. Then
the standard map $M_\liep(w'\cdot \tilde\lambda)\to M_\liep(w\cdot \tilde\lambda)$ is nonzero.
\end{theorem}

\begin{remark}
In \cite{Lepowski}, the theorem is given only for $\tilde\lambda\in P^{++}$ but the proof
works for non-integral $\tilde\lambda$ as well. Note, that for non-integral (and not $\liep$-dominant)
$\tilde\lambda$, $w\cdot\tilde\lambda$ may still be $\liep$-dominant and $\liep$-integral.

\end{remark}

\begin{proof}
Assume that the standard homomorphism is zero. Then
it maps $M(w'\cdot\tilde\lambda)$ into $K_{w\cdot\tilde\lambda}=
\sum_{\alpha\in S} M(s_\alpha w\cdot\tilde\lambda)).$
Lemma \ref{zeromap}  states that this happens exactly if there exists  
$\alpha\in S$ so that $M(w'\cdot\tilde\lambda)\subset M(s_\alpha w \cdot \tilde\lambda)$.
In such case, it follows from Theorem \ref{truevermamodulesmap} that there is sequence of submodules 
$$M(s_\alpha w\cdot\tilde\lambda)\supset 
M(s_{\gamma_1} s_\alpha w\cdot\tilde\lambda)\supset\ldots\supset M(s_{\gamma_j}\ldots 
s_{\gamma_1}s_\alpha w\cdot\tilde\lambda)=M(w'\cdot\tilde\lambda)$$
and $(s_{\gamma_{i-1}}\ldots s_{\gamma_1}s_\alpha w\cdot\tilde\lambda)(H_{\gamma_i})\in\N$.
It follows from lemma \ref{W2lieh} that $$s_{\gamma_j}\ldots s_{\gamma_1}s_\alpha w\geq s_{\alpha}w$$ 
in the Weyl group. Because $\tilde\lambda+\delta$ is strictly dominant (not in the interior 
of the fundamental Weyl chamber) the equality $w_1\cdot\tilde\lambda=w_2\cdot\tilde\lambda$ implies 
$w_1=w_2$, especially $w'=s_{\gamma_j}\ldots s_{\gamma_1} s_\alpha w$. So, $w'\geq s_\alpha w$.
Clearly, $w'\neq s_\alpha w$ and so $l(w')>l(s_\alpha w)$. This contradicts
$l(w')=l(w)+1=l(s_\alpha w)$ (lemma \ref{s_alpha w}).
\end{proof}

In most cases, the assumption $M(w'\cdot\tilde\lambda)\subset M(w\cdot\tilde\lambda)$ is superfluous.
For this to be satisfied, it suffices to show $(w\cdot\tilde\lambda)(H_\gamma)\in\N$, where $\gamma$
is the root so that $w'=s_\gamma w$.  Suppose that $w\cdot\tilde\lambda$ and $w'\cdot\tilde\lambda$
are both $\liep$-dominant. Then for any $\alpha\in S$\,\,\, 
$(w'\cdot\tilde\lambda)(H_\alpha)=
(w\cdot\tilde\lambda)(H_\alpha)-(w\cdot\tilde\lambda)(H_\gamma)\gamma(H_\alpha)\in \N$.
The first term $(w\cdot\tilde\lambda)(H_\alpha)$ is integral and to show that $(w\cdot\tilde\lambda)(H_\gamma)$
is integral as well, it suffices to show that there exists $\alpha\in S$ such that $\gamma(H_\alpha)\in\{-1,0,1\}$
(note that $\gamma\notin S$ because $s_\gamma\notin W_p$). In the orthogonal Lie algebras, this is 
always satisfied except the case $\lieg=B_n$, $S=\{\alpha_n\}$. Our cases of interest will be mainly the algebras
$B_n$ and $D_n$ with parabolic subalgebras determined by $\Sigma=\{\alpha_k\}$, i.e. $S=\Delta-\{\alpha_k\}$ (so, we should
be careful in the case $\lieg=B_2, \Sigma=\{\alpha_1\}$).

If $\tilde\lambda+\delta$ is on the wall of the fundamental Weyl chamber, we say that the Verma modules 
$M(w\cdot\tilde\lambda)$ have {\it singular character}. Unfortunately, Theorem \ref{genlepowski} cannot 
be generalized to this case. We will give a counterexample here.
Let us consider $\lieg=B_4$, $\Sigma=\{\alpha_2\}$. We will represent elements
of the Weyl group $w$ by the weight $w\delta$. 
Let $w\delta=\half[3,-1|7,5]$ and $w'\delta=(s_\gamma w)\delta=\half[1,-3|7,5]$ ($\gamma=[1,1|0,0]$). 
The grading element evaluation is $w\delta(E)=3/2-1/2=1$ and $w'\delta(E)=1/2-3/2=-1$. They are both strictly
$\liep$-dominant, so $w,w'\in W^p$.
The grading element
evaluation $\tilde{w}\delta(E)$ cannot be $0$ for any $\tilde{w}\in W^p$ in this case, so it follows from lemma \ref{parabolicpath}
and Theorem \ref{grading2length} that $l(w')=l(w)+1$. Let $w''=\half[-1,-3|7,5]$. We see that 
$w''=s_\beta w'$ for $\beta=[1,0|0,0]$ and similarly we can show that $l(w'')=l(w')+1$. Choosing
$\alpha=\alpha_1=[1,-1|0,0]$ and using lemma \ref{s_alpha w}, we see that the full Hasse graph of $\lieg$ 
contains the following square ($\eta=[0,1|0,0]$):

\begin{center}
\includegraphics{square.1}
\end{center}

Now let as choose a singular weight $\tilde\lambda+\delta=[3,2|1,0]$. 
Then $w(\tilde\lambda+\delta)=[1,0|3,2]$ and $w'(\tilde\lambda+\delta)=[0,-1|3,2]$
are both strictly $\liep$-dominant. But the $s_\beta$ fixes $w'(\tilde\lambda+\delta)$ and
so $w'(\tilde\lambda+\delta)=w''(\tilde\lambda+\delta)$: 
\begin{center}
\includegraphics{square.2}
\end{center}
Using Theorem \ref{truevermamodulesmap} we see that there exists a homomorphism of true Verma 
modules $M(w'\cdot\tilde\lambda)\to M(w\cdot\tilde\lambda)$ (note that $(w\cdot\tilde\lambda)(H_\gamma)=1$).
Similarly, we see that 
$M(w'\cdot\tilde\lambda)=M(w''\cdot\tilde\lambda)\subset M(s_\alpha w\cdot\tilde\lambda)$
and, applying Theorem \ref{zeromap}, the standard homomorphism of generalized Verma modules 
$M_\liep(w'\cdot\tilde\lambda)\to M_\liep(w\cdot\tilde\lambda)$ is zero.

However, there are indications that in such cases there still exists a homomorphism $M_\liep(w'\cdot\tilde\lambda)\to
M_\liep(w\cdot\tilde\lambda)$ but it is not standard. We 
can state the following conjecture, which holds in the known cases:

\begin{hypothesis}
\label{hypo}
Let $w\to w'$ in $W^p$, $M(w'\cdot\tilde\lambda)\subset M(w\cdot\tilde\lambda),\,\, w\cdot\tilde\lambda, 
w'\cdot\tilde\lambda\in P_\liep^{++}.$
Then there exist a (standard or nonstandard) nonzero
homomorphism $M_\liep(w'\cdot\tilde\lambda)\to
M_\liep(w\cdot\tilde\lambda)$.
\end{hypothesis}

\subsection{BGG graph}

\begin{definition}
\label{singularbgggraph}
Let $\lieg$ be a semisimple Lie algebra, $\liep$ a parabolic subalgebra, $\lambda\in\lieh^*$. We will define
the BGG graph for $(\lieg,\liep,\lambda)$ as an oriented graph, where the vertices are 
$\liep$-dominant and $\liep$-integral weights on the affine Weyl orbit of $\lambda$ and there is an arrow 
$\mu\to\nu$ if and only if the following conditions are satisfied:
\begin{enumerate}
\item{There exists a nonzero homomorphism $M_\liep(\nu)\to M_\liep(\mu)$}
\item{If there exist nontrivial homomorphisms 
$$M_\liep(\nu)=M_\liep(\xi_0)\to M_\liep(\xi_1)\to \ldots \to M_\liep(\xi_{j-1})\to M_\liep(\xi_j)=M_\liep(\mu)$$ then $j=1$ 
(by nontrivial we mean that they are nonzero and $\xi_i\neq \xi_{i+1}$).}
\end{enumerate}
\end{definition}

We know from section \ref{invdifop} that an arrow $\mu\to\nu$ in the BGG graph describes an invariant differential operator
acting between sections of associated vector bundles, induced by representations dual to $\V_\mu, \, \V_\nu$.

\begin{theorem}
\label{sing2regular}
If $\tilde\lambda\in P^{++}$ (in other words, $\tilde\lambda+\delta$ is integral and strictly dominant), 
then the BGG graph associated to $(\lieg,\liep,\tilde\lambda)$
is isomorphic, as a graph, to the Hasse graph for $(\lieg,\liep)$.
\end{theorem}
\begin{proof}
First, note that for $\tilde\lambda\in P^{++}$, $\tilde\lambda+\delta$ is in the interior
of the fundamental Weyl chamber and the map $w\in W^p\mapsto w\cdot\tilde\lambda$ is a bijection
from $W^p$ onto the set of $\liep$-dominant weights on the affine Weyl orbit of
$\tilde\lambda$. All these weights are integral. 
We will show that it is an isomorphism of the Hasse graph and
the BGG graph as well.

Let $\tilde\lambda\in P^{++}$ and $w\to w'$ in the  Hasse graph for $(\lieg, \liep)$.
Then $\mu:=w\cdot\tilde\lambda$
and $\nu:=w'\cdot\tilde\lambda$ are both $\liep$-dominant and theorem
\ref{genlepowski} states that there is a nonzero standard homomorphism 
$M_\liep(\nu)\to M_\liep(\mu)$.
Assume that there exist nonzero nontrivial homomorphisms 
$f_i: M_\liep(\xi_i)\to M_\liep(\xi_{i+1})$ for $i=0,\ldots,j-1$, $\xi_0=\nu$, $\xi_j=\mu$. 
Let $v_{\xi_i}$ be the
highest weight vector in $\V_{\xi_i}$. Then $1\otimes v_{\xi_i}$ is the highest weight 
vector in $M_\liep(\xi_i)$ and $f_i(1\otimes v_{\xi_i})$ is a weight vector in $M_\liep(\xi_{i+1})$
of weight $\xi_i$. Because $M_\liep(\xi_{i+1})$ is a highest weight module, it is easy to
see that each weight vector of weight $\beta$ fulfills $\beta\leq \xi_{i+1}$. Especially,
$\xi_{i}\leq \xi_{i+1}$. Because $\tilde\lambda\in P^{++}$, there exist a unique $w_i$ such that
$\xi_i=w_i\cdot\tilde\lambda$. The inequality $w_{i}\cdot\tilde\lambda\leq w_{i+1}\cdot\tilde\lambda$ implies
$w_i\geq w_{i+1}$ in the Bruhat ordering (lemma \ref{bruhat2weight}).
Because $\xi_i\neq \xi_{i+1}$, we get $w_i\neq w_{i+1}$ and
$w'>w_1>\ldots >w_j=w$. The assumption $w\to w'$ implies $j=1$. So, both conditions $(1)$ and $(2)$ are satisfied 
and there is an arrow $w\cdot\tilde\lambda\to w'\cdot\tilde\lambda$ in the BGG graph.



Assume, on the other side, that there is an arrow $\nu\to\mu$ in the BGG
graph. Let $w,w'\in W^p$ be the unique elements so that $\mu=w\cdot\tilde\lambda$ and $\nu=w'\cdot\tilde\lambda$.
Because there is a nonzero homomorphism $M_\liep(\nu)\to M_\liep(\mu)$, we again observe
that $\nu\leq\lambda$ and using $\ref{bruhat2weight}$, we see that $w\leq w'$. 
If $l(w')=l(w)+1$, then 
$w\to w'$ and we are done. Assume, for contradiction, that $l(w')>l(w)+1$. We know from
lemma \ref{parabolicpath} that there exist $w_j\in W^p$ such that $w=w_0\to w_1\to\ldots w_j=w'$.
Choosing the weighs $\xi_i=w_i\cdot\tilde\lambda$, we have
$$M(\nu)=M(\xi_0)\subsetneq M(\xi_1)\subsetneq\ldots\subsetneq M(\xi_j)=M(\mu),$$
$j>1$. Theorem \ref{genlepowski} implies that there are nonzero standard homomorphism of generalized
Verma modules $M_\liep(\xi_i)\to M_\liep(\xi_{i+1})$ for all $i$ and we obtain, from the definition
of singular Hasse graph that $j=1$ what contradicts $j>1$.

So, we proved that $w\to w'$ in $W^p$ if and only if 
$w\cdot\tilde\lambda\to w'\cdot\tilde\lambda$ in the BGG graph. 
\end{proof}

%

\begin{definition}
The {\it singular Hasse graph} for $(\lieg, \liep, \lambda)$ is a graph where vertices are $\liep$-dominant
and $\liep$-integral weights on the affine Weyl orbit of $\lambda$ and there is an arrow $\mu\to \nu$
if there exist $w,w'\in W^p$ so that $\mu=w\cdot\tilde\lambda$, $\nu=w'\cdot\tilde\lambda$ and $w\to w'$
in $W^p$ ($\tilde\lambda+\delta$ is some $\lieg$-dominant weight on the orbit of $\lambda+\delta).$
\end{definition}

In the regular case ($\tilde\lambda\in P^{++}$), the singular Hasse graph is clearly isomorphic to the Hasse graph.
Usually, we assume that $\lambda$ has singular character (i.e. its affine Weyl orbit does not contain any $\lieg$-dominant
weight) when we consider the singular Hasse graph. 
If $\tilde\lambda+\delta$ is on the wall of the fundamental Weyl chamber, then the weight $\mu$ on the affine orbit of $\lambda$ 
does not determine a unique $w$ such that $\mu=w\cdot\tilde\lambda$.
The conjecture \ref{hypo} says that the singular Hasse graph is a subgraph of the BGG graph. 

We will give an example of singular Hasse and BGG graph computation that will turn out to be important later. Let
$\lieg=D_4=\so(8,\C)$ and $\Sigma=\{\alpha_1\}$ (in the language of Dynkin diagrams, 
$\liep\simeq \dynkin\noroot{}\link\wr\whiterootdownright{}\whiterootupright{} \enddynkin$, i.e.
$\liep$ consists of the Cartan subalgebra and those root spaces, the determining roots of which could be written
as a linear combination of simple roots having nonnegative coefficient in the first simple
root $\alpha_1$).

Let us compute the Hasse graph for $(\lieg,\liep)$. If we represent a Weyl group element $w\in W^p$ by the
weight $w\delta$, we know that $W^p$ consists exactly of elements $w\in W$ such that $w\delta$ is strictly 
$\liep$-dominant. In this case, $\delta=[3|2,1,0]$ and a weight $[a|b_1,b_2,b_3]$ is strictly $\liep$-dominant 
if and only if $b_1>b_2>|b_3|$. It follows from lemma \ref{gradingelement} that the grading element 
evaluation is $[a|b_1,b_2,b_3](E)=a$. We state that (if we identify $w\sim w\delta)$ 
the regular Hasse graph for $(\lieg,\liep)$ has the following form:

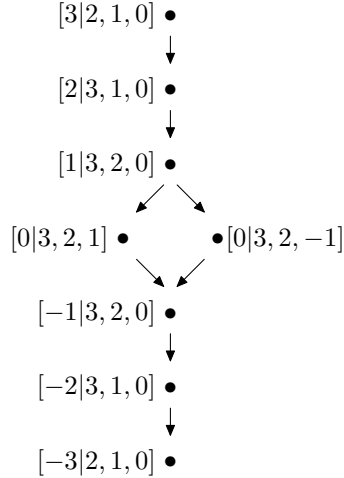
\begin{center}
\includegraphics{bgg.11}
\end{center}

We see that each arrow in this graph corresponds to a root reflection. Moreover, the grading element evaluation 
decreases by one in each arrow $w\delta\to w'\delta$ in this graph, so it follows from lemma \ref{grading2length} that $l(w')=l(w)+1$. 
These are all strictly $\liep$-dominant weights on the Weyl orbit of $\delta$, so we constructed the regular
Hasse graph for $(\lieg, \liep)$.

For any $\tilde\lambda\in P^{++}$, the same graph with
opposite arrows describes the structure of standard homomorphisms of generalized Verma modules on the affine
orbit of $\tilde\lambda$.

Further, let as consider a weight 
$\tilde\lambda:=\dynkin\noroot{0}\link\whiteroot{0}\whiterootdownright{-1}\whiterootupright{0} \enddynkin=\half[-1|-1,-1,1]$. 
Then $\tilde\lambda+\delta=\dynkin\noroot{1}\link\whiteroot{1}\whiterootdownright{0}\whiterootupright{1}\enddynkin=\half[5|3,1,1]$ 
is on the wall of fundamental Weyl chamber, so it has a singular affine orbit.

The $W^p$-orbit of $\tilde\lambda+\delta$ has the following form:

\begin{center}
\includegraphics{bgg.12}
\end{center}

The crosses $\times$ correspond to weights that are not strictly $\liep$-dominant, so, after subtracting $\delta$, 
the weights are not $\liep$-dominant and there are no associated generalized Verma modules for them. 
The nodes $\bullet$ are strictly $\liep$-dominant and the encircled weights coincide in this case.

We see a general fact that the affine orbit of a singular weight $\tilde\lambda$ (i.e. $\tilde\lambda+\delta$ is on the wall of 
the fundamental Weyl chamber) is smaller than the regular one: some weights
are ``glued together'' and some are not $\liep$-dominant. In our case, the only $\liep$-dominant weights on the affine Weyl
orbit of $\tilde\lambda$ are 
$\lambda=\dynkin\noroot{-3}\link\whiteroot{0}\whiterootdownright{0}\whiterootupright{1} \enddynkin=
\half[1|5,3,1]-\delta=\half[-5|1,1,1]$ and 
$\mu=\dynkin\noroot{-4}\link\whiteroot{0}\whiterootdownright{1}\whiterootupright{0} \enddynkin=\half[-1|5,3,-1]-\delta=\half[-7|1,1,-1]$.
We see immediately from the last diagram that the singular Hasse graph is $\lambda\to\mu$ in this case.

We will show that the standard homomorphism \gvmhom is nonzero. If it is zero, then lemma \ref{zeromap} says that
$M(\mu)\subset M(s_\alpha\cdot\lambda)$ for some $\alpha\in S$. So, there is a sequence $M(\mu)=M(\mu_0)\subset M(\mu_1)\subset\ldots
\subset M(\mu_k)=M(s_\alpha\cdot\lambda)\subset M(\lambda)$, the weights $\mu_j$ are increasing and connected by affine root 
reflections. 
Let $j$ be the largest number such that $\mu_0+\delta, \ldots, \mu_j+\delta$ have $-\half$ on the
first position.  So, $\mu_0+\delta,\ldots, \mu_j+\delta$ contain
a sign-permutation of $\half\{5,3,1\}$ on the $2$--$4$ positions. 
If $j>0$, $M(\mu)=M(w_p\cdot\mu_j)$ for $w_p\in W_p$ and $\mu_j+\delta=w_p^{-1}(\mu+\delta)$. 
For the $\liep$-dominant weight $\mu+\delta$, $\mu+\delta\geq w_p^{-1}(\mu+\delta)=\mu_j+\delta$ for any $w_p\in W_p$, but
this contradicts
the condition $\mu\leq\mu_j$. Therefore, $j=0$ and $(\mu_1+\delta)(E)=-\half+k$ for some $k\in\N$ (lemma \ref{grading2length}). 
The inequality $\mu_k+\delta\geq \mu_1+\delta$ implies $\half=(\mu_k+\delta)(E)\geq (\mu_1+\delta)(E)=-\half+k$.
This implies $k=1$ and $(\mu_1+\delta)(E)=\half$. Therefore, $\mu_1+\delta$ has $\half$ on the first position.
The only possible root reflection $s_\gamma$ such that
$s_\gamma(\half[-1|5,3,-1])=\half[1|\mathrm{something}]$ is clearly $\gamma=[1,0,0,1]$ and 
$s_\gamma(\mu+\delta)=\mu_1+\delta=\lambda+\delta$. This contradicts $M(\mu_1)\subset M(s_\alpha\cdot\lambda)$.

We showed that the standard map \gvmhom is nonzero and
the BGG graph for $(\lieg,\liep,\lambda)$ is $\lambda\to\mu$ as well.

\newpage
\section{Dirac operator in orthogonal parabolic geometry}
\label{diropinpargeo}

\subsection{Singular orbit of a particular weight}
We will generalize the example  
given at the end of chapter \ref{gvm}.

\begin{lemma}
\label{1dirgvmeven}
Let $\lieg=D_{n+1}=\so(2n+2,\C)$, $n\geq 2$ and
let $\liep$ be the parabolic subalgebra corresponding to $\Sigma=\{\alpha_1\}$, i.e.
$\liep=\dynkin\noroot{}\link\wr\link\ldots\link\wr\whiterootdownright{}\whiterootupright{} \enddynkin$. Then
choosing
$$\lambda=\frac{1}{2}[-2n+1|1,1,\ldots,1]\quad{\mathrm{and}}\quad \mu=\frac{1}{2}[-2n-1|1,1,\ldots,1,-1],$$
or, in the language of Dynkin diagrams,
$$\lambda=\dynkin\noroot{-n}\link\whiteroot{0}\link\ldots\link\whiteroot{0}\whiterootdownright{0}\whiterootupright{1} \enddynkin,\quad
\mu=\dynkin\noroot{-n-1}\link\whiteroot{0}\link\ldots\link\whiteroot{0}\whiterootdownright{1}\whiterootupright{0} \enddynkin,$$
there exists a unique (up to a multiple) nonzero standard homomorphisms
of generalized Verma modules $$M_\liep(\mu)\to M_\liep(\lambda)$$ and the BGG graph for $(\lieg,\liep,\lambda)$ is
$\lambda\to\mu$.

Similarly, for the weights
$$\lambda'=\frac{1}{2}[-2n+1|1,1,\ldots,1,-1]\quad{\mathrm and}\quad \mu'=\frac{1}{2}[-2n-1|1,1,\ldots,1,1],$$
there also exists a unique (up to multiple) nonzero homomorphisms
of generalized Verma modules $$M_\liep(\mu')\to M_\liep(\lambda')$$ and the BGG graph is $\lambda'\to\mu'$.
\end{lemma}

\begin{proof}
First, we find all strictly $\liep$-dominant weights on the Weyl orbit of 
$\lambda+\delta=\half[-2n+1|1,\ldots,1]+[n|n-1.\ldots,1,0]=\half[1|2n-1,\ldots,3,1]$. 
The condition of strict $\liep$-dominance implies that the weights must be of the form 
$[a|b_1,\ldots,b_n]$, where $b_1>\ldots>b_{n-1}>|b_n|$ and $\{a,b_1,\ldots,b_n\}$ is a sign-permutation of
$\{2n-1,\ldots,3,1,1\}$ with an even number of minuses, 
so $\lambda+\delta=\half[1|2n-1,\ldots,3,1]$ and $\mu+\delta=\half[-1|2n-1,\ldots,3,-1]$ 
are the only possibilities. 

We see that $\mu=s_\gamma\lambda=\lambda-\gamma$ for $\gamma=[1,0\ldots,0,1]$ and we know from theorem
\ref{truevermamodulesmap} that for true Verma modules, $M(\mu)\subset M(\lambda)$. Suppose that the
standard map \gvmhom is zero, so there exists $\alpha\in S$ so that $M(\mu)\subset M(s_\alpha\cdot\lambda)$.
So, there is a sequence $M(\mu)=M(\mu_0)\subset M(\mu_1)\subset\ldots
\subset M(\mu_k)=M(s_\alpha\cdot\lambda)\subset M(\lambda)$, the weights $\mu_j$ are increasing and connected 
by affine root reflections. 
Let $j$ be the largest number such that $(\mu_0+\delta)(E)=-\half, \ldots, (\mu_j+\delta)(E)=-\half$.
So, $\mu_0+\delta,\ldots, \mu_j+\delta$ are of the form $\half[-1|b_1,\ldots,b_n]$, where $\{b_1,\ldots,b_n\}$ is
a sign-permutation of $\{2n-1,\ldots,3,-1\}$. It follows that these weights are on the $W_p$-orbit of 
$\mu$. But $\mu$ is $\liep$-dominant, so all the other weights on its $W_p$-orbit are smaller. Therefore,
$j=0$. Further, $(\mu_1+\delta)(E)$ must be $\half$, because $(\lambda+\delta)(E)$ is $\half$ and the weights
$\mu_j$ are increasing with increasing $j$. The only root reflection 
which maps $\half[1|2n-1,\ldots,3,1]$ to $\half[-1|b_1,\ldots,b_n]$ 
is clearly the reflection $s_\gamma$ and $\mu_1=\lambda$. We obtained
$$M(\mu_1)\subsetneq M(s_\alpha\cdot\lambda)\subsetneq M(\lambda)=M(\mu_1)$$
which is a contradiction. Therefore, the standard homomorphism \gvmhom is nonzero.

The proof for $\lambda', \mu'$ is similar.
\end{proof}

\begin{lemma}
\label{1dirgvmodd}
Let $\lieg=B_{n+1}=\so(2n+3,\C)$, $n\geq 1$ and
$\liep$ be the parabolic subalgebra corresponding to $\Sigma=\{\alpha_1\}$,
$\liep=\dynkin\noroot{}\link\wr\link\ldots\link\wr\llink>\wr\enddynkin$. Then
choosing
$$\lambda=[-n|\half,\ldots,\half]\quad{\mathrm{and}}\quad \mu=\frac{1}{2}[-n-1|\half,\ldots,\half],$$
or, in the language of Dynkin diagrams,
$$\lambda=\dynkin\noroot{-n-{1\over 2}}\link\whiteroot{0}\link\ldots\link\whiteroot{0}\llink>\whiteroot{1}\enddynkin,\quad
\mu=\dynkin\noroot{-n-{3\over 2}}\link\whiteroot{0}\link\ldots\link\whiteroot{0}\llink>\whiteroot{1}\enddynkin,$$
there exists a unique (up to a multiple) nonzero standard homomorphisms
of generalized Verma modules $$M_\liep(\mu)\to M_\liep(\lambda)$$ and the BGG graph for $(\lieg,\liep,\lambda)$ is
$\lambda\to\mu$.
\end{lemma}

\begin{proof}
In this case, $\delta=\half[2n+1,\ldots,3,1]$ and the weight 
$\lambda+\delta=[\half|n,\ldots,2,1]$ is on the orbit of $\tilde\lambda+\delta=[n|n-1,\ldots,2,1,\half]$
that is regular but not integral. A weight $[a|b_1,\ldots,b_k]$ is strictly $\liep$-dominant and $\liep$-integral
if and only if $b_1>\ldots>b_n>0$ and the $b_i$'s are all integers or all half-integers.
 The only strictly $\liep$-dominant and $\liep$-integral weights on its orbit are $\lambda+\delta$ and 
$\mu+\delta=[-\half|n,\ldots,2,1]$. They are connected by  reflection with respect to the root $\epsilon_1$. There is a unique
$w$ taking $\tilde\lambda+\delta$ to $\lambda+\delta$, namely $w\delta=\half[1|2n-1,\ldots,5,3]$ and a unique 
$w'$ taking $\tilde\lambda+\delta$ to $\mu+\delta$, namely $w'\delta=\half[-1|2n-1,\ldots,5,3]$. The grading element evaluation
differs by $1$ and $w'=s_{\epsilon_1}w$, so it follows from lemma \ref{grading2length} that 
there is an arrow $w\to w'$ in the regular Hasse graph for $(\lieg, \liep)$. 
Because $\mu=\lambda-\epsilon_1$, it follows from Theorem \ref{truevermamodulesmap} that 
$M(\mu)\subset M(\lambda)$. All the conditions of Theorem \ref{genlepowski} are satisfied, so
there exists a nonzero standard homomorphism $M_\liep(\mu)\to M_\liep(\lambda)$.
\end{proof}

\subsection{The real version} 
\label{therealversion}
Let us now suppose that $\lieg=\so(n+1,1;\R)$ is the real Lie algebra consisting
of matrices invariant with respect to the quadratic form  $x_0x_{n+1}+\sum_{j=1}^{n} x_j^2$ and $\liep$ is
the (real) parabolic subalgebra killing a chosen line in the null-cone. In matrices, elements of $\lieg$ are represented as
$$\left(
\begin{tabular}{c|c|c}
$\R$ & $\lieg_{1}$ & $0$ \\
\hline
$\lieg_{-1}$ & $\so(n)$ & $\lieg_{1}$ \\
\hline
$0$ & $\lieg_{-1}$ & $\R$
\end{tabular}
\right)
$$
The negative part $\lieg_{-1}\simeq \R^{n}$ is the fundamental defining representation of $\so(n)\subset\lieg_0$ via the adjoint action 
and $\lieg_0=\so(n)\oplus \R$.

We assume that $\lieg$ is naturally embedded into its
complexification $\lieg^c=\so(n+2,\C)$ and that the Cartan subalgebra, positive roots and fundamental weights of the complexification 
are given like before. The complexification of $\liep$ is exactly the parabolic subalgebra corresponding to 
$$\dynkin\noroot{}\link\wr\link\ldots\link\wr\whiterootdownright{}\whiterootupright{} \enddynkin$$ in case $n$ is even and
$$\dynkin\noroot{}\link\wr\link\ldots\link\wr\llink>\wr\enddynkin$$ in case $n$ is odd.

Let $\mu, \lambda$ be elements of $\lieh^c$ from lemma \ref{1dirgvmeven} or \ref{1dirgvmodd},
$\V_\lambda$, resp. $\V_\mu$, be representation of $\liep^c$ with highest weight 
$\lambda$, resp. $\mu$. Via restriction, they are (complex) 
representations of the real form $\liep$ as well. 

As vector spaces, the generalized Verma modules for real Lie algebras and complex inducing representation are 
isomorphic to the generalized Verma modules for the complex Lie algebras:
\begin{equation}
\label{realvermamodule}
M_{\liep^c}(\mu):=\univ(\lieg^c)\otimes_{\univ(\liep^c)}\V_\mu\simeq
\univ(\lieg)\otimes_{\univ(\liep)}\V_\mu=:M_\liep(\mu).
\end{equation}
The first product is over the complex universal enveloping algebra $\univ(\liep)$ 
(that is a factor of the complex tensor algebra)
and the second is over the real universal enveloping algebra $\univ(\liep)$.

As vector spaces, $$M_{\liep^c}(\mu)\simeq\univ(\lieg^c_-)\otimes_\C \V_\mu\simeq \univ(\lieg_-)\otimes_\R \V_\mu\simeq M_\liep(\mu)$$
and the middle isomorphism is given by
$$((y_1\otimes c_1)(y_2\otimes c_2)\ldots (y_k\otimes c_k))\otimes v \,\,\mapsto \,\, (y_1 y_2\ldots y_k)\otimes (c_1 c_2 \ldots c_k)v$$
for $y_j\in\lieg_-$, $c_j\in\C$ and $v\in\V_\mu$.

This vector space homomorphism is compatible with the action of $\lieg\subset\lieg^c$ on both spaces, i.e. it is a
$\lieg$-isomorphism. The same is true for any $\liep^c$-dominant and $\liep^c$-integral weight $\mu'$, 
in particular $M_{\liep^c}(\lambda)\simeq M_\liep(\lambda)$.

Because we know from the previous section that there exists a unique (up to multiple) $\lieg^c$-homomorphism
$M_{\liep^c}(\mu)\to M_{\liep^c}(\lambda)$, it follows that there exists a unique (up to multiple) nonzero homomorphism
of the real generalized Verma modules $M_\liep(\mu)\to M_\liep(\lambda)$ in this case as well.

\subsection{Description of the differential operator}\label{onedirak} Let $\lieg,\liep$ be as in the last section, 
$G=\Spin(n+1,1)$ the real Lie group with Lie algebra $\lieg$, $P$ the parabolic subgroup of $G$ fixing a 
line in the null-cone so that the Lie algebra of $P$ is $\liep$.  Let $\V_\lambda$ and $\V_\mu$ be representations 
of $\liep^c$ as in lemma \ref{1dirgvmeven} (in case $n$ is even) or lemma \ref{1dirgvmodd} (in case $n$ is odd). 
By restriction, they are complex representations of the real Lie algebra $\liep$ and of the real Lie group $P$ as well
($P$ is simply connected). The duality between homomorphisms of generalized Verma modules and invariant 
differential operators yields a nonzero invariant 
differential operator $D:\Gamma(G\times_P \V_{\lambda}^*)\to \Gamma(G\times \V_{\mu}^*)$ acting between the
spaces of smooth sections.

\begin{lemma}
The operator $D$ is of first order.
\end{lemma}
\begin{proof}
We know from section \ref{parapara} that for a weight $\lambda=[a|b_1,\ldots,b_l]$ the evaluation on the grading element
is $\lambda(E)=a$. Using Theorem \ref{grading2degree} and the definition of $\mu, \lambda$ yields the order to be
$$(\half[2k+1|\ldots,3,1]-\half[2k-1|\ldots,3,-1])(E)=\frac{2k+1}{2}-\frac{2k-1}{2}=1$$ in the even case $n=2k$ and 
$$([\half|2n-1,\ldots,2,1]-[-\half|2n-1,\ldots,2,1])(E)=\half-(-\half)=1$$
in the odd case.
\end{proof}

To describe the operator $D$, we will embed the vector space $\lieg_-$ into $G/P$ by the exponential 
$$i: \lieg_- \to G/P,\quad y\mapsto \exp(y)P.$$
It is known that $i$ is injective and the image $i(\lieg_-)$ is an open dense subspace in $G/P$.

We will identify $\lieg_-$ with its image under $i$. To any section $s\in\Gamma(G\times_P \V)$
such that $s(gP)=[g,v]_P$ we can assign a $\V$-valued function $f$ on $\lieg_-$ defined by
\begin{equation}
\label{identification}
f: \lieg_-\to \V, y\mapsto v, \quad\hbox{where}\quad s(i(y))=[\exp(y),v]_P.
\end{equation}
On the other hand, to any such function with compact support we can assign the corresponding section $s$ by
$s(\exp(y)P)=[\exp(y), f(y)]_P$ for $y\in\lieg_-$ and $s=0$ on $G/P-\exp(\lieg_-)P$.

The space $\lieg_-$ is endowed with a basis 
$$\left(
\begin{tabular}{c|c|c}
$0$ & $0$ & $0$ \\
\hline
$e_j$ & $0$ & $0$ \\
\hline
$0$ & $-e_j^T$ & $0$
\end{tabular}
\right)
$$
where $e_j=(0,\ldots,0,1,0,\ldots,0)^T$ is the $j$-th vector of the standard basis of $\R^{n}$.

In the even case, the $\V_\lambda$ and $\V_\mu$, as a representation of $\lieg_0^{ss}=\so(n)$, 
are the fundamental spinor representation $S^+$ and $S^-$.
Assume that  $S=S^+\oplus S^-$ is realized as a subspace of 
the Clifford algebra $\Cliff(n,\C)$, so that $S$ is a minimal left ideal and 
$S^+$ and $S^-$ are its $\so(n,\C)$-invariant subspaces.
In case $n$ is odd, $\V_\mu\simeq \V_\lambda\simeq S$ as a $\lieg_0^{ss}$-module, where 
$S$ is the unique fundamental spinor representation of $\so(n,\C)$, again
realized as a minimal left ideal in $\Cliff(n,\C)$. In both cases, the 
Clifford multiplication $\R^{n}\otimes S\to S$ is well defined.

Let us denote by $\cV_{\lambda^*}$ the spinor bundle $G\times_P (\V_\lambda^*\oplus \V_{\lambda'}^*)$ 
in case $n$ is even (see lemma \ref{1dirgvmeven} for definition of $\lambda, \lambda'$) or 
$G\times_P \V_{\lambda}^*$ in case $n$ is odd (see \ref{1dirgvmodd} for the definition of $\lambda$ in this case). 
Similarly, define $\cV_{\mu^*}$ to be the vector bundle $G\times_P (\V_\mu^*\oplus\V_{\mu'}^*)$ for $n$ even
and $G\times_P \V_{\mu}^*$ for $n$ odd. We can consider the operator $D$ as acting between $\Gamma(\cV_{\lambda^*})$
and $\Gamma(\cV_{\mu^*})$.
To each section $s$ of $\cV_{\lambda^*}$, resp. $\cV_{\mu^*}$, we assign a function $f:\lieg_-\to S$ defined by $(\ref{identification})$.
\begin{theorem}
\label{whydirakisdirak}
Let $s\in \Gamma(\cV_{\lambda^*})$ and $s'\in\Gamma(\cV_{\mu^*})$ 
are sections and $f,f':\lieg_-\to S$ the spinor valued
functions corresponding to $s$ and $s'$ under the above identification. Assume that $s'=Ds$.
Then $f'=\sum_j e_j \cdot \partial_{j}f$, $e_j$ being the standard basis of $\R^{n}$.
\end{theorem}

{\bf Remark.} Therefore, we call the operator $D$ ``Dirac operator''.
\begin{proof}
Recall that $\lieg_0\simeq \so(n)\oplus \R$ and the group $G_0$ corresponding to $\lieg_0$
contains $\Spin(n)$. Let us choose $g\in \Spin(n)\subset G_0\subset P\subset G$ and $y\in\lieg_-$.
The identity $g\exp(y)g^{-1}=\exp(\Ad(g)y)$ yields $g\exp(y)P=\exp(\Ad(g)y)P$ and therefore
the action of $\Spin(n)$ on $\lieg_-$ (identified with a subset of $G/P$) is the 
adjoint action. It is easy to check that $\lieg_-$ with the adjoint action of $\so(n)\subset\lieg_0$
is the fundamental defining representation of $\so(n)$ and so $\lieg_-\simeq \R^n$ 
is the fundamental defining representation of $\Spin(n)$ as well.

We know from construction that the operator $D$ is $G$-invariant, it means it commutes with the action of $G$
on sections. The action of $g\in G$ on $s$ can be described by
$$g\cdot s(x):=g(s(g^{-1}\cdot x)),$$
where on the right hand side, the actions of $g$ are on $G/P$ and on the spinor bundle $G\times_P S$.
So, if $s$ is described by a function $f$ on $\lieg_-$ and $g\in\Spin(n)$, the 
equation reads  
\begin{eqnarray*}
&&\mmmm (g\cdot s)(\exp(y) P)= g([\exp(\Ad(g^{-1})(y))P, f(\Ad(g^{-1})(y))]_P)=\\
&&\mmmm =[g\exp(\Ad(g^{-1})(y))P, f(\Ad(g^{-1})(y))]_P=\\
&&\mmmm =[g\exp(\Ad(g^{-1})(y))g^{-1}, g(f(\Ad(g^{-1})(y)))]_P=\\
&&\mmmm =[\exp(y), g(f(\Ad(g^{-1})(y)))]_P
\end{eqnarray*}
so we see that $(g\cdot f)(y)=g(f(g^{-1}\cdot y))$
is the usual action of $\Spin(n)$ on spinor valued functions. So, 
$D$ is invariant (acting on spinor valued functions on $\lieg_-\simeq \R^n$) with respect to the 
spin group (considering $\R^n$ and $S$ as the defining and spinor representations).

Further, choose $g=\exp(y_0)$ for some $y_0\in \lieg_-$. Then $\exp(y_0) (\exp(y))P=\exp(y+y_0)P$ because $\lieg_-$ is
commutative. We see that the action of $\exp(y_0)$ on $\lieg_-$ (identified with a subset of $G/P$) is just $y_0$-addition.
Further, the action of $\exp(y_0)$ on $f$ is $(\exp(y_0) (f))(y)=f(y-y_0)$, because, if $f$ represents the section
$s$, 
\begin{eqnarray*}
&&\hspace{-20pt} (\exp(y_0) s)(\exp(y)P)=\exp(y_0)(s(\exp(y-y_0)))=\exp(y_0)[\exp(y-y_0),f(y-y_0)]_P=\\
&&\hspace{-20pt} =[\exp(y),f(y-y_0)]_P.
\end{eqnarray*}

This means that the operator $D$ is invariant with respect to $\Spin(n)$ and with respect to $\R^{n}$ (translation)
as well. 

The invariance with respect to $\Spin$ implies that 
$Df(0)=\pi(\sum_j \epsilon_j \otimes \frac{\partial f}{\partial x_j}(0))$
for some $\Spin(n)$-homomorphism $$\pi: (\R^{n})^*\otimes S\to S$$ 
where $\{\epsilon_j\}$ is the dual basis to the standard basis $\{e_j\}$ of $\R^n$ and $\{x_j\}$ are the
coordinates on $\R^n$ with respect to the basis $e_j$.

The invariance with respect to translations
implies that for each $y\in \R^{n}$ $Df(y)=\pi(\sum_j \epsilon_j\otimes \frac{\partial f(y)}{\partial x_j})$.

We know from \cite{sb} that such a homomorphism $\pi$ is unique (up to multiple), namely 
$\pi(\sum_j \epsilon_j\otimes v_j)=e_j\cdot v_j$ ($(\cdot )$ is the Clifford multiplication),  so  $D$ is the Dirac operator
$D=\sum_j e_j\cdot \partial_{j}$.
\end{proof}

Note, that the operator $D$ acts between sections of vector bundles associated to representations 
$$\V_\lambda^*\oplus\V_{\lambda'}^*\simeq \dynkin\noroot{\frac{n}{2}-1}\link\whiteroot{0}\link\ldots\link\whiteroot{0}\whiterootdownright{1}\whiterootupright{0} \enddynkin\oplus
\dynkin\noroot{\frac{n}{2}-1}\link\whiteroot{0}\link\ldots\link\whiteroot{0}\whiterootdownright{0}\whiterootupright{1} \enddynkin$$
and 
$$
\V_\mu^*\oplus\V_{\mu'}^*\simeq \dynkin\noroot{\frac{n}{2}}\link\whiteroot{0}\link\ldots\link\whiteroot{0}\whiterootdownright{0}\whiterootupright{1} \enddynkin\oplus
\dynkin\noroot{\frac{n}{2}}\link\whiteroot{0}\link\ldots\link\whiteroot{0}\whiterootdownright{1}\whiterootupright{0} \enddynkin
$$
in the even case and between sections of vector bundles associated to
$$\V_\lambda^*\simeq \dynkin\noroot{\frac{n}{2}-1}\link\whiteroot{0}\link\ldots\link\whiteroot{0}\llink>\whiteroot{1}\enddynkin\quad\and\quad
\V_\mu^*=\dynkin\noroot{\frac{n}{2}}\link\whiteroot{0}\link\ldots\link\whiteroot{0}\llink>\whiteroot{1}\enddynkin$$
in the odd case (this can be derived easily using the definition of the weights $\lambda, \mu$ in \ref{1dirgvmeven},
resp. \ref{1dirgvmodd}, considering that here, the rank is not $n$ but $n/2$, resp. $n/2-1/2$, and using the fact
that the action of the center of $\lieg_0$ on the dual representation $\V^*$ is just the negative of its action on $\V$).

Further, note that the bundle $G\to G/P$ can be reduced to the standard spin structure  $\Spin(n+1)\to \Spin(n+1)/\Spin(n)$ 
on the sphere by choosing a Weyl structure on $G\to G/P$ (see
\cite{Cap-Weyl}) and factorizing the center of $G_0$. 
Under this reduction, the operator $D: \Gamma(\Spin(n+1)\times_{\Spin(n)} S)\to \Gamma(\Spin(n+1)\times_{\Spin(n)} S)$
is the Dirac operator on the sphere, associated to the usual Euclidean metric, 
see \cite{Friedrich} for details. The Dirac operator $D$ can be defined also on curved analogues of this and 
is still conformally invariant.

\subsection{More Dirac operators}

Consider now a pair of real Lie algebras ($\lieg,\liep$) described by the Dynkin diagram 
$$\dynkin\whiteroot{}\link\ldots\link\whiteroot{}\link\whiteroot{}\link\noroot{}\link\whiteroot{}\link\ldots\link\whiteroot{}\whiterootdownright{}\whiterootupright{}\enddynkin$$
or, equivalently, $\lieg=D_{k+n}=\so(2(k+n))$, $\Sigma=\{\alpha_k\}$.
Using the formalism of chapter \ref{intro}, the matrices in $\lieg$ are graded like this:
$$\left(
\begin{tabular}{c|c|c}
$\lieg_0$ & $\lieg_{1}$ & $\lieg_{2}$ \\
\hline
$\lieg_{-1}$ & $\lieg_0$ & $\lieg_{1}$ \\
\hline
$\lieg_{-2}$ & $\lieg_{-1}$ & $\lieg_0$
\end{tabular}
\right)
$$
\begin{theorem} 
\label{kdirgvmeven}
Choosing
$$\lambda=\frac{1}{2}[-2n+1,\ldots,-2n+1|1,\ldots,1]$$ and $$\mu=\frac{1}{2}[-2n+1,\ldots,-2n+1,-2n-1|1,\ldots,1,-1],$$
or, in the language of Dynkin diagrams,
\begin{equation}
\label{lambdaweight}
\lambda=\dynkin\whiteroot{0}\link\ldots\link\whiteroot{0}\link\noroot{-n}\link\whiteroot{0}\link\ldots\link\whiteroot{0}\whiterootdownright{0}\whiterootupright{1} \enddynkin
\end{equation}
and
\begin{equation}
\label{muweight}
\mu=\dynkin\whiteroot{0}\link\ldots\link\whiteroot{0}\link\whiteroot{1}\link\noroot{-n-1}\link\whiteroot{0}\link\ldots\link\whiteroot{0}\whiterootdownright{1}\whiterootupright{0} \enddynkin
\end{equation}
there exists a unique (up to a multiple) nonzero standard homomorphisms
of generalized Verma modules $$M_\liep(\mu)\to M_\liep(\lambda)$$ 
and the corresponding (dual) differential operator is of first order. Analogous statement holds for
the weights $\mu'$ and $\lambda'$ having interchanged $0$ and $1$ over the last positions in the Dynkin diagram.
\end{theorem}

\begin{proof}
Notice that $\lambda+\delta=\half[2k-1,\ldots,3,1|2n-1,\ldots,3,1]$ and $\mu+\delta=\half[2k-1,\ldots,3,-1|2n-1,\ldots,3,1]$.
Clearly, they are both strictly $\liep$-dominant and $\mu=s_\gamma\cdot\lambda=\lambda-\gamma$ for 
$\gamma=[1,0,\ldots,0,1]$, so $M(\mu)\subset M(\lambda)$. Suppose, for contradiction, that the standard map \gvmhom is zero.
Again we find a sequence of submodules $M(\mu)=M(\mu_0)\subset\ldots\subset M(s_\alpha\cdot\lambda)$. Similarly as in the 
proof of \ref{1dirgvmeven}, we can show that $\mu_1(E)=\mu_0(E)+1=\lambda(E)$. Let $\gamma_1$ be the reflection so that
$\mu_1+\delta=s_{\gamma_1}(\mu_0+\delta)$. Grading element evaluation is the sum of the first $k$ coordinates in the
expression of the weights, so $s_{\gamma_1}$, applied to $\half[2k-1,\ldots,3,-1|2n-1,\ldots,3,-1]$ increases one
of the first $k$ coordinates by one. If $\gamma_1=\gamma$, it follows $\mu_1=\lambda$ and we obtain a contradiction
with $M(\mu_1)\subset M(s_\alpha\cdot\lambda)$, similarly as in the proof of \ref{1dirgvmeven}. 
The only other possibility is that some $2j-1$ will be replaced
with $2j+1$, i.e. $\mu_1+\delta=\half[2k-1,\ldots,2j+1,2j+1,2j-3,\ldots,3,-1|2n-1,\ldots,2j+3,2j-1,2j-1,\ldots,3,-1]$.
This weight should be strictly smaller than $\lambda$, but again, this is a contradiction, because the weight difference
$\mu_1-\lambda=[0,\ldots,0,1,\mathrm{something}]$ cannot be expressed as a sum of negative roots (negative roots are just
$[\ldots,-1,\ldots,1\ldots]$ and $[\ldots,-1,\ldots,-1\ldots]$).
\end{proof}

A similar theorem holds for the odd case: for $\lieg=B_{k+n}$, $\Sigma=\{\alpha_k\}$, we can choose 
$\lambda=[-n,\ldots,-n|1/2,\ldots,1/2]$ and $\mu=[-n,\ldots,-n-1|1/2\ldots,1/2]$. To prove the
existence of a nonzero \gvmhom, it suffices to use Theorem \ref{genlepowski} in this case (the proof
is absolutely analogous to the proof of lemma \ref{1dirgvmodd} and will be omitted).

\subsection{Description of the operator}
\label{morediracs}
Let us choose a real form $\lieg=\so(n+k,k)$ of the Lie algebra from the previous section generating the simply connected Lie group 
$G=\Spin(n+k,k)$ that
fixes the inner product $$\sum_{i=1}^k x_i x_{n+2k+1-i}+\sum_{j=1}^{n} x_{k+j}^2$$ 
(it has signature $(n+k,k)$).
Let $P$ be the parabolic subgroup and $\liep$ its Lie algebra, so that the complexification $(\lieg^c, \liep^c)$
is isomorphic to the pair from the last section ($\lieg=D_{k+n/2}$ in case $n$ is even and $B_{k+(n-1)/2}$, if $n$ is odd).
The reductive part is $\lieg_0=\sl(k,\R)\oplus\so(n)\oplus \R E$ and, as a $\lieg_0$-module,
$\lieg_{-1}\simeq ((\R^k)^*\otimes \R^{n})$,
the product of dual resp. defining representations of $\sl(k,\R)$, resp. $\so(n)$.
The $\lieg_{-2}$ component is commutative.
We will consider the case $n$ is even, the odd case is analogous.
Let $\mu, \lambda$ be weights like before and consider $\V_\mu$ and $\V_\lambda$ to be complex representations of
the real Lie algebra $\liep$ (given by restriction of the complex representations of $\liep^c$ to the real form)
and of the Lie group $P$. 
From the expression of $\mu, \lambda$ in $(\ref{muweight})$ and $(\ref{lambdaweight})$
we see that, as
a $\lieg_0^{ss}$-module, $\V_\mu\simeq {\C^k}^* \otimes S^-$ and
$\V_\lambda=\C\otimes S^+$ where $\C^k$ resp. $\C$ are the defining resp. trivial representation of $\sl(k,\R)$.

We know from the isomorphism $(\ref{realvermamodule})$ and Theorem $\ref{kdirgvmeven}$ that there is a nonzero
homomorphism of generalized Verma modules \gvmhom in this case as well.

The corresponding differential operator acts between sections of dual representation:
$$D: \Gamma(G\times_P (\C\otimes S^-))\to \Gamma(G\times_P (\C^k\otimes S^+)),$$
resp.
$$D: \Gamma(G\times_P (\C\otimes S^+))\to \Gamma(G\times_P (\C^k\otimes S^-)),$$
where we identified $(S^\pm)^*\simeq S^\mp$, resp. $(S^\pm)^*\simeq S^\pm$, depending on the parity of $n/2$ 
(see section \ref{diracintro}).  We will assume that $\V_\lambda^*\simeq S^-$, the other case is similar.

Assume that $s$ is a section of $G\times_P (\C\otimes S^-)$ and $f$ the associated $\C\otimes S^-\simeq S^-$-valued 
function on $\lieg_-$, defined by $(\ref{identification})$.
The coordinates on $\lieg_{-1}$ can be chosen to be $y_{11},\ldots,y_{1n},\ldots,y_{k1},\ldots,y_{kn}$ and on 
$\lieg_{-2}$ $y_1,\ldots, y_l$. To the section $Ds$ we assign a function 
$Df: \lieg_- \to \C^k\otimes S^+$ which can be naturally identified with $k$ $S^+$-valued function 
$D_1(f),\ldots,D_k(f)$.

Assume that $f$ is constant in the $\lieg_{-2}$ variables $y_1,\ldots,y_l$, so, it can be considered as a function of 
$y_{i,j}$ only.
The same argument as in the proof of Theorem \ref{whydirakisdirak} shows that the action of $g\in G_0$ on
$i(y)=\exp(y)P$ for $y\in\lieg_{-1}$ is $i(y)\mapsto i(\Ad(g)y)$. We know that $\lieg_{-1}\simeq (\R^k)^*\otimes \R^n$
as the adjoint representation of $G_0$ and we easily derive that the action of $g\in\Spin(n)$ on $f$ is
$(g\cdot f)(y_1,\ldots,y_k)=g(f(g^{-1}\cdot y_1,\ldots,g^{-1}\cdot y_k))$ and the action of $g\in\Sl(k)$ is
$(g\cdot f)(y_1,\ldots,y_k)=f((y_1,\ldots,y_k)g)$ ($\Sl(k)$ ``mixes'' the variables and the dual defining action of $g^{-1}$
on $(y_1,\ldots, y_k)$ is $(y_1,\ldots,y_k)g$). 

\begin{lemma}
The function $D_i(f)(y_1,\ldots,y_k)$ depends only on $$f(y_1,\ldots,y_{i-1},y,y_{i+1},\ldots,y_k), \quad y\in\R^{n}$$
\end{lemma}
\begin{proof}
Due to the fact that $D$ is a first order differential operators, the value $$D_i(f)(y_1,\ldots,y_k)$$ depends only on
$f(y_1,\ldots,y_k)$ and the derivatives $\frac{\partial}{\partial y_{ij}}f(y_1,\ldots,y_k)$, $j=1,\ldots,n$.
It suffice to prove that $D_1(f)(y_1,y_2)$ does not depend on $\frac{\partial}{\partial y_{2j}} f$ 
for $k=2$, the general case follows easily. 
Let us choose 
$g=\left(
\begin{tabular}{cc}
$K$ & $0$ \\
$0$ & $K^{-1}$
\end{tabular}\right)\in\Sl(2,\R)
$ for some $K\neq 1$.
The invariance of $D$ implies 
\begin{equation}
\label{invar}
(g\cdot (Df))(y_1,y_2)=(D(g\cdot f))(y_1,y_2)
\end{equation} 
The function $(g\cdot f)(y_1, y_2)=f((y_1, y_2)g)=f(Ky_1, K^{-1}y_2)$. The left hand side of $(\ref{invar})$
is 
$$
\left(
\begin{tabular}{cc}
$K$ & $0$ \\
$0$ & $K^{-1}$
\end{tabular}\right)
{D_1 f(Ky_1, K^{-1}y_2)\choose D_2 f(Ky_1, K^{-1}y_2)}={K (D_1 f)(Ky_1, K^{-1}y_2)\choose K^{-1}(D_2 f)(Ky_1, K^{-1}y_2)}
$$
Because $D$ is linear and first order, $D_1 f(Ky_1, K^{-1}y_2)$ depends linearly on 
$f, \frac{\partial f}{\partial y_1}$ and
$\frac{\partial f}{\partial y_2}$:
$$D_1 f=A f + B \frac{\partial f}{\partial y_1} + C \frac{\partial f}{\partial y_2}$$
in the point $(Ky_1, K^{-1}y_2)$. Comparing the left and right side of the first component of $(\ref{invar})$, we obtain
that 
$$K (A f+ B \frac{\partial f}{\partial y_1} + C \frac{\partial f}{\partial y_2})=f+ K A \frac{\partial f}{\partial y_1}+
K^{-1} C\frac{\partial f}{\partial y_2}$$
This holds for any $f$ and we see that $A=C=0$ and $D_1 f(Ky_1, K^{-1}y_2)$ depends only on 
$\frac{\partial f}{\partial y_1}$.
\end{proof}

\begin{theorem}
The operator $D$, considered as acting on functions on $\lieg_{-1}$ (i.e. on functions being constant on $\lieg_{-2}$) has
the following form:
$$D_i(f)=\sum_j e_j\cdot\partial_{i,j} f$$
Therefore, we call $D$ the ``Dirac operator in $k$ variables''.
\end{theorem}

\begin{proof}
Again, we can restrict to $k=2$ for simplicity. 
From the previous lemma it follows that $D_1 (f)(y_1,y_2)=D_1 (\tilde{f})(y_1,y_2)$ where $\tilde{f}(y_1',y_2'):=f(y_1',y_2)$
does not depend on the second variable and can be represented by a function $\tilde{\tilde{f}}$ of one (vector) variable $y_1'$: 
$\tilde{\tilde{f}}(y_1'):=f(y_1',y_2)$.
The group $\Spin(n)\subset G_0$ acts on $\tilde{\tilde{f}}$ by 
$$(g\cdot \tilde{\tilde{f}}(y_1'))=(g\tilde{f})(y_1',y_2)=g(\tilde{f}(g^{-1}y_1',g^{-1}y_2))=g(\tilde{\tilde{f}}(g^{-1}y_1'))$$
where $\Spin(n)$ acts on $\R^{n}$ by the fundamental representation. From the previous lemma 
$D_1\tilde{f}(y_1', y_2')=D_1\tilde{f}(y_1',y_2)$
and so, it can be represented by a function $(D_1 \tilde{\tilde{f}})(y_1')$ of one variable.
So, $D_1$ can be considered as  acting on functions $\tilde{\tilde{f}}$ of one variable and 
the invariance of $D$ with respect to $g\in\Spin(n)\subset G_0\subset G$ 
$$g ((D_i f)(g^{-1}y_1, g^{-1}y_2))=D_i (g\cdot f)(y_1, y_2),\quad i=1,2$$
implies that $D_1$, acting on spinor valued functions in one variable $y_1'$, is $\Spin(n)$-invariant as well.

Let $y\in\lieg_{-1}$ be a
vector with entries in the first column, $y=(\tilde{y}_1,0)$. 
The action of $g=\exp(y)$ on $s$ is
\begin{eqnarray*}
&&\mmmm (\exp(y) s)(\exp{(y_1,y_2)}P)=\exp(y)(s(\exp(-y)\exp{(y_1,y_2)}P))=\\
&&\mmmm \exp(y)[\exp((-\tilde{y}_1,0)+(y_1,y_2)+y_3)P, f((-\tilde{y}_1+y_1, y_2)+ y_3)]_P
\end{eqnarray*}
where $y_3\in\lieg_{-2}$. Now we use the assumption that $f$ does not depend on $\lieg_{-2}$ and 
$f((-\tilde{y}_1+y_1, y_2)+y_3)=f(-\tilde{y}_1+y_1, y_2)$. The last equation is equal to
$[\exp(y_1, y_2)P, f(-\tilde{y_1}+y_1, y_2)]_P$ and if we represent the section $s$ by a function $f$,
the action of $g$ is just 
$(g\cdot f)(y_1,y_2)=f(-\tilde{y}_1+y_1,y_2)$. But this means that $D_1$, acting on a function of one variable $\tilde{\tilde{f}}$, 
is invariant with respect to translations and $\Spin(n)$ group as well, so similarly as 
in the proof of Theorem \ref{whydirakisdirak}, we obtain that 
$$D_1(f)(y_1, y_2)=D_1(\tilde{f})(y_1,y_2)=\sum_j e_j \partial_{1,j} \tilde{f}(y_1,y_2)=\sum_j e_j \partial_{1,j} f(y_1, y_2)$$ 
Similarly for $D_2$.
\end{proof}

Note that the operator $D$ acts between sections of dual representations $\V_\lambda^*, \V_\mu^*$,
with highest weights 
\begin{eqnarray*}
&&\mmmm \V_\lambda^*\simeq \dynkin\whiteroot{0}\link\ldots\link\whiteroot{0}\link\noroot{\frac{n}{2}-1}\link\whiteroot{0}\link\ldots\link\whiteroot{0}\whiterootdownright{1}\whiterootupright{0} \enddynkin\\
&&\mmmm \V_\mu^*\simeq \dynkin\whiteroot{1}\link\whiteroot{0}\ldots\link\whiteroot{0}\link\noroot{\frac{n}{2}-1}\link\whiteroot{0}\link\ldots\link\whiteroot{0}\whiterootdownright{0}\whiterootupright{1} \enddynkin
\end{eqnarray*}
in the even case (the $0$ and $1$ on the last positions may be reversed, depending on the parity of $n/2$) and
\begin{eqnarray*}
&&\mmmm \V_\lambda^*\simeq \dynkin\whiteroot{0}\link\ldots\link\whiteroot{0}\link\noroot{\frac{n}{2}-1}\link\whiteroot{0}\link\ldots\link\whiteroot{0}\llink>\whiteroot{1}\enddynkin\\
&&\mmmm \V_\mu^*=\dynkin\whiteroot{1}\link\whiteroot{0}\link\ldots\link\whiteroot{0}\link\noroot{\frac{n}{2}-1}\link\whiteroot{0}\link\ldots\link\whiteroot{0}\llink>\whiteroot{1}\enddynkin
\end{eqnarray*}
in the odd case (note that the weight coordinate over the crossed node does not change here).

\newpage
\section{Singular orbit corresponding to $k$ Dirac operators}

\subsection{Even dimension, $k\leq n$}\label{evenorbits}
Let $\lieg=D_{k+n}=\so(2(n+k),\C)$, $\liep$ its parabolic subalgebra
corresponding to $$\dynkin\whiteroot{}\link\ldots\link\whiteroot{}\link\noroot{}\link\whiteroot{}\link\ldots\link\whiteroot{}\whiterootdownright{}\whiterootupright{} \enddynkin$$
where the $k$-th node is crossed ($\Sigma=\{\alpha_k\}$). We have seen in chapter \ref{diropinpargeo} that there exists a nonzero homomorphism
of generalized Verma modules $M_\liep(\mu)\to M_\liep(\lambda)$ where $\mu+\delta=\frac{1}{2}[\ldots,3,-1|\ldots,5,3,-1]$ and
$\lambda+\delta=\frac{1}{2}[\ldots,3,1|\ldots,5,3,1]$. 
We will investigate now the singular orbit of $\lambda$ and its
singular Hasse graph in case $k\leq n$.

\begin{lemma}
\label{2dirgvmeven}
Let $k=2$, $n\geq 2$. Then the singular Hasse graph for $(\lieg,\liep,\lambda)$ is
\begin{center}
\includegraphics{bgg.6}
\end{center}
\end{lemma}
\begin{proof}
The weights $\mu+\delta, \lambda+\delta$ are on the Weyl orbit of the $\lieg$-dominant weight 
$\tilde\lambda+\delta=\half[2n-1,\ldots,7,5,3,3,1,1]$ that is on the wall of the fundamental Weyl chamber,
so we see that they have singular character.
If $\nu$ is on the affine Weyl orbit of $\lambda$ and $\liep$-dominant, it means that
$\nu+\delta$ is on the Weyl orbit of $\lambda+\delta$ and strictly $\liep$-dominant.
In coordinates, $\nu+\delta$ consists of sign-permutation of $\half \{\ldots,7,5,3,3,1,1\}$.
A weight $[a_1,a_2|b_1,\ldots,b_n]$ is strictly $\liep$-dominant, iff $a_1>a_2$ and $b_1>\ldots,b_{n-1}>|b_n|$.
Therefore, $\nu+\delta=\half [a_1,a_2|\ldots,7,5,3,\pm 1]$ where $(a_1,a_2)$ is some decreasing 
sign-permutation of $(3,1)$
and the number of negative coordinates is even. There are only $4$ such possibilities: 
$$\half[3,1|\ldots,3,1],\,
\half[3,-1]\ldots,3,-1],\, \half[1,-3|\ldots,3,-1],\,\half[-1,-3|\ldots,3,1].$$ 
%
%

To show the existence of $\half[3,1|\ldots,3,1]-\delta\to\half[3,-1|\ldots,3,-1]-\delta$, consider $w\in W^p$
taking $\delta=[\ldots,2,1,0]$ to $w\delta=[2,0|\ldots,4,3,1]$ and $w'$ taking $\delta$ to $w'\delta=[2,-1|\ldots,4,3,0]$.
It is easy to see that $w(\tilde\lambda+\delta)=\half[3,1|\ldots,3,1]$,  $w'(\tilde\lambda+\delta)=\half[3,-1|\ldots,3,-1]$
and $w'=s_\gamma w$ for $\gamma=[0,1,0,\ldots,0,1]$. We know from lemma \ref{parabolicpath} that
there exists $w_j\in W^p$ such that $w\to w_1\to\ldots w_j=w'$.  
Further, if $E$ is the grading element, we see that 
$(w'\delta-w\delta)(E)=(2+0)-(2-1)=1$ and it follows from lemma \ref{grading2length} that $l(w')=l(w)+1$. So,
$w\to w'$ in the Hasse graph and $\half[3,1|\ldots,3,1]-\delta\to \half[3,-1|\ldots,3,-1]-\delta$ in the singular Hasse
graph.

To show the second arrow, let
$w\in W^p$ be the element that takes $\delta=[\ldots,3,2,1,0]$ to $[2,-1|\ldots,4,3,0]$ 
and $w'$ takes $\delta$ to $[1,-2|\ldots,4,3,0]$. We see that $w'=s_\gamma w$, where $\gamma$ is the 
positive root $[1,1,0,\ldots,0]$. 
In this case, $w\delta(E)=1$ and $w'\delta (E)=-1$, so the difference is $2$. But the grading element evaluation 
$w''\delta(E)$ cannot be zero for any $w''\in W^p$ in this case, because it is the sum of $2$ different integers. Therefore, 
$w\to w'$ in the parabolic Hasse graph 
and there is an arrow $[3,-1|\ldots]-\delta\to [1,-3|\ldots]-\delta$ 
in the singular Hasse graph for $(\lieg,\liep,\lambda)$.

Similarly, one can show that there is an arrow $\half[1,-3|\ldots,3,-1]-\delta\to \half[-1,-3|\ldots,3,1]-\delta$
by choosing $w$ taking $\delta$ to $[0,-2|\ldots,4,3,-1]$ and $w'$ taking $\delta$ to $[-1,-2|\ldots,4,3,0]$.
\end{proof}

We know that for a generalized Verma module homomorphism $M_\liep(\mu)\to M_\liep(\lambda)$ the difference
$(\lambda-\mu)(E)$ is the order of the dual differential operator, if it is one or two. 
Conjecture \ref{hypo} implies that
each arrow in the singular Hasse graph corresponds to a nonzero invariant differential operator between
sections of bundles associated to dual representations. In the following pictures, we will draw in the 
singular Hasse diagrams $q$ lines between $\lambda$ and $\mu$, if $(\lambda-\mu)(E)=q\in\{1,2\}$, i.e. 
the dual operator, if it exists, corresponds to the operator of order $q$ (Theorem \ref{grading2degree}). In particular,
the middle operator in our case is of second order, so we can draw the singular Hasse graph
\begin{center}
\includegraphics{bgg.9}
\end{center}
in this case.

\begin{theorem}
\label{kdirgvmevenorb}
Let $(\lieg,\liep,\lambda)$ be like at the beginning of this section, assume $k\geq 3$, $n\geq k$ and let  
$S_{k,n}$ be the singular Hasse graph associated to it. Then $S_{k,n}$ contains
two disjoint subgraphs $S^1$ and $S^2$, both isomorphic to $S_{k-1,n}$.

It follows that $S^1$ contains $S^{1,1}$ and $S^{1,2}$, both copies of $S_{k-2,n}$.
Similarly,
$S^2$ contains $S^{2,1}$ and $S^{2,2}$. Let $\phi_1, (\phi_2): S_{k-2,n}\to S^{1,2}\,(S^{2,1})$
be the isomorphisms, respectively. Then each element $\phi_1(\mu)$ of $S^{1,2}$ is connected to the
corresponding element $\phi_2(\mu)$ in $S^{2,1}$ by an arrow, that corresponds to a 
second order differential operator.

This describes all the singular Hasse graph of $(\lieg,\liep,\lambda)$.

Graphically, the $S_{k,n}$ has the following fractal-shape:

\begin{center}
\includegraphics{singorb.1}
\end{center}

\end{theorem}

\begin{definition}
\label{Sk}
From the lemma \ref{2dirgvmeven} and Theorem \ref{kdirgvmevenorb} it follows in case $k\leq n$, the shape of the graph $S_{k,n}$ 
depends only on $k$. Therefore, we can denote this graph by $S_k$. 
We also define $S_0$ to be a one-point graph.
\end{definition}

These pictures show $S_k$ for $k=3, 4$ (the arrows goes from up to down and from right to left):
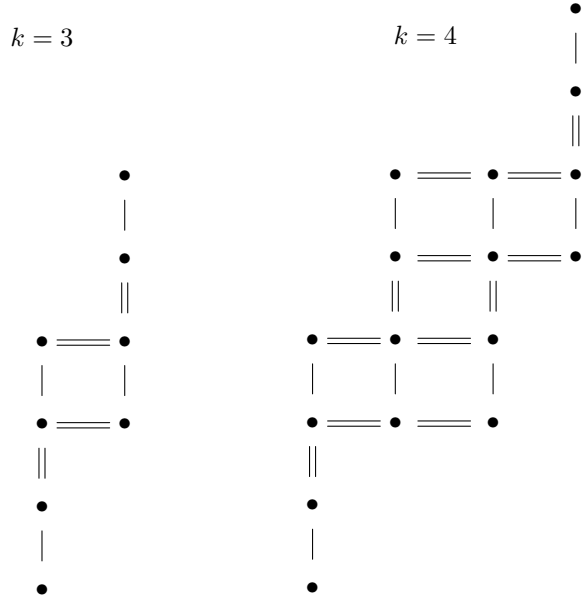
\begin{center}
\includegraphics{bgg.7}
\hspace{2cm}
\includegraphics{bgg.8}
\end{center}

\begin{proof}
Let us denote $\lieg_{k,n}=\so(2(n+k))$, $\liep_{k,n}$ the parabolic subalgebra from the beginning of this section,
and similarly,
$\delta_{k,n}=[n+k-1,\ldots,1,0]$,
$$\tilde\lambda_{k,n}+\delta_{k,n}=\half[2n-1,2n-3,\ldots,2k-1,2k-1,2k-3,2k-3,\ldots,3,3,1,1]$$
elements of $\lieh_{k,n}^*$, the $(k+n)$-dimensional Cartan subalgebra. 

We find all the weights $\nu$ so that $\nu+\delta$ is on the Weyl orbit of 
$\half[2k-1,\ldots,3,1|2n-1,\ldots,3,1]$ and strictly $\liep$-dominant. The condition
$k\leq n$ implies that $\nu+\delta=\half[\ldots|2n-1,\ldots,5,3,\pm 1]$ and on the first $k$
positions there is a strictly decreasing  sign-permutation of $\half\{2k-1,\ldots,5,3,1\}$.
Let $S^1$ be the set of weights of the type $\half[2k-1,\ldots|2n-1,\ldots,\pm 1]$ and 
$S^2$ the set of the weights $\half[\ldots,-(2k-1)|2n-1,\ldots,\pm 1]$. Clearly, $S_{k,n}$
is a disjoint union of $S^1$ and $S^2$.

We define a map 
$i: S_{k-1,n}\to S^1$, 
$$\half[a_1,\ldots,a_{k-1}|\ldots,3,\pm 1]-\delta_{k-1,n}\mapsto 
\half[2k-1,a_1,\ldots,a_{k-1}|\ldots,3,\pm 1]-\delta_{k,n}$$ We will show that it is a graph isomorphism. 
First, note that it is a bijection on the set of vertices
because any decreasing sign-permutation  $\half(a_1,\ldots,a_{k-1})$ of $\half\{2k-3,\ldots,3,1\}$ 
is strictly decreasing if and only if $\half(2k-1,a_1,\ldots,a_{k-1})$ is strictly
decreasing.

Suppose that there is an arrow $$\half[a_1,\ldots,a_{k-1}|b_1,\ldots,b_n]-\delta \to \half[a_1',\ldots,a_{k-1}'|b_1',\ldots,b_n']$$
in the singular Hasse graph for $(\lieg_{k-1,n}, \liep_{k-1,n}, \tilde\lambda_{k-1,n})$. 
The first weight, $\half[a|b]=w(\tilde\lambda_{k-1,n}+\delta_{k-1,n})$ and the second weight
$\half[a'|b']=w'(\tilde\lambda_{k-1,n}+\delta_{k-1,n})$ for some $w\to w'$. Let $w\delta=[x_1,\ldots,x_{k-1}|y_1,\ldots,y_n]$ and
$w'\delta=[x_1',\ldots,x_{k-1}'|y_1',\ldots,y_n']$. Note that all $x_i$ and $x_j'$ are strictly smaller then $2(k-1)$, because
$$w(\tilde\lambda_{k-1}+\delta_{k-1})=w\half[2(n+k)-3,\ldots,2k+1,2k-1,2k-3,2k-3,\ldots,3,3,1,1]$$ contains a number 
smaller or equal than $2k-3$ on the first position (which is at a position less than $2(k-1)$ from the right).

Let $i(w), i(w')$ be elements of $W_{k,n}^p$ satisfying $i(w)\delta_{k,n}=[2k-1,x_1,\ldots,x_{k-1}|\ldots]$ and 
$i(w')\delta_{k,n}=[2k-1,x_1',\ldots,x_{k-1}'|\ldots]$. The dots $\ldots$ don't mean $y_j$'s but something uniquely determined by
the first $k$ coordinates (the first $k$ coordinates
determine uniquely the last $n$ coordinates on the $W^p$-orbit of $\delta_{k,n}$).
We will show that there is an arrow $i(w)\to i(w')$ in $W^p$. Clearly, if $w'=s_\gamma w$, then $i(w')=s_{\gamma'} i(w)$ for 
the root $\gamma'=[0,\gamma]$. 

Further, we claim that the grading element evaluation $(w\delta-w'\delta)(E_{k-1,n})$ is $1$ or $2$. We will prove this by induction:
for the case $k=2$ it holds, as we saw in the proof of lemma \ref{2dirgvmeven} and we will show by induction that it holds in
general.

Because $i(w)\delta_{k,n}(E_{k,n})=(2k-1)+w\delta_{k-1,n}(E_{k-1,n})$, we see that 
$$(i(w)\delta_{k,n}-i(w')\delta_{k,n})(E_{k,n})$$
is $1$ or $2$ in this case as well. We see that $i(w)\leq i(w')$ and it follows from lemma \ref{grading2length} that 
$l(i(w'))-l(i(w))\leq 2$. But the length difference cannot be $2$, because $i(w)$ and $i(w')$ are connected by a root 
reflection, so the length difference must be an odd number. Therefore, it must be one and $i(w)\to i(w')$ in $W^p$.
It is easy to see that 
$$(i(w))(\tilde\lambda_{k,n}+\delta_{k,n})-\delta_{k,n}=i(w (\tilde\lambda_{k-1,n}+\delta_{k-1,n})-\delta)$$
$$(i(w'))(\tilde\lambda_{k,n}+\delta_{k,n})-\delta_{k,n}=i(w' (\tilde\lambda_{k-1,n}+\delta_{k-1,n})-\delta)$$
so there is an arrow $\half[2k-1,a|b]-\delta \to \half[2k-1,a'|b']-\delta$
in the singular Hasse graph and the map $i: S_{k-1,n}\to S_{k,n}$ preserves the arrows.

On the other hand, let $\half[2k-1,a_1,\ldots,a_{k-1}|b_1,\ldots,b_n]-\delta_{k,n}\to 
\half[2k-1,a_1',\ldots,a_{k-1}'|b_1',\ldots,b_n']-\delta_{k,n}$ 
be an arrow
in $S^1$. The weights are represented by Weyl group elements $w,w'$ so that $w\to w'$
and $\half[2k-1,a|b]=w(\tilde\lambda_{k,n}+\delta_{k,n})$ and $\half[2k-1,a'|b']=w'(\tilde\lambda_{k,n}+\delta_{k,n})$.
It follows that $w\delta=[x_1,\ldots,x_k|y_1,\ldots,y_n]$, where $x_1$ is either $2k-1$ or $2k-2$
and $x_2,\ldots,x_k$ are smaller. Similarly, $w'\delta=[x_1',\ldots,x_k'|y_1',\ldots,y_n']$ and 
$x_1=x_1'$ (because they are connected by a root reflection that does not fix $w(\tilde\lambda+\delta)$).
We define $\iota(w)=[x_2,\ldots,x_k|\ldots]$ and $\iota(w')=[x_2',\ldots,x_k'|\ldots]$ to be elements in
$W_{k-1,n}^p$, where the last $n$ coordinates are again uniquely determined. It is easy to check that
$\iota(w)\to \iota(w')$ in $W_{k-1,n}^p$ and that 
$\iota(w)(\tilde\lambda_{k-1,n}+\delta_{k-1,n})=[a_2,\ldots,a_k|b_1,\ldots,b_n]$ and 
$\iota(w')(\tilde\lambda_{k-1,n}+\delta_{k-1,n})=[a_2',\ldots,a_k'|b_1',\ldots,b_n']$.
So, there is an arrow $\half[a|b]-\delta\to \half[a'|b']-\delta$ in $S_{k-1,n}$ 
if and only if there is an arrow $\half[2k-1,a|b]-\delta\to \half[2k-1,a'|b']-\delta$
in $S^1$.

Similarly, we define the map $j:S_{k-1,n}\to S^2$ by
$\half[a_1,\ldots,a_{k-1}|b_1,\ldots,b_n]-\delta_{k-1,n}\mapsto \half[a_1,\ldots,a_{k-1},-(2k-1)|b_1,\ldots,b_n]-\delta_{k,n}$
and we can check that it maps arrows to arrows.

%

It remains to determine, which elements of $S^1$ are connected with arrows to elements of $S^2$. If a weight
$\half[2k-1,\ldots|\ldots]$ is connected by some root reflection $s_\gamma$ to $\half[\ldots,-(2k-1)|\ldots]$ 
so the only possible root is $\gamma=[1,0,\ldots,0,1|0,\ldots,0]$. Simple combinatorics implies that
the only possibility for the root reflection to exist is 
$\mu+\delta=\half[2k-3,a_2,\ldots,a_{k-1},-(2k-1)|\ldots]= 
s_\gamma \half[2k-1,a_2,\ldots,a_{k-1},-(2k-3)|\ldots]=s_\gamma (\lambda+\delta)$. 
The first weight is in $S^{2,1}$ and the second in $S^{1,2}$. 
To show that there is an arrow $\lambda\to\mu$, consider $w\in W^p$ taking $\delta$ to
$[2k-2,a_2,\ldots,a_{k-1},-(2k-3)|2n-2, 2n-4,2n-6,\ldots,2k-1,2k-4,\ldots,4,2,0]$ and 
$w'\in W^p$ taking $\delta$ to 
$[2k-3,a_2,\ldots,a_{k-1},-(2k-2)|2n-2,2n-4,2n-6,\ldots,2k-1,2k-4,\ldots,4,2,0]$, 
where $(a_2,\ldots,a_{k-1})$ 
is decreasing so that $w(\tilde\lambda+\delta)=\lambda+\delta$ and 
$w'(\tilde\lambda+\delta)=\mu+\delta$. 
(For example, the weight $\half [5,1,-3|5,3,-1]$ is represented with 
$w\delta=[4,1,-3|5,2,0]$ and $\half[3,1,-5|5,3,-1]$ is represented with
$w'\delta=[3,1,-4|5,2,0]$)
Now $w'=s_\gamma w$ implies $w\leq w'$ or $w'\leq w$. But for $E$ the grading element, 
$w\delta(E)=2k-2+a_2+\ldots+a_{k-1}-(2k-3)=\sum a_j+1$ and $w'\delta(E)=\sum a_j-1$, 
so $w\leq w'$. The difference of
the grading element evaluation is $2$, so the length difference is at most $2$, but it must be odd, because
$w'=s_\gamma w$ and therefore
$w\to w'$ in $W^p$.
So, there is an arrow 
$w\cdot \tilde\lambda\to w'\cdot\tilde\lambda$ in the singular Hasse graph. 

We see that $(\lambda-\mu)(E)=2$ and this also proves the
assumption that all the arrows $\lambda\to\mu$ fulfill $(\lambda-\mu)(E)\leq 2$. If there exists
the homomorphism of the generalized Verma modules, the corresponding dual differential operator
is of order $2$, therefore we draw $2$ lines in the picture.

Because we saw in lemma \ref{2dirgvmeven} that the form of the graph $S_{2,n}$ does not depend on $n$, 
so $S_{k,n}\simeq S_{k,n'}$ for $n,n'\geq k$ as well.
\end{proof}

\bigskip

\begin{theorem}
\label{standard-nonstandard}
All the arrows $\mu\to\nu$ in the singular Hasse graph from the previous theorem such that $(\mu-\nu)(E)=1$ 
are in the BGG graph 
and there is a standard homomorphism $M_\liep(\nu)\to M_\liep(\mu)$. If $\lambda\to\mu$ in the previous theorem and $(\lambda-\mu)(E)=2$
($2$ lines in the picture), the standard homomorphism $M_\liep(\mu)\to M_\liep(\lambda)$ is zero. So, under the assumption of 
conjecture \ref{hypo}, the homomorphisms corresponding to second order operators are all nonstandard.
\end{theorem}

\begin{proof}
We will prove the first part of the theorem by induction. In case $k=2$ we have seen in Theorem \ref{kdirgvmeven} that
there exists a standard homomorphism $[3,1|\ldots,3,1]-\delta\to [3,-1|\ldots,3,-1]$. Similarly, we could show that there 
exist a standard homomorphism $[1,-3|\ldots,3,-1]-\delta\to[-1,-3|\ldots,3,1]-\delta$.

The BGG graph has the same set of vertices as the singular Hasse graph, so we can define the sets $S_{k,n}, S^1, S^2$
the same way as in Theorem \ref{kdirgvmevenorb}
and we define the maps $i: S_{k-1,n}\to S^1$ and $j: S_{k-1,n}\to S^2$ by
$$i:[a_1,\ldots,a_{k-1}|b_1,\ldots,b_n]-\delta_{k-1,n}\mapsto [2k-1,a_1,\ldots,a_{k-1}|b_1,\ldots,b_n]-\delta_{k,n}$$
$$j:[a_1,\ldots,a_{k-1}|b_1,\ldots,b_n]-\delta_{k-1,n}\mapsto [a_1,\ldots,a_{k-1},-(2k-1)|b_1,\ldots,b_n]-\delta_{k,n}$$

We will show that there exists a nonzero standard homomorphism of true Verma modules 
$M_{\lieb_{k-1,n}}(\alpha)\to M_{\lieb_{k-1,n}}(\beta)$ ($\alpha,\beta\in S_{k-1,n}$)
if and only if there exists a nonzero standard homomorphism of true Verma modules
$M_{\lieb_{k,n}}(i(\alpha))\to M_{\lieb_{k,n}}(i(\beta))$. 
To see this, note that it follows from Theorem \ref{truevermamodulesmap} that for integral $\alpha, \beta\in\lieh_{k-1,n}^*$
there is a nonzero homomorphisms 
of true Verma modules $M(\alpha)\to M(\beta)$ if and only if there exist root reflections $\gamma_i$ so that 
$$\beta+\delta\geq s_{\gamma_1} (\beta+\delta)\geq s_{\gamma_2}s_{\gamma_1}(\beta+\delta)\geq\ldots\geq \alpha+\delta$$ 
But this is exactly if 
$$i(\beta)+\delta\geq s_{\tilde{i}(\gamma_1)} (i(\beta)+\delta)\geq s_{\tilde{i}(\gamma_2)}s_{\tilde{i}(\gamma_1)} (i(\beta)+\delta)\geq\ldots\geq (i(\alpha)+\delta)$$
where $\tilde{i}: \lieh_{k-1,n}\hookrightarrow\lieh_{k,n}$ is defined by $[\gamma]\mapsto [0,\gamma]$, especially 
$\alpha_i\in\Delta_{k-1,n}\mapsto \alpha_{i+1}\in\Delta_{k,n}$. 
So, nonzero homomorphism $M(\alpha)\to M(\beta)$ implies
nonzero homomorphism $M(i(\alpha))\to M(i(\beta))$. 
On the other hand, if $M(i(\alpha))\to M(i(\beta))$ is nonzero, there exists $\tilde\gamma_j$, so that
$$i(\beta)+\delta\geq s_{\tilde\gamma_1} (i(\beta)+\delta)\geq s_{\tilde\gamma_2}s_{\tilde\gamma_1} (i(\beta)+\delta)\geq\ldots\geq i(\alpha)+\delta.$$ 
All the $\tilde\gamma_i$ fix the first coordinate $(2k-1)/2$, because the first coordinate in the root expression cannot increase,
if the weights are decreasing, but the first coordinate of $i(\alpha)+\delta$ and $i(\beta)+\delta$ are both equal to $(2k-1)/2$. 
Therefore, all the $\tilde\gamma_j$ have pre-images $\gamma_j$ in $\lieh_{k-1,n}$ so that $\tilde{i}(\gamma_j)=\tilde\gamma_j$.
It follows that
$$\beta+\delta\geq s_{\gamma_1} (\beta+\delta)\geq s_{\gamma_2}s_{\gamma_1} (\beta+\delta)\geq\ldots\geq (\alpha+\delta)$$ 
and there exists a nonzero map $M(\alpha)\to M(\beta)$. So, there exists a nonzero homomorphism $M(\alpha)\to M(\beta)$
if and only if there exists a nonzero homomorphism $M(i(\alpha))\to M(i(\beta))$.

We will prove the following statement: there exists a nonzero standard homomorphism 
$$M_{\liep_{k-1,n}}(\alpha)\to M_{\liep_{k-1,n}}(\beta)$$ for $\alpha,\beta\in S_{k-1,n}$
if and only if there exists a nonzero standard homomorphism
$$M_{\liep_{k,n}}(i(\alpha))\to M_{\liep_{k,n}}(i(\beta))$$
It follows from Theorem \ref{zeromap} that the standard homomorphism  
$$M_{\liep_{k-1,n}}(\alpha)\to M_{\liep_{k-1,n}}(\beta)$$ is zero
if and only if $M(\alpha)\subset M(s_{\alpha_j}\cdot\beta)$ for some parabolic simple root
$\alpha_j\neq \alpha_{k-1}$. If this is the case, then $M(i(\alpha))\subset M(s_{\alpha_{j+1}}\cdot i(\beta))$
follows from the previous paragraph and the map $M_\liep(i(\alpha))\to M_\liep(i(\beta))$ is zero as well. On the other hand,
if $M_\liep(i(\alpha))\to M_\liep(i(\beta))$ is zero, then $M(i(\alpha))\subset M(s_{\alpha_i}\cdot i(\beta))$
for some parabolic simple root $\alpha_i\neq \alpha_k$. If $i=1$, then $M(i(\alpha))\subset M(s_{\alpha_1}\cdot i(\beta))$ implies
$i(\alpha)+\delta\leq s_{\alpha_1}(i(\beta)+\delta)$. But $i(\alpha)+\delta$ contains $(2k-1)/2$ on the first position and 
$s_{\alpha_1}(i\beta)+\delta)$ contains
a number strictly smaller then $(2k-1)/2$ on the first position, what is a contradiction.
Therefore, $i>1$ and that implies $M(\alpha)\subset M(s_{\alpha_{i-1}}\cdot\beta)$,
so the map $M_\liep(\alpha)\to M_\liep(\beta)$ is zero as well. 

%
Similarly, we can show that the map $j: S_{k-1,n}\to S^2$ given by 
is a graph isomorphism. 

To complete the proof, we have to show that the standard homomorphism 
$$M_\liep(\half[2k-1,a_2,\ldots,a_{k-1},-(2k-3)|2n-1,\ldots,3,\pm 1]-\delta)\to $$
$$\to M_\liep(\half[2k-3,a_2,\ldots,a_{k-1},-(2k-1)|2n-1,\ldots,3,\pm 1]-\delta)$$
is zero (this are exactly the arrows represented by $2$ lines in the singular Hasse graph, see theorem
\ref{kdirgvmevenorb}).

In case $k=2$, we observe that
\begin{eqnarray*}
\half[3,-1|\ldots,3,-1]&&\mmmm \geq s_\alpha \half[3,-1|\ldots,3,-1]=\half[3,-1|\ldots,1,-3]\geq\\ 
&&\mmmm \geq \half[1,-1|\ldots,3,-3]\geq \half[1,-3|\ldots,3,-1]
\end{eqnarray*}
The first reflection is with respect to the parabolic simple root $\alpha=[0,\ldots,1,1]$, the second interchanges
$1$ and $3$, the reflection being with respect to $\gamma_1=[1,0,\ldots,0,-1,0]$ and the last interchanges
$-1$ and $-3$, the reflection being with respect to the root $\gamma_2=[0,1,0\ldots,0,-1]$. So, it follows
from Theorem \ref{truevermamodulesmap} that $M(\half[1,-3|\ldots]-\delta)\subset M(s_\alpha\cdot(\half[3,-1|\ldots,3,-1]-\delta))$
and we see from Theorem \ref{zeromap} that the standard homomorphism 
$M_\liep(\half[1,-3|\ldots]-\delta)\to M_\liep(\half[3,-1|\ldots]-\delta)$
is zero.

Similarly, for $k\geq 2$, we observe that

\begin{eqnarray*}
&& \mmmm     \half[2k-1,a_2,\ldots,a_{k-1},-(2k-3)|2n-1,\ldots,3,\pm 1]\geq \\
&& \mmmm\geq \half[2k-1,a_2,\ldots,a_{k-1},-(2k-3)|2n-1,\ldots,2k-3,2k-1,\ldots,3,\pm 1]\\
&& \mmmm\geq \half[2k-3,a_2,\ldots,a_{k-1},-(2k-3)|2n-1,\ldots,2k-1,2k-1,\ldots,3,\pm 1]\\
&& \mmmm\geq \half[2k-3,a_2,\ldots,a_{k-1},-(2k-1)|2n-1,\ldots,2k-1,2k-3,\ldots,3, \pm 1]
\end{eqnarray*}
The first reflection is with respect to the simple root $\alpha=[0,\ldots,0,1,-1,0\ldots]$ interchanging 
$2k-1$ and $2k-3$ on the $(n+2)$-nd and $(n+3)$-th position, 
the second reflection interchanges the $2k-1$ on the first position with the $2k-3$ on
the $(n+2)$-th position, and the last reflection sign-interchanges the $-(2k-3)$ on the $k$-th position
and the $2k-1$ on the $(n+2)$-nd position. It is easy to check that 
$$M(\half[2k-3,\ldots,-(2k-1)|\ldots]-\delta)\subset M(s_\alpha\cdot(\half[2k-1,\ldots,-(2k-3)|\ldots]-\delta))$$
and it follows from \ref{zeromap} that the standard homomorphism is zero.

\end{proof}

\begin{remark}
Note, however, that for $k>n$ the singular orbits is larger as the one described in the theorem. We will show this in more details 
in section \ref{strangeorbits}.
\end{remark}

\begin{remark}
If we choose the weight $\lambda+\delta=\half[\ldots,3,1|\ldots,5,3,-1]$ (corresponding to the other spinor representation), the singular 
Hasse graph associated to it has the same structure and there is no real difference in the proof.
\end{remark}


\subsection{Even dimension, $k>n$}
\label{strangeorbits} Let $\lieg_{k,n}=\so(2(n+k),\C))$, $\liep$ the parabolic subalgebra
corresponding to the $k$-th node crossed and let $\lambda$ be as in section \ref{evenorbits}.
We saw that the shape of the singular Hasse graph associated to $(\lieg,\liep,\lambda)$
does not depend on $n$ for $n\geq k$ and all the weights $\nu$ on the affine Weyl
orbit of $\lambda$ are of the form $\nu+\delta=\half[\ldots|2n-1,\ldots,5,3,\pm 1]$.
In case
$n<k$, there exists other $\liep$-dominant weights on the affine Weyl orbit of
$\lambda$ and the singular Hasse graph is larger. We will illustrate it on the simplest example:

\begin{example}
\label{exak=3n=2}
Let $k=3, n=2.$ Then the singular Hasse graph has to following form:
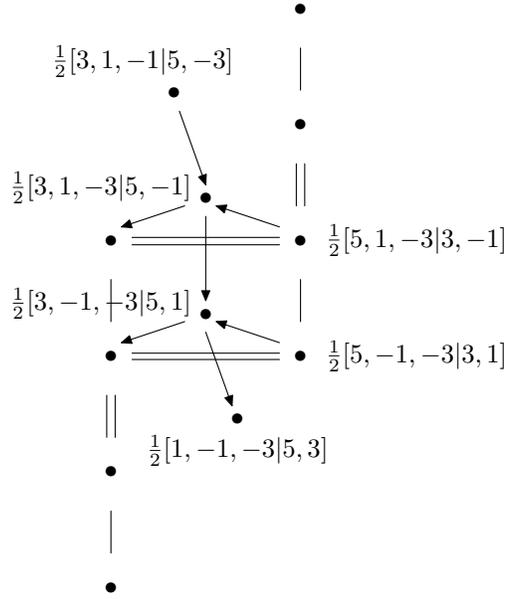
\begin{center}
\includegraphics{singorb.2}
\end{center}
\end{example}
In this diagram, we should subtract $\delta$ from each of the weights to obtain
the singular Hasse graph for $(\lieg,\liep,\lambda)$.

We see that there are four ``new'' strictly $\liep$-dominant weights 
of type $$\half[\ldots|5,\pm 1] \quad\mathrm{and}\quad \half[\ldots|5,\pm 3]$$ 
on the
Weyl orbit of $\lambda+\delta$. Each of the new arrows corresponds clearly
to a root reflection. To see that each of the new arrows is in the
singular Hasse graph, consider, for example, 
the arrow $(\half[5,1,-3|3,-1]-\delta) \to (\half[3,1,-3|5,-1]-\delta)$.
Let $\tilde\lambda+\delta=\half[5,3,3,1,1]$ 
be the dominant weight on the orbit of $\lambda+\delta$.
Choose $w, w'$  be elements of $W^p$ that take $\delta=[4,3,2,1,0]$ to 
$w\delta=[4,2,-1|3 -0]$ and $w'\delta=[3,2,-1|4, -0]$. Because $w'=s_\gamma w$
for $\gamma=[1,0,0,-1,0]$ and $(w\delta-w'\delta)(E)=1$ (E is the grading element),
it follow $w\to w'$ in $W^p$. The weights $w(\tilde\lambda+\delta)=\half[5,1,-3|3,-1]$
and $w'(\tilde\lambda+\delta)=\half[3,1,-3|5,-1]$. 
The other arrows in the above diagram can be shown in an analogous way.

Moreover, we claim that all these new arrows are in the BGG graph as well, representing standard
homomorphisms of generalized Verma modules. To show, for example, that there is a nonzero standard 
homomorphism
$$M_\liep(\half[3,1,-3|5,-1]-\delta)\to M_\liep(\half[5,1,-3|3,-1]-\delta)$$
assume that it is zero, then there is a chain of weights 
$$\half[3,1,-3|5,-1]=\mu_0\leq \mu_1\leq\ldots\leq \mu_j= s_\alpha \half[5,1,-3|3,-1]$$
for some parabolic simple root $\alpha\neq \alpha_3$, connected by root reflections.

Because $\mu_0$ is $\liep$-dominant and the only elements fixing $\mu_0(E)$ are from $W_p$,
it follows that $\mu_1(E)=\mu_0(E)+1$. If $\mu_1=\half[3,1,-1|\ldots]$, $\mu_1$ would be $\liep$-dominant
and to leave the $W_p$-orbit of this weight, one would have to increase the grading element evaluation
once more, which is not possible. If $\mu_1=\half[5,1,-3|\ldots]$, the only possibility is 
$\mu_1=\half[5,1,-3|3,-1]$ which is larger than any $s_\alpha \half[5,1,-3|3,-1]$. In any case, we get a 
contradiction.

Let us denote by $S_{k,n}$ the singular Hasse graph associated to $(\lieg,\liep,\lambda)$
in general.
We will show that it has similar structure as the picture above
if $n=k-1$. The weights of type $\half[\ldots|2n-1,\ldots,3,\pm 1]$
and the arrows between them are a copy of the graph $S_{k}$ from definition
\ref{Sk}. 
We can define the subgraphs $S^{1,1}$, $S^{1,2}$, $S^{2,1}$ and $S^{2,2}$
similarly as in Theorem \ref{kdirgvmevenorb}. For example, $S^{1,2}$ consists
of weights $\half[2k-1,\ldots,-(2k-3)|2n-1,\ldots,3,\pm 1]$. We denote
by $K$ the subgraph of $S_{k,n}$ consisting of weights not in $S^1$, $S^2$.

\begin{theorem}
Let $n=k-1$, $n\geq 2$.
The subgraph $K$ of $S_{k,n}$ contains a copy of $S_{k-2}$ that consists of 
weights $\half[2k-3,\ldots,-(2k-3)|2k-1,2k-5,2k-7,\ldots,3,\pm 1]$ (Remember
that $2k-3=2n-1$). We
denote this subgraph by $K_{k-2}$. The graph isomorphisms $\psi: S_{k-2}\to K_{k-2}$
is given by 
\begin{eqnarray*}
&&(\half[a_1,\ldots,a_{k-2}|\ldots,3,\pm 1]-\delta_{k-2,n}) \mapsto \\ 
&&\mapsto (\half[2k-3,a_1,\ldots,a_{k-2},-(2k-3)|2k-1,2k-5, \ldots,3,\pm 1]-\delta_{k,n})\\
\end{eqnarray*}
(\{$a_1,\ldots,a_{k-2}$\} is a decreasing sign-permutation of $\half\{(2k-5),\ldots,3,1\}$).

Let $\phi_1\, (\phi_2): S_{k-2}\to S^{1,2} \,\, (S^{2,1})$ be the isomorphism from
Theorem \ref{kdirgvmevenorb}, respectively. Then for each $\mu\in S_{k-2}$, there
is an arrow $\phi(\mu)\to \psi(\mu)$ and an arrow $\psi(\mu)\to \phi_2(\mu)$ in 
$S_{k,n}$. Graphically, it is described by the following picture:
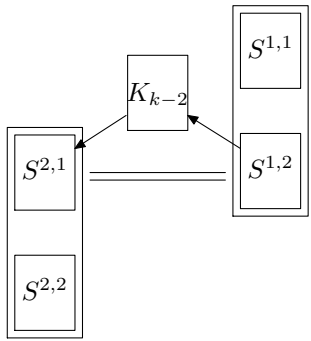
\begin{center}
\includegraphics{singorb.3}
\end{center}
The arrows connecting $K_{k-2}$ with $S$ correspond to operators of first order.

Again, we can divide $K_{k-2}$ into two parts $K_{k-2}^1$ and $K_{k-2}^2$, each
isomorphic to $S_{k-3}$ as a graph. Then there are two copies $K_{k-3}^{(1)}$ and
$K_{k-3}^{(2)}$ of $S_{k-3}$ in $K-K_{k-2}$ connected to $K_{k-2}$ in a natural way:
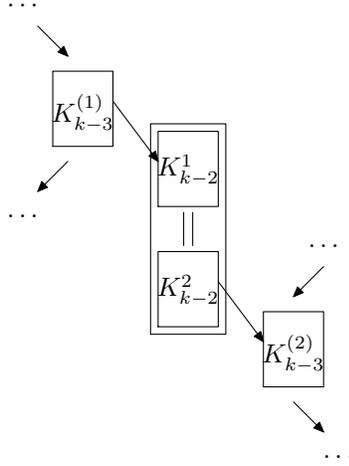
\begin{center}
\includegraphics{singorb.4}
\end{center}
The arrows connecting $K_{k-3}^{(i)}$ with $K_{k-2}^j$ correspond to first order differential operators. 
In a similar way, $K_{k-3}^{(i)}$ is connected with $K_{k-4}^{(i)(1)}$ and
$K_{k-4}^{(i)(2)}$, each being isomorphic to $S_{k-4}$.
So we can continue until we come to $S_{0}$, which is a one-point graph.
This describes all the arrows and vertices in the singular Hasse graph.

All the arrows $\mu\to\nu$ such that $(\mu-\nu)(E)=1$ are in the BGG graph as well and there
exists a nonzero standard homomorphism $M_\liep(\nu)\to M_\liep(\mu)$ in such case.
\end{theorem}

\begin{proof}
The technique of the proof is similar to the proofs of Theorem \ref{kdirgvmevenorb} and example \ref{exak=3n=2}.
We will just outline the basic steps.

First, note that if a weight $\nu$ on the affine Weyl orbit of $\lambda$ is $\liep$-dominant,
so $\nu+\delta$ is of the form $\half[\ldots|b_1,\ldots,b_n]$ where $(b_1,\ldots,b_n)$ 
is strictly decreasing. We can prove in the same way as in
Theorem \ref{kdirgvmevenorb} that the weights of type $(b_1,\ldots,b_n)=
(2n-1,\ldots,3,\pm 1)$ and arrows between them are a copy of $S_{k}$. Let us
denote by $K_{k-2}$ the set of weights of type $(b_1,\ldots,b_n)=(2n+1,2n-3,\ldots)$.
The map $\psi$ is a graph isomorphism $S_{k-2}\simeq K_{k-2}$. We will omit the 
subtraction of $\delta$ in the expression of the following weights.
An element of $K_{k-2}^1$
is of the form $\half[2k-3,2k-5,c_1,\ldots,c_{k-3},-(2k-3)|2k-1,2k-5,\ldots]$ and it is connected
with root reflection to $\half[2k-3,2k-5,c_1,\ldots,c_{k-3},-(2k-5)|2k-1,2k-3,2k-7,\ldots]$.
These are exactly elements of $K_{k-3}^{(1)}$. Similarly, an element of $K_{k-2}^2$
is of the form $\half[2k-3,c_1,\ldots,c_{k-3},-(2k-5),-(2k-3)|2k-1,2k-5,\ldots]$ and it is connected
with root reflection to $\half[2k-5,c_1,\ldots,c_{k-3},-(2k-5), -(2k-3)|2k-1,2k-3,2k-7,\ldots]$.
This are exactly elements of $K_{k-3}^{(2)}$. We can continue similarly and see that
elements of $K_{k-4}^{(i)(j)}$ are weights of the form 
$\half[\ldots|2k-1,2k-3,2k-5,2k-9,\ldots,3,\pm 1]$. Finally, $K_{0}^{(i_1)\ldots (i_{k-2})}$
are weights of the form $\half[\ldots|2k-1,2k-3,\ldots,7,5,\pm 3]$.
\end{proof}

For general $k>n$, the situation is even more complicated. 
%


\subsection{Odd dimension} 
Let $\lieg=B_{k+n}=\so(2(n+k)+1,\C)$, $\liep$ its parabolic subalgebra
corresponding to $$\dynkin\whiteroot{}\link\ldots\link\whiteroot{}\link\noroot{}\link\whiteroot{}\link\ldots\link\whiteroot{}\llink>\whiteroot{}\enddynkin$$
where the $k$-th node is crossed ($\Sigma=\{\alpha_k\}$). We have seen in chapter \ref{diropinpargeo} that there exists a nonzero homomorphism
of generalized Verma modules $M_\liep(\mu)\to M_\liep(\lambda)$ where $\mu+\delta=[\ldots,3/2,-1/2|\ldots,3,2,1]$ and
$\lambda+\delta=[\ldots,3/2,1/2|\ldots,3,2,1]$.
We will show that the BGG graph associated to $(\lieg,\liep,\lambda)$ is isomorphic to
the singular Hasse graph $S_{k,n}$ from Theorem \ref{kdirgvmevenorb}. On the other hand, the singular Hasse graph does
not contain all the arrows in this case. Moreover, the BGG graph does not depend on $n$ at all (for a fixed $k$).
We will start with the simplest case $k=2$.

\begin{lemma}
\label{2diraksodd}
Let $k=2$. Then there exist three nonzero generalized Verma module homomorphisms on the affine orbit of $\lambda$
described by the following diagram (the middle homomorphism corresponds to an operator of second order):
\begin{center}
\includegraphics{bgg.10}
\end{center}
\end{lemma}

\begin{proof}
The existence of the homomorphism $M_\liep(\mu)\to M_\liep(\lambda)$ was shown in chapter 
\ref{diropinpargeo}. 

We see that $[1/2,-3/2|\ldots,2,1]$ and $[3/2,-1/2|\ldots,2,1]$ are connected by root reflection $s_\gamma$,
where $\gamma=[1,1|0,\ldots,0]$ and $[3/2,-1/2|\ldots,2,1](H_\gamma)=3/2-1/2=1$ is a nonnegative integer.
It follows from Theorem \ref{truevermamodulesmap} that there is a true Verma module homomorphism 
$i: M(\nu)\to M(\mu)$. 
We will show that the standard map of the generalized Verma modules is nonzero.

The weight $\lambda+\delta=[3/2,1/2|\ldots,2,1]$ is on the Weyl orbit of the dominant (but non-integral) weight 
$\tilde\lambda+\delta=[\ldots,2,3/2,1,1/2]$. This weight is non-singular, so there exist unique elements 
$w, w'\in W$ so that $w\cdot \tilde\lambda=\mu$ and $w'\cdot\tilde\lambda=\nu$. In this case, $\delta=\half[\ldots,7,5,3,1]$
and it is easy to check that $w\delta=\half[5,-1|\ldots,11,9,7,3]$ and $w'\delta=\half[1,-5|\ldots,11,9,7,3]$. Because
$w'=s_\gamma w$ and the evaluation on the grading element is 
$(w\delta-w'\delta)(E)=(5/2-1/2)-(1/2-5/2)=4$, we see that $w\leq w'$ and either $w\to w'$ or $w\to w_1\to w_2 \to w'$
in the (Borel) Hasse graph (there must be an odd number of root reflections, if their composition is a root reflection and
it cannot be $5$, because the difference of the grading element evaluation is 4). If $w\to w'$, we can use Theorem 
\ref{genlepowski} and we are done. Assume that $w\to w_1\to w_2\to w'$ and that the standard map 
$M_\liep(\nu)\to M_\liep(\mu)$ is zero. Theorem \ref{zeromap} says that
\begin{equation}
\label{condition}
M(\nu)\subset M(s_\alpha\cdot\mu)
\end{equation}
for some parabolic simple root $\alpha\neq \alpha_k$. 
We know that for $\liep$-dominant $\mu$, $M(s_\alpha \cdot \mu)\subsetneq M(\mu)$.
The weight $s_\alpha(\mu+\delta)$ is one of the following types:
\begin{enumerate}
\item{$[-1/2,3/2|\ldots,3,2,1]$ if $\alpha=\alpha_1$}
\item{$[3/2,-1/2|n,n-1,\ldots,l-1,l,\ldots,2,1]$ is $\alpha=\alpha_i$, $1<i<n+k$} 
\item{$[3/2,-1/2|\ldots,3,2,-1]$ if $\alpha=\alpha_{n+k}$} 
\end{enumerate}
First we show that $\alpha\neq\alpha_1$. If $\alpha=\alpha_1$, (\ref{condition}) implies 
$\nu+\delta \leq s_{\alpha_1}(\mu+\delta)$, i.e. $[1/2,-3/2|\ldots,2,1]\leq [-1/2,3/2|\ldots,2,1]$. Subtracting the second
weight from the first we get $[1,-3|0,\ldots,0]$, which cannot be obtained as a sum of negative roots,
because none of them has a positive coefficient on the first position.

Now assume that $s_\alpha(\mu+\delta)$ is of type $(2)$. Because
$$
M(w'\cdot\tilde\lambda)=M(\nu)\subsetneq M(s_\alpha\cdot\mu)\subsetneq M(\mu)=M(w\cdot\tilde\lambda),
$$
$l(w')-l(w)=3$ and $\nu+\delta$ is not connected with $s_\alpha(\mu+\delta)$ with any root reflection, it follows from
Theorem \ref{truevermamodulesmap} that there must be $\beta_1, \beta_2$ so that
\begin{equation}
\label{vermasubsets}
M(\nu)\subsetneq M(s_{\beta_1}\cdot\nu)=M(s_{\beta_2} s_\alpha \cdot\mu)\subsetneq M(s_\alpha\cdot\mu).
\end{equation}

Note, that the weights are $s_\alpha (\mu+\delta)=[3/2,-1/2|\ldots,l-1,l,\ldots,2,1]$ and 
$\nu+\beta=s_{\beta_1} s_{\beta_2} s_\alpha(\mu+\delta)=[1/2,-3/2|\ldots,2,1]$. In coordinates, $s_{\beta_j}$ cannot be a (sign)-transposition
interchanging an integer and a half-integer, because of the conditions $s_\alpha (\mu+\delta)(H_{\beta_2})\in\N$ and 
$s_{\beta_2}s_\alpha(\mu+\delta)(H_{\beta_1})\in\N$.
So, exactly one of these reflections
interchanges $(3/2,-1/2)$ to $(1/2,-3/2)$ and the other one interchanges $(n,n-1,\ldots,l-1,l,\ldots,2,1)$ to 
$(n,n-1\ldots,l,l-1,\ldots,2,1)$. So either $s_{\beta_2}s_\alpha (\mu+\delta)=[1/2,-3/2|\ldots,l-1,l\ldots]$ or
$s_{\beta_2}s_\alpha (\mu+\delta)=[3/2,-1/2|\ldots,l,l-1,\ldots]$. In the first case, 
$s_{\beta_2} s_\alpha(\mu+\delta)=s_{\alpha}(\nu+\delta)<(\nu+\delta)$ ($\nu$ is $\liep$-dominant) which contradicts 
$(\ref{vermasubsets})$. In the second case, $s_{\beta_2}s_\alpha (\mu+\delta)=\mu+\delta>s_\alpha(\mu+\delta)$ 
which also contradicts  $(\ref{vermasubsets})$. 
So $s_{\alpha}(\mu+\delta)$ cannot be of type $(2)$.

Similarly, we can show that $s_\alpha (\mu+\delta)$ cannot be of type $(3)$. But this means that $(\ref{condition})$
does not hold and the standard map $M_\liep(\nu)\to M_\liep(\mu)$ is nonzero.

Finally, to show that there is a nonzero standard homomorphisms $M(\xi)\to M(\nu)$, note that 
$\nu+\delta=w(\tilde\lambda+\delta)$ and $\xi+\delta=w'(\tilde\lambda+\delta)$ where $w\delta=\half[1,-5|\ldots]$
and $w'\delta=\half[-1,-5|\ldots]$. Evaluating $w\delta$ and $w'\delta$ on the grading element, we see that
the difference is $1$ and $w\to w'$ (lemma \ref{grading2length}). Finally,
$\nu-\xi=[1,0,\ldots,0]$ is an integral multiple of a root, so $M(\xi)\subset M(\nu)$. 
Applying Theorem \ref{genlepowski}, we see that the standard homomorphism $M_\liep(\xi)\to M_\liep(\nu)$
is nonzero as well.
\end{proof}

Note, however, that in the case $k=n=2$, there are other strictly $\liep$-dominant weights on the orbit of $\lambda+\delta$. 
Namely, the weights 
\begin{equation*}
[2,1|3/2,1/2], [2,-1|3/2,1/2], [1,-2|3/2,1/2] \quad\and\quad [-1,-2|3/2,1/2]. 
\end{equation*}
Neither of them
is connected with a root reflection to any weight from the last theorem. It is easy to see that
they are connected with standard homomorphism analogous as the weight in the last theorem. 
So, in that case, the BGG graph
consists of $2$ connected parts,
each of them being a copy of $S_2$ (cf. Theorem \ref{kdirgvmevenorb}). 
However, the homomorphisms of generalized Verma modules appearing in the other component correspond to
differential operators that are all of second order (because, intechanging $1$ with a $-1$ decreases 
the grading element evaluation by $2$, so we can use Theorem \ref{grading2degree}).

\begin{remark}
The singular Hasse graph, however, does not contain the arrow $\mu\to\nu$.
\end{remark}

\begin{proof}
For $\tilde\lambda+\delta=[\ldots,3,2,3/2,1,1/2]$, there is a unique $w$ resp. $w'$ 
taking $\tilde\lambda+\delta$ to $\mu+\delta$ resp. $\nu+\delta$. 
As we saw in the proof of \ref{2diraksodd}, $w\delta=\half[5,-1|\ldots,7,3]$ and 
$w'\delta=\half[1,-5|\ldots,7,3]$. The length difference $l(w')-l(w)$ is not $1$, because
(identifying $w\leftrightarrow w\delta$ for $w\in W^p$)
\begin{eqnarray*}
\half[5,-1|\ldots,7,3]&& \mmmm\to\half[3,-1|\ldots,7,5]\to\half[1,-3|\ldots,7,5]\to \\
&& \mmmm\to\half[1,-5|\ldots,7,-3]
\end{eqnarray*}
\end{proof}
\begin{theorem}
\label{2diroddcomplex}
The BGG graph from lemma \ref{2diraksodd} is a complex.
\end{theorem}

\begin{proof}
We want to show that the standard homomorphism $M_\liep(\nu)\to M_\liep(\lambda)$ is zero. This can be see from
the following sequence of weights connected by reflections:
$$[1/2,-3/2|\ldots, 2,1]\leq \half[1/2,3/2|\ldots,2,1]=s_{\alpha_1} [3/2,1/2|\ldots,2,1]$$
Similarly, the standard homomorphism $M_\liep(\xi)\to M_\liep(\mu)$ is zero because
$$[-1/2,-3/2|\ldots,2,1]\leq [-1/2,3/2|\ldots,2,1]=s_{\alpha_1}[3/2,-1/2|\ldots,2,1]$$
\end{proof}

The situation in $k>2$ is analogous:

\begin{theorem}
\label{kdirgvmodd}
For any $k,n\geq 2$, $k\neq n$, the graph $S_k$ from Theorem \ref{kdirgvmevenorb} 
describes the BGG graph associated to $(\lieg,\liep,\lambda)$.
The double-arrows describe homomorphisms, whose corresponding invariant differential operators are of second order.

All homomorphism in the graph are standard.

However, the arrows corresponding to second order operators are not in the singular Hasse graph. 

In case $k=n$, the orbit has two connected components. One of them is a copy of $S_k$ and the other 
one is similar except that all of its homomorphisms correspond  to second order differential operators.
\end{theorem}

\begin{proof}
Let $\tilde\lambda+\delta=[\ldots,5/2,2,3/2,1,1/2]$ be the $\lieg$-dominant weight so that $\lambda+\delta$ 
is on its Weyl orbit. 
Recall that there are $n$ integers and $k$ half-integers in the expression of $\tilde\lambda+\delta$. Assume that $k\neq n$.

The condition on a weight $\nu+\delta=[a_1,\ldots,a_k|b_1,\ldots,b_n]$ to be strictly $\liep$-dominant and $\liep$-integral is $a_1>\ldots>a_k$,
$b_1>\ldots>b_n>0$, $a_i-a_j\in\Z$, $b_i-b_j\in\Z$ and the $b_i$'s are all
integers or all half-integers. Simple combinatorics implies that, if $\nu$ is on the affine orbit of $\lambda$ and
$k\neq n$, the only possibility is $\nu+\delta=[a_1,\ldots,a_k|n,n-1,\ldots,2,1]$, where $(a_1,\ldots,a_k)$
is some strictly decreasing sign-permutation of $((2k-1)/2,\ldots,3/2,1/2)$.

Let $S_{k,n}$ be the BGG graph for $(\lieg,\liep,\lambda)$.
In the same way as in Theorem \ref{kdirgvmevenorb} we can define $S^1$ to be 
the subgraph consisting of weights $[(2k-1)/2,\ldots|\ldots]$ and
$S^2$ the subgraph consisting of weights $[\ldots,-(2k-1)/2|\ldots]$. Clearly, the set of all vertices in the BGG graph 
is a disjoint union of vertices in $S^1$ and vertices in $S^2$. Similarly, we can define $S^{1,1}, S^{1,2}$ and $S^{2,1}, S^{2,2}$.

We shall show that the map $i: S_{k-1,n}\to S^1$ given by 
$([a_1,\ldots,a_{k-1}|\ldots]-\delta)\mapsto ([(2k-1)/2,a_1,\ldots,a_{k-1}|\ldots]-\delta)$
is a graph isomorphism. The weight $[a_1,\ldots,a_{k-1}|\ldots]$ is strictly 
$\liep_{k-1,n}$-dominant and $\liep_{k-1,n}$-integral
if and only if $[(2k-1)/2,a_1,\ldots,a_{k-1}|\ldots]$ is strictly $\liep_{k,n}$-dominant and
$\liep_{k,n}$-integral. Therefore, $i$ is a bijection on vertices.

Similarly as in the proof of Theorem \ref{kdirgvmevenorb}, it can be shown that
there is a nonzero standard homomorphism 
$M_{\liep_{k-1,n}}(\beta_1)\to M_{\liep_{k-1,n}}(\beta_2)$ if and only if
there is a nonzero standard homomorphism
$M_{\liep_{k,n}}(i(\beta_1))\to M_{\liep_{k,n}}(i(\beta_2))$. Therefore, 
$i: S_{k-1,n}\to S^1$ is a graph isomorphism.

Similarly, the map $j: S_{k-1,n}\to S^2$ given by 
$$([a_1,\ldots,a_{k-1}|\ldots]-\delta)\mapsto ([a_1,\ldots,a_{k-1},-(2k-1)/2|\ldots]-\delta)$$
is a graph isomorphism.

We will show now that all the arrows connecting elements of $S^{1,2}$ with elements of $S^{2,1}$ in $S_k$ (defined
analogously as in \ref{Sk})
are in the BGG graph for $(\lieg,\liep,\lambda)$. 

%
%
Let us denote $\mu+\delta=[(2k-1)/2,a_2,\ldots,a_{k-1},-(2k-3)/2|\ldots]$ and 
$\nu+\delta=[(2k-3)/2,a_2,\ldots,a_{k-1},-(2k-1)/2|\ldots]$. We will show that there is an arrow 
$\mu\to\nu$ in the BGG graph. First, we show that there is a standard homomorphism
$M_\liep(\nu)\to M_\liep(\mu)$. For true Verma modules, $M(\nu)\subset M(\mu)$, because
$\nu+\delta=s_\gamma(\mu+\delta)$ and $\nu=\mu-\gamma$ for $\gamma=[1,0\ldots,0,1|0,\ldots,0]$.

Let $w\in W$ be the Weyl group element that takes $\delta=\half[\ldots,5,3,1]$ to 
$\half[4k-3,b_2,\ldots,b_{k-1},-(4k-7)|\ldots]$ where $(b_2,\ldots,b_{k-1})$ is some decreasing sign-permutation
of $((4k-11)/2,\ldots,5/2,1/2)$ and $w'$ takes $\delta$ to $\half[4k-7,b_2,\ldots,b_{k-1},-(4k-3)|\ldots]$ so that
$w(\tilde\lambda+\delta)=\mu+\delta$ and 
$w'(\tilde\lambda+\delta)=\nu+\delta$. The difference of the grading element evaluation
is $(w\delta-w'\delta)(E)=4$ and, similarly as in the proof of lemma \ref{2diraksodd}, either $w\to w'$ or
$w\to w_1\to w_2\to w'$. If $w\to w'$, we apply Theorem \ref{genlepowski} and see that there is an arrow $\mu\to\nu$. 
Let $w\to w_1\to w_2\to w'$ and
assume, for the sake of contradiction, that the standard map
$M_\liep(w'\cdot \tilde\lambda)\to M_{\liep}(w\cdot\tilde\lambda)$ is zero. Therefore,
\begin{equation}
\label{cond2}
M(\nu)=M(w'\cdot\tilde\lambda)\subset M(s_\alpha w\cdot\tilde\lambda)=M(s_\alpha\cdot \mu)
\end{equation}
for some $\alpha\in S$.

The weight $s_\alpha(\mu+\delta)$ is one of the following types:
\begin{enumerate}
\item{$[a_2,(2k-1)/2,\ldots|\ldots,3,2,1]$ if $\alpha=\alpha_1$}
\item{$[(2k-1)/2,\ldots,a_{l},a_{l-1},\ldots,-(2k-3)/2|\ldots]$} 
\item{$[(2k-1)/2,\ldots,-(2k-3)/2,a_{k-1}|\ldots]$}
\item{$[(2k-1)/2,\ldots,-(2k-3)/2|n,n-1,\ldots,l-1,l,\ldots,2,1]$} 
\item{$[(2k-1)/2,\ldots,-(2k-3)/2|\ldots,3,2,-1]$ if $\alpha=\alpha_{n+k}$} 
\end{enumerate}

First we show that it is not of type $(1)$. If $\alpha=\alpha_1$, (\ref{cond2}) implies 
$\nu+\delta \leq s_{\alpha_1}(\mu+\delta)$, i.e. 
$$[(2k-3)/2,a_2,\ldots,-(2k-1)/2|\ldots]\leq [a_2, (2k-1)/2,\ldots,-(2k-3)/2|\ldots]$$ 
where $a_2\leq(2k-5)/2$.
Subtracting the first
weight from the second we get $[a_2-(2k-3)/2,\ldots]$
which cannot be obtained as a sum of positive roots, because it contains a negative number on the 
first position.

Now assume that $s_\alpha(\mu+\delta)$ is of type $(2)-(5)$. 
Because
$$
M(w'\cdot\tilde\lambda)=M(\nu)\subsetneq M(s_\alpha\cdot\mu)\subsetneq M(\mu)=M(w\cdot\tilde\lambda),
$$
$l(w')-l(w)=3$ and $\nu+\delta$ is not connected to $s_\alpha(\mu+\delta)$ with any root reflection, it follows from
Theorem \ref{truevermamodulesmap} that there must be $\beta_1, \beta_2$ so that
\begin{equation}
\label{subsetsgenodd}
M(\nu)\subsetneq M(s_{\beta_1}\cdot\nu)=M(s_{\beta_2} s_\alpha \cdot\mu)\subsetneq M(s_\alpha\cdot\mu)
\end{equation}
Similarly as in the proof of lemma \ref{2diraksodd}, we will show that $\alpha$ cannot be
of type $(2)-(5)$ leading to a contradiction. Let $\alpha$ be of type $(2)$, i.e. 
\begin{eqnarray*}
s_\alpha(\mu+\delta)&=&[(2k-1)/2,\ldots,a_{l},a_{l-1},\ldots,-(2k-3)/2|\ldots],\\
\nu+\delta&=&[(2k-3)/2,\ldots,a_{l-1},a_l,\ldots,-(2k-1)/2)|\ldots].
\end{eqnarray*}
The root reflections $s_{\beta_1}$ and $s_{\beta_2}$ cannot interchange an integer with a half-integer, because
of the integrality conditions $s_\alpha(\mu+\delta)(H_{\beta_2})\in\N$ and 
$s_{\beta_2}s_\alpha(\mu+\delta)(H_{\beta_1})\in\N$.
There are two possibilities: either $s_{\beta_2}$ interchanges $a_l$ with $a_{l-1}$ and $s_{\beta_1}$ interchanges
$((2k-1)/2,-(2k-3)/2)$ with $((2k-3)/2,-(2k-1)/2)$ on the particular positions, or $s_{\beta_2}$ interchanges
$((2k-1)/2,-(2k-3)/2)$ with $((2k-3)/2,-(2k-1)/2)$ and
$s_{\beta_1}$ interchanges $a_l$ with $a_{l-1}$.
In the first case, $\beta_2=\alpha$ and (\ref{subsetsgenodd}) implies $M(\mu)\subsetneq M(s_\alpha\cdot\mu)$, 
which contradicts the fact that $M(s_\alpha\cdot\mu)\subsetneq M(\mu)$
for a parabolic simple root $\alpha$ and $\mu\in P_\liep^{++}$.
In the second case, $\beta_1=\alpha$ and (\ref{subsetsgenodd}) implies
$M(\nu)\subsetneq M(s_\alpha\cdot\nu)$, which also contradicts $M(s_\alpha\cdot\nu)\subsetneq M(\nu)$.

Let $\alpha$ be of type $(3)$, i.e.
\begin{eqnarray*}
s_\alpha(\mu+\delta)&=&[(2k-1)/2,\ldots,-(2k-3)/2,a_{k-1}|\ldots],\\
\nu+\delta&=&[(2k-3)/2,\ldots,a_{k-1},-(2k-1)/2)|\ldots]
\end{eqnarray*}
If either $\beta_1=\alpha$ or $\beta_2=\alpha$, we get contradiction similarly as in case $(2)$. 
But there is no other possibility, because the $a_{k-1}$ on the $k$-th position has to move somehow
to the $(k-1)$-th position: if $\beta_2$ would fix it, then $\beta_1=\alpha$, if $\beta_2$ would take it
to the $(k-1)$-th position, then $\beta_2=\alpha$ and if $\beta_2$ would take it (possibly with a minus sign) 
to the $l$-th position for $l\neq k,k-1,1$, then $\beta_2$ has to (sign-) interchange the $l$-th and $(k-1)$-th
position, so $s_{\beta_1} s_{\beta_2}$ would fix the $(2k-1)/2$ on the first position, which is impossible.
The last possibility is $l=1$: this would mean that $\beta_2$ takes $a_{k-1}$ to the first position (possibly with a minus sign),
but $a_{k-1}<(2k-3)/2$ would imply that $s_{\beta_2}s_\alpha(\mu+\delta)$ has a smaller number on the first position
as $\nu+\delta$ which contradicts $\nu+\delta\leq s_{\beta_2}s_\alpha(\mu+\delta)$.


In case $(4)$, we have 
\begin{eqnarray*}
s_\alpha(\mu+\delta)&=&[(2k-1)/2,\ldots,-(2k-3)/2|n,\ldots,l-1,l,\ldots,2,1]\\
\nu+\delta&=&[(2k-3)/2,\ldots,-(2k-1)/2)|n,\ldots,l,l-1,\ldots,2,1]
\end{eqnarray*}
Because the reflections with respect to $\beta_1,\beta_2$ cannot interchange an integer and a 
half-integer, it follows that one of them interchanges $l$
with $l-1$, so either $\beta_1=\alpha$ or $\beta_2=\alpha$ and we get a contradiction as in case $(2)$. The same happens in case 
$(5)$.

In either case, we get a contradiction, so
the standard map
$M_\liep(\nu)\to M_\liep(\mu)$
is nonzero. 

We will show that the condition $(2)$ from the definition of the BGG graph (\ref{singularbgggraph}) is satisfied as well.
Let us suppose that there exists nontrivial homomorphisms 
$$M_\liep(\nu)=M_\liep(\xi_0)\to M_\liep(\xi_1)\to\ldots\to M_\liep(\xi_j)=M_\liep(\mu)$$
for $j>1$.
The weights $\xi_i$ are increasing and $\liep$-dominant $\liep$-integral, so $\xi_i-\xi_{i-1}$ cannot be written down 
as a sum of positive roots in $\lieg_0$: therefore, the grading element evaluation $\xi_i(E)$ is
strictly increasing.  The difference $(\mu-\nu)(E)=2$, so it follows $j=2$ and
$$\nu+\delta\leq\xi_1+\delta=\nu+\delta+\gamma_1\leq\nu+\delta+\gamma_1+\gamma_2=\mu+\delta$$
for some $\gamma_1,\gamma_2\in\lieg_1$ (grading element evaluation on roots from $\lieg_1$ is $1$). 
But $\gamma_1\in\lieg_1$ implies that $\gamma_1$ is of the form
$$[0,\ldots,1,\ldots,0|0,\ldots,\pm 1,0\ldots,0],$$ and so the $\gamma_1$-addition changes $$[a_1,\ldots,a_k|b_1,\ldots,b_n]$$
to $$[a_1,\ldots,a_l+1,\ldots|b_1,\ldots,b_{l'}-1,\ldots,b_n],$$ 
but this cannot be strictly $\liep$-dominant and $\liep$-integral
when $(b_1,\ldots,b_n)=(n,\ldots,2,1)$. So, there is really an arrow $\mu\to\nu$ in the BGG graph and 
we leave to the reader to check that there are no other arrows in the BGG graph.



For $k=n$, there are other $\liep$-integral and strictly $\liep$-dominant weights on the 
orbit of $\lambda$: the weights of type $[a_1,\ldots,a_k|(2k-1)/2,\ldots,3/2,1/2]$
where $(a_1,\ldots,a_k)$ is some strictly decreasing sign-permutation of $(k,k-1,\ldots,1)$.
We could again, define $S'^1$ as a set of weights $[k,\ldots|\ldots]-\delta$ and $S'^2$ as a 
set of weights $[\ldots,-k|\ldots]-\delta$ and show that there are arrows 
$$([k,\ldots,-(k-1)|\ldots]-\delta) \to ([k-1,\ldots,-k]-\delta),$$ similarly as in the
first part of the proof. However, the order of the dual operators is always $2$ in this case, because
changing $1$ to $-1$ decreeses the grading element evaluation by $2$.
\end{proof}

\newpage
\section{Calculus of extremal vectors}
\label{singularvectors}

\subsection{Computation of extremal vectors}
We saw in the last chapter that in the even orthogonal case, the
singular Hasse graph contains arrows that may correspond to nonstandard homomorphisms
and the dual operators are of second order. One way to prove that existence of a nonstandard 
homomorphism \gvmhom is to compute the so called {\it extremal vector}. It is a vector
$v\in M_\liep(\lambda)$ of weight $\mu$ such that $X\cdot v=0$ for any
positive root space generator $X$ in $\lieg$. If we find such a vector, the homomorphism
can be defined by $1\otimes v_\mu\mapsto v$ ($1\otimes v_\mu$ is the highest weight vector
in $M_\liep(\mu)$) and consequently
$$y_1\ldots y_k\otimes Y_1\ldots Y_l v_\mu\mapsto y_1\ldots y_k Y_1\ldots Y_l v$$
where $y_i$ are generators of some negative root spaces in $\lieg_-$, 
$Y_j$ are generators of some negative root spaces in $\lieg_0$ and the right hand side
of the last equation is the result of the action of $y_1\ldots y_k Y_1\ldots Y_l$ on $v$
in $M_\liep(\lambda)$.

Consider the Lie algebra $\lieg=D_4=\so(8,\C)$, $\liep$ the parabolic subalgebra determined by $\Sigma=\{\alpha_1\}$ 
(first node crossed in the Dynkin diagram), and the weights $\lambda=\half[-5|1,1,1]$
and $\mu=\half[-7|1,1,-1]$ (this is a special case of lemma \ref{1dirgvmeven}). We will compute the
extremal vector corresponding to \gvmhom. 

Let us represent the elements of $\so(8,\C)$ as matrices antisymmetric with respect 
to the anti-diagonal, as in section \ref{intro}. 

Let $y_{i,j}$ resp. $Y_{i,j}$ be a matrix $E_{i,j}-E_{9-j,9-i}$ so that $y_{i,j}\in\lieg_-$ and $Y_{i,j}\in {\lieg_0}_-$ 
($E_{i,j}$ is a matrix having $1$ in $i$-th row and $j$-th column and $0$ on other places). These are exactly the generators
of negative root spaces in $\lieg$. Similarly, we denote the generators of positive root spaces by $x_{i,j}$ and $X_{i,j}$ and 
the generators of the Cartan subalgebra by $h_{i}=E_{i,i}-E_{9-i,9-i}$:

\begin{equation}
\left(
\begin{tabular}{c|cccccc|c}
$h_{1}$ & $x_{12}$ & $x_{13}$ & $x_{14}$ & $x_{15}$ & $x_{16}$ & $x_{17}$ & 0 \\
\hline
$y_{21}$ & $h_{2}$ & $X_{23}$ & $X_{24}$ & $X_{25}$ & $X_{26}$ & $0$ & $-x_{17}$\\
$y_{31}$ & $Y_{32}$ & $h_{3}$ & $X_{34}$ & $X_{35}$ & $0$  & $-X_{26}$ & $-x_{16}$\\
$y_{41}$ & $Y_{42}$ & $Y_{43}$ & $h_{4}$ & $0$ & $-X_{35}$  & $-X_{25}$ & $-x_{15}$\\
$y_{51}$ & $Y_{52}$ & $Y_{53}$ & $0$ & $-h_4$ & $-X_{34}$  & $-X_{24}$ & $-x_{14}$\\
$y_{61}$ & $Y_{62}$ & $0$& $-Y_{53}$ & $-Y_{43}$ & $-h_3$ & $-X_{23}$ & $-x_{13}$\\
$y_{71}$ & $0$& $-Y_{62}$ & $-Y_{52}$ & $-Y_{42}$ & $-Y_{32}$ & $-h_2$ & $-x_{12}$\\
\hline
$0$ & $-y_{71}$ & $-y_{61}$ &  $-y_{51}$ &  $-y_{41}$ &  $-y_{31}$ &  $-y_{21}$ & $-h_{1}$ \\
\end{tabular}\right)
\label{matrix}
\end{equation}

\bigskip
%
\begin{lemma}
\label{1dirgvmspecialextremal}
There is exactly one vector (up to multiple) in $M_\liep(\lambda)$ of weight $\mu$ that is extremal, i.e. 
annihilated by all positive root spaces in $\lieg$,
namely the vector
$$y_{5,1}\otimes v_\lambda - y_{31} \otimes Y_{53}v_\lambda - y_{21}\otimes Y_{52} v_\lambda$$
(under the identification $M_\liep(\lambda)\simeq \univ(\lieg_-)\otimes \V_\lambda$).
\end{lemma}

\begin{proof}
We know from \ref{vermamodules-general} that the vector 
\begin{equation}
\label{basis}
y_{i_1,j_1}\ldots y_{i_n,j_n} \otimes Y_{k_1,l_1}\ldots Y_{k_m,l_m} v_\lambda
\end{equation}
is a weight vector with weight
$\lambda-\sum_k {\root}\, (y_{i_k,j_k}) - \sum_{k'}{\root} \,(Y_{i_{k'}, j_{k'}})$ (if it is nonzero).
The difference $\mu-\lambda=[-1|0,0,-1]$ in our case, so the $\mu$-weight space in $M_\liep(\lambda)$ is generated by
vectors of type $(\ref{basis})$, where the sum 
$$\sum_k {\root}\, (y_{i_k,j_k}) + \sum_{k'}{\root} \,(Y_{i_{k'}, j_{k'}})=[-1|0,0,-1].$$
There are only $4$ possibilities how to obtain $[-1|0,0,-1]$ as a sum of negative roots in $\lieg$:
\begin{itemize}
\item{$[-1|0,0,-1]$ itself -- corresponds to $y_{51}$, so the weight vector is $y_{51}\otimes v_\lambda$}
\item{$[0|-1,0,-1] + [-1|1,0,0]$ -- weight vector  $y_{21} \otimes Y_{52} v_\lambda$}
\item{$[0|0,-1,-1] + [-1|0,1,0]$ -- weight vector $y_{31} \otimes Y_{53} v_\lambda$}
\item{$[0|0,-1,-1] + [0|-1,1,0] + [-1|1,0,0,]$ -- weight vector $y_{21} \otimes Y_{53}Y_{32} v_\lambda$.}
\end{itemize}

The last vector is zero because $\lambda=\frac{1}{2}[-5|1,1,1]$, $Y_{32}$ is the negative root space of the root 
$\beta=[0,-1,1,0]$ 
and the copy of $\sl(2,\C)$ in $\lieg$ generated by $H_\beta, X_{23}, Y_{32}$ acts trivial on $v_\lambda$, 
because $H_\beta(v_\lambda)=\lambda(H_\beta)v_\lambda=\lambda(h_2-h_3)(1-1)v_\lambda=0$ and therefore,
this submodule generated by $v_\lambda$ is an irreducible $\sl(2,\C)$-module with highest weight $0$,
so it is trivial. From the same reason, $Y_{43}v_\lambda=Y_{42}v_\lambda=0$. 
On the other hand, $Y_{52}v_\lambda\neq 0$ and $Y_{53}v_\lambda\neq 0$, because this root spaces correspond
to the coroot $h_1+h_3$ resp. $h_2+h_3$ and the action of $\lambda$ on these coroots is nonzero.

We have identified a $3$-dimensional $\mu$-weight space in $M_\liep(\lambda)$ and are looking for a vector in this space that is 
extremal, i.e. annihilated by all positive root spaces in $\lieg$. 
The action of the positive root spaces can be computed using just the commutation
relations in $\univ(\lieg)$ and the fact that we know the action of $\liep$ on $v_\lambda$. 

In fact, it suffices to find a vector in this weight space that is annihilated by 
$x_{12}, X_{23}, X_{34}, X_{35}, X_{26}$ and $x_{17}$,
because the other positive root spaces generators can be obtained by commuting those.
We compute the action of $x_{12}$ on the three vectors:
\begin{eqnarray*}
&& x_{12} (y_{51}\otimes v_\lambda)=y_{51} x_{12}\otimes v_\lambda + [x_{12},y_{51}] \otimes v_\lambda=\\
&& y_{51}\otimes x_{12}v_\lambda + [x_{12},y_{51}]\otimes v_\lambda=
0 + (-Y_{52})\otimes v_\lambda\\
&&=1\otimes (-Y_{52} v_\lambda),\\
&&\\
&& x_{12}(y_{21}\otimes Y_{52}v_{\lambda})=y_{21}x_{12}\otimes Y_{52}v_\lambda + [x_{12},y_{21}]\otimes Y_{52}v_\lambda =\\
&& =y_{21}\otimes x_{12}Y_{52}v_\lambda + (h_1-h_2)\otimes Y_{52}v_\lambda=y_{21}\otimes Y_{52}x_{12} v_\lambda + \\ 
&& + y_{21} \otimes [x_{12},Y_{52}]v_\lambda+ 1\otimes (h_1-h_2) Y_{52 }v_\lambda = 0 + 0 +\\
&& + 1\otimes Y_{52}(h_1-h_2)v_\lambda + 1\otimes [h_1-h_2,Y_{52}]= \\ 
&& = 1\otimes (-\frac{5}{2}-\frac{1}{2}) v_\lambda + 1\otimes Y_{52}v_\lambda=-2\otimes Y_{52}v_\lambda, \\
&& \\
&& x_{12} (y_{31}\otimes Y_{53}v_\lambda)=y_{31}\otimes x_{12}Y_{53} v_\lambda + [x_{12},y_{31}]\otimes Y_{53}v_\lambda=\\
&& = y_{31}\otimes Y_{53}x_{12} v_\lambda + y_{31}\otimes [x_{12},Y_{53}] v_\lambda + (-Y_{32}) \otimes Y_{53}v_\lambda=\\
&& =0+0-1\otimes Y_{32}Y_{53}v_\lambda=-Y_{53}Y_{32}v_\lambda-[Y_{32},Y_{53}] v_\lambda=\\
&& = 0 - 1\otimes (-Y_{52})v_\lambda=1\otimes Y_{52}v_\lambda,
\end{eqnarray*}
where $\otimes$ means product over $\univ(\liep)$.

Similarly, we compute the action of the other positive root spaces on each of the $3$ nonzero
vectors of weight $\mu$. We can write the result into a table of actions on vectors:

\begin{center}
\begin{tabular}{|c|c|c|c|}
\hline
action & \multicolumn{3}{c|}{vector}\\
\hline
 & $v_1=y_{51}\otimes v_\lambda$ & $v_2=y_{21}\otimes Y_{52}v_{\lambda}$ & $v_3=y_{31}\otimes Y_{53}v_\lambda$\\
\hline
\hline
$x_{12}$ & $1\otimes (-Y_{52} v_\lambda)$ & $-2\otimes Y_{52}v_\lambda$ & $1\otimes Y_{52}v_\lambda$\\ 
\hline
$X_{23}$ & $0$ & $-y_{21}\otimes Y_{53}v_\lambda$ & $y_{21}\otimes Y_{53}v_\lambda$\\
\hline
$X_{34}$ & $0$ & $0$ & $y_{21}\otimes Y_{34}v_\lambda=0$\\
\hline
$X_{35}$ & $y_{31}\otimes v_\lambda$ & $y_{21}\otimes Y_{32}v_\lambda=0$ & $y_{31}\otimes v_\lambda$\\
\hline
$X_{26}$ & $0$ & $0$ & $0$ \\
\hline
$x_{17}$ & $0$ & $0$ & $0$ \\
\hline
\end{tabular}
\end{center}

We want to find some combination $a v_1+b v_2+c v_3$ of the vectors $v_1, v_2, v_3$ such that the actions on
this are zero. We see from the table that $a-2b+c=0$, $b=c$ and $a=-c$, so the solution is one-dimensional
$(a,b,c)\in\br{(1,-1,-1)}$.
\end{proof}

We know from \ref{1dirgvmeven} that there exist a standard homomorphism \gvmhom, but now we see that
it is the only one (up to multiple), and there is no nonstandard homomorphism \gvmhom.

In a similar way, we could compute that taking $\lambda'=\half[-5|1,1,-1]$ and $\mu'=\half[-7|1,1,1]$
there exists a unique homomorphism $M_\liep(\mu')\to M_\liep(\lambda')$ and the extremal vector is
$$y_{41}\otimes v_{\lambda'}-y_{31}\otimes Y_{43}v_{\lambda'}-y_{21}\otimes Y_{42}v_{\lambda'}.$$

This example can be generalized:

\begin{lemma}
\label{1dirgvmevenextremal}
Let $\lieg=D_{n+k}$ and $\Sigma=\{\alpha_k\}$, $\lambda+\delta=[2k-1,\ldots,3,1|\ldots,3,1]$ and 
$\mu+\delta=[2k-1,\ldots,3,-1|\ldots,3,-1]$, as in \ref{kdirgvmeven}. Let us represent elements of
$\lieg$ as matrices, the same way as in \ref{intro}.
Then the Verma module homomorphism $M_\liep(\mu)\to M_{\liep}(\lambda)$ is unique (up to multiple) 
and is described by the extremal vector
\begin{eqnarray}
&& \mmmm \label{extrdirgvmeven_gen_kn} y_{k+n+1,k}\otimes v_\lambda - y_{k+n-1,k} \otimes Y_{k+n+1,k+n-1}v_\lambda - \\
\nonumber && \mmmm -y_{k+n-2,k}\otimes Y_{k+n+1,k+n-2} v_\lambda -\ldots - y_{k+1,k}\otimes Y_{k+n+1,k+1} v_\lambda.
\end{eqnarray}
Similarly, for $\lambda'+\delta=[2k-1,\ldots,3,1|\ldots,3,-1]$ and 
$\mu'+\delta=[2k-1,\ldots,3,-1|\ldots,3,1]$, the extremal vector is 
\begin{eqnarray*}
&& \mmmm y_{k+n,k}\otimes v_\lambda - y_{k+n-1,k} \otimes Y_{k+n,k+n-1}v_\lambda - y_{k+n-2,k}\otimes Y_{k+n,k+n-2} v_\lambda -\ldots\\
&& \mmmm \ldots -y_{k+1,k}\otimes Y_{k+n,2} v_\lambda.
\end{eqnarray*}
\end{lemma}

\begin{proof}
The technique of the proof is essentially the same as in the previous lemma. The weight difference 
$\mu-\lambda=[0,\ldots,-1|0,\ldots,-1]$ which can be obtained either directly (it is the root so 
that its root space generator is $y_{k+n+1,k}$), or as a sum
$$[0,\ldots,-1|,\ldots,1,0,\ldots,0]+[\ldots|\ldots,-1,\ldots,-1].$$
The weight vector corresponding to this decomposition is $$y_{k+j,k}\otimes Y_{k+n+1,k+j}v_\lambda.$$
Other decompositions of type $(\ldots+[0,\ldots,0|\ldots,-1,\ldots,1\ldots])$ do not occur, because the
root space generator associated to this root is $Y_{k+l,k+m}$ for $l<m$ and $Y_{k+l,k+m}v_\lambda=0$,
because $\lambda(h_{k+l}-h_{k+m})=\half-\half=0$ and we can use the same argument as in the proof of 
lemma \ref{1dirgvmspecialextremal}.

Writing down the actions of the positive root spaces on this weight vectors, we obtain that the 
only weight vector annihilated by them is the vector from the lemma. It suffices to take the action
of $X_{12},X_{23},\ldots, X_{k-1,k}$, $x_{k,k+1}$, $X_{k+1, k+2}$, $\ldots$, 
$X_{k+n-1,k+n}$, $X_{k+n-1,k+n+1}$,$\ldots$, $X_{k+1,k+2n-1}$, $x_{k, k+2n}$, $\ldots, x_{1,2k+2n-1}$.
Recall that 
$$\lambda=\half[-(2n-1),\ldots,-(2n-1)|1,1,\ldots,1].$$ We compute the action of $X_{12}$ on the 
summands of $(\ref{extrdirgvmeven_gen_kn})$ in case $k>1$:
\begin{eqnarray*}
&& \mmmm X_{12} (y_{k+n+1,k}\otimes_{\univ(\liep)} v_\lambda)=y_{k+n+1,k}\otimes_{\univ(\liep))} 
X_{12}v_\lambda+\\
&& \mmmm +[X_{12},y_{k+n+1,k}]\otimes_{\univ(\liep)} v_\lambda=0,\\
&& \mmmm \\
&& \mmmm X_{12} (y_{k+j,k}\otimes Y_{k+n+1,k+j}v_\lambda)=y_{k+j,k}X_{12}\otimes_{\univ(\liep)}Y_{k+n+1,k+j}v_\lambda+\\
&& \mmmm +[X_{12},y_{k+j,k}]\otimes_{\univ(\liep)}v_\lambda=y_{k+j,k}\otimes_{\univ(\liep)}Y_{k+n+1,k+j} X_{12}v_\lambda+\\
&& \mmmm +y_{k+j,k}\otimes_{\univ(\liep)}[X_{12},Y_{k+n+1,k+j}]v_\lambda+0=0,
\end{eqnarray*}
for $j=1,\ldots,n-1$, because the commutators are zero. 
Similarly, the action of $X_{23},\ldots, X_{k-1,k}$ are zero on all the vectors we consider. For $x_{k,k+1}$,
however (if $k=1$, we start here), there are nonzero commutators $[x_{k,k+1}, y_{k+n+1,k}]=-Y_{k+n+1,k+1}$, 
$[x_{k,k+1},y_{k+j,k}]=-Y_{k+j,k+1}$ for $j=2,\ldots,n-1$ and $[x_{k,k+1}, y_{k+1,k}]=h_k-h_{k+1}$. Note that
$\lambda(h_k-h_{k+1})=-n$. We obtain
\begin{eqnarray*}
&& \mmmm x_{k,k+1}(y_{k+n+1,k}\otimes_{\univ(\liep)} v_\lambda)=-1\otimes Y_{k+n+1,k+1}v_\lambda,\\
&& \mmmm \\
&& \mmmm x_{k,k+1}(y_{k+j,k}\otimes_{\univ(\liep)} Y_{k+n+1,k+j} v_\lambda)=-1\otimes Y_{k+j,k+1} Y_{k+n+1,k+j}v_\lambda=\\
&& \mmmm =-1\otimes (Y_{k+n+1,k+j} Y_{k+j, k+1} v_\lambda + [Y_{k+j,k+1},Y_{k+n+1,k+j}]v_\lambda)=\\
&& \mmmm =1\otimes Y_{k+n+1,k+1}v_\lambda, 
\end{eqnarray*}
for $j=2,\ldots,n-1$, because $Y_{k+j,k+1}v_\lambda=0$.
Finally, 
\begin{eqnarray*}
&& \mmmm x_{k,k+1} (y_{k+1,k}\otimes Y_{k+n+1,k+1})=(h_k-h_{k+1})\otimes_{\univ(\liep)}Y_{k+n+1, k+1}v_\lambda=\\
&& \mmmm =1\otimes (Y_{k+n+1,k+1} (h_k-h_{k+1})+[h_{k}-h_{k+1},Y_{k+n+1,k+1}])v_\lambda=\\
&& \mmmm =1\otimes (-n+1) Y_{k+n+1,k+1})v_\lambda=(-n+1)\otimes Y_{k+n+1,k+1}v_\lambda.
\end{eqnarray*}
We see that the action of $x_{k,k+1}$ on the vector from the lemma is 
$$1\otimes Y_{k+n+1,k+1}(-1-(n-2)+n-1)=0.$$
Further, the action of $X_{k+i,k+i+1}$ for some fixed $i\in \{1,\ldots,n-1\}$ is zero on 
$y_{k+n+1,k}\otimes_{\univ(\liep)} v_\lambda$ and the only nonzero terms come from the commutators
$[X_{k+i,k+i+1},y_{k+i+1,k}]=y_{k+i,k}$ for $1\leq i\leq n-2$ and $[X_{k+i,k+i+1}, Y_{k+n+1, k+i}]=-Y_{k+n+1,k+i+1}$.
We obtain:
\begin{eqnarray*}
&& \mmmm X_{k+i,k+i+1} (y_{k+i+1,k}\otimes Y_{k+n+1,k+i+1}v_\lambda)=y_{k+i,k} \otimes Y_{k+n+1,k+i+1} v_\lambda
\end{eqnarray*}
and
\begin{eqnarray*}
&& \mmmm X_{k+i,k+i+1} (y_{k+i,k}\otimes Y_{k+n+1,k+i}v_\lambda)=-y_{k+i,k}\otimes Y_{k+n+1,k+i+1} v_\lambda.
\end{eqnarray*}
Therefore, the coefficient of the term $y_{k+i,k}\otimes Y_{k+n+1,k+i}v_\lambda$ in the extremal vector 
is the same as the coefficient of the term $y_{k+i+1,k}\otimes Y_{k+n+1, k+i+1}v_\lambda$. This already
implies the extremal vector from the lemma and the reader can check that the action of the other 
positive root spaces is zero as well.
\end{proof}

\begin{lemma}
\label{2dirgvmspecialextremal}
Let $\lieg=D_4,\,\, \Sigma=\{\alpha_2\}$ and $\lambda+\delta=\half[3,-1|3,-1]$, $\mu+\delta=\half[1,-3|3,-1]$
are the weights representing the double-arrow from lemma \ref{2dirgvmeven}. Then there is a unique 
(up to multiple)
nonzero nonstandard homomorphism \gvmhom, i.e. the arrow $\lambda\to\mu$ is in the BGG graph.
\end{lemma}

\begin{proof}
We compute the extremal vector the same way as in the previous lemmas. 
The matrices are $8\times 8$ with the following gradation:
\begin{equation}
\label{matrix2grad22}
\left(
\begin{tabular}{cc|cccc|cc}
$h_{1}$ & $X_{12}$ & $x_{13}$ & $x_{14}$ & $x_{15}$ & $x_{16}$ & $x_{17}$ & 0 \\
$Y_{21}$ & $h_{2}$ & $x_{23}$ & $x_{24}$ & $x_{25}$ & $x_{26}$ & $0$ & $-x_{17}$\\
\hline
$y_{31}$ & $y_{32}$ & $h_{3}$ & $X_{34}$ & $X_{35}$ & $0$  & $-x_{26}$ & $-x_{16}$\\
$y_{41}$ & $y_{42}$ & $Y_{43}$ & $h_{4}$ & $0$ & $-X_{35}$  & $-x_{25}$ & $-x_{15}$\\
$y_{51}$ & $y_{52}$ & $Y_{53}$ & $0$ & $-h_4$ & $-X_{34}$  & $-x_{24}$ & $-x_{14}$\\
$y_{61}$ & $y_{62}$ & $0$& $-Y_{53}$ & $-Y_{43}$ & $-h_3$ & $-x_{23}$ & $-x_{13}$\\
\hline
$y_{71}$ & $0$& $-y_{62}$ & $-y_{52}$ & $-y_{42}$ & $-y_{32}$ & $-h_2$ & $-X_{12}$\\
$0$ & $-y_{71}$ & $-y_{61}$ &  $-y_{51}$ &  $-y_{41}$ &  $-y_{31}$ &  $-Y_{21}$ & $-h_{1}$ \\
\end{tabular}\right)
\end{equation}

The computation is quite long and 
technical but the result is that the extremal vector is
\begin{eqnarray*}
&& \mmmm v_{ext}=y_{51}y_{42}\otimes v_\lambda+y_{52}y_{42}\otimes Y_{21}v_\lambda+y_{62}y_{31}\otimes v_\lambda+y_{62}y_{32}\otimes Y_{21}v_\lambda+\\
&& \mmmm + y_{52}y_{31}\otimes Y_{43}v_\lambda-y_{51}y_{32}\otimes Y_{43}v_\lambda.
\end{eqnarray*}
For simplicity, we will write only $y\in\univ(\lieg)$ instead of $y\otimes_{\univ(\liep)}v_\lambda$ 
for the extremal vector.
To make the expression
of $y$ more unique, we will write it as a sum of elements of type 
$y_{i_1,j_1}\ldots y_{i_l,j_l} Y_{u_1,v_1}\ldots, y_{u_k,v_k}$
so that $$((i_1,j_1),\ldots,(i_l,j_l))\quad \and \quad ((u_1,v_1),\ldots, (u_k,v_k))$$ are ordered lexicographically.
Of course, $yx=0$ for any positive root space generator $x$ in this formalism, $Y^k$ may be zero for large $k$ and
$yh=\lambda(h)y$ for $h\in\lieh$.

So, we write $v_{ext}=y_{51}y_{42}+\ldots +(-y_{51}y_{32}Y_{43})$. The reader may easily verify
that all the summands have weights $\mu=\lambda+[-1,-1,|0,0]$. For example, the first summand $y_{51}y_{42}$ corresponds
to the decomposition $[-1,0|0,-1]+[0,-1|0,1]$ of $[-1,-1,0,0,]$.

We will show that the action of $X_{12}$ on this is zero.
\begin{eqnarray*}
&& \mmmm X_{12}y_{51}y_{42}=y_{51}x_{12}y_{42}+[X_{12},y_{51}]y_{42}=y_{51}y_{42}X_{12}+y_{51}[X_{12},y_{42}]-y_{52}y_{42}=\\
&& \mmmm =-y_{52}y_{42},\\
&& \mmmm \\
&& \mmmm X_{12}y_{52}y_{42}Y_{21}=y_{52}X_{12}y_{42}Y_{21}+[X_{12},y_{52}]\ldots =y_{52}y_{42}X_{12}Y_{21}+y_{52}[X_{12},y_{42}]\ldots=\\
&& \mmmm =y_{52} y_{42} Y_{21}X_{12}+y_{52}y_{42}[X_{12}, Y_{21}]=y_{52}y_{42}(h_1-h_2)=y_{52}y_{42},\\
&& \mmmm \quad\quad \mathrm{(because }\,\, \lambda=\half[-3/2,-5/2|1/2,-1/2]\mathrm{)}\\
&& \mmmm \\
&& \mmmm X_{12} y_{62} y_{31}=y_{62}X_{12} y_{31}=y_{62}y_{31}X_{12}+y_{62}[X_{12},y_{31}]=-y_{62}y_{32},\\
&& \mmmm \\
&& \mmmm X_{12}y_{62}y_{32}Y_{21}=y_{62}X_{12}y_{32}Y_{21}=y_{62}y_{32}X_{12}Y_{21}=y_{62}y_{32}Y_{21}X_{12}+\\
&& \mmmm +y_{62}y_{32}[X_{12},Y_{21}]=y_{62}y_{32}(h_1-h_2)=y_{62}y_{32},\\
&& \mmmm \\
&& \mmmm X_{12} y_{52}y_{31}Y_{43}=y_{52}X_{12}y_{31}Y_{43}=y_{52}y_{31}X_{12}Y_{43}+y_{52}[X_{12},y_{31}]Y_{43}=\\
&& \mmmm =y_{52}y_{31}Y_{43}X_{21}-y_{52}y_{32}Y_{43}=-y_{52}y_{32}Y_{43},\\
&& \mmmm \\
&& \mmmm X_{12}(-y_{51}y_{32}Y_{43})=-y_{51}X_{12}y_{32}Y_{43}-[X_{12},y_{51}]y_{32}Y_{43}=-y_{51}y_{32}X_{12}Y_{43}+\\
&& \mmmm +y_{52}y_{32}Y_{43}=-y_{51}y_{32}Y_{43}X_{21}+y_{52}y_{32}Y_{43}=y_{52}y_{32}Y_{43}.
\end{eqnarray*}
Summing up all the results, we get zero. The same can be checked for any positive root space $X$ 
and it is left to the reader. To prove the uniqueness
of $v_{ext}$, one has to write down the basis of the $\mu$-weight space in $M_\liep(\lambda)$ and compute the actions 
of positive root spaces on them: solving this is very technical but straightforward.
\end{proof}

This extremal vector represented by $y_{ext}\in\univ(\lieg)$ can be rewritten in an easier way:
$$y_{ext}=(y_{51}-y_{31}Y_{53})(y_{42}-y_{32}Y_{43})+(y_{52}-y_{32}Y_{53})(y_{42}-y_{32}Y_{43})Y_{21}-y_{71}.$$
To check that it is the same, multiply the brackets, use commutation relations in $\lieg$ and the facts
that $yY_{53}=0$ (because $Y_{53}v_\lambda=0$).

The vectors in the bracket look very similar to those from lemma \ref{1dirgvmspecialextremal} and \ref{1dirgvmevenextremal}.
The following theorem says that this holds in general.

\begin{theorem}
\label{2dirgvmevenextremal}
Let $\lieg=D_{2+n}$, $\Sigma=\{\alpha_2\}$, $\lambda+\delta=\half[3,-1|\ldots,3,-1]$ and $\mu+\delta=[1,-3|\ldots,3,-1]$.
Then there exists a nonzero homomorphism \gvmhom and the extremal vector is 
\begin{equation}
\label{extrvect}
y_{ext}=D_1^+ D_2^- + D_2^+ D_2^- Y_{21} - y_{2n+3,1}
\end{equation}
($y_{2n+3,1}$ is the generator of the root space in $\lieg_{-2}$), 
where 
\begin{equation}
\label{D+}
D_i^+=y_{n+3,i}-y_{n+1,i}Y_{n+3,n+1}-y_{n,i}Y_{n+3,n}-\ldots -y_{3,i}Y_{n+3,3}
\end{equation}
and
\begin{equation}
\label{D-}
D_i^-=y_{n+2,i}-y_{n+1,i}Y_{n+2,n+1}-y_{n,i}Y_{n+2,n}-\ldots -y_{3,i}Y_{n+2,3}
\end{equation}
for $i=1,2$. 
\end{theorem}

\begin{remark}
Note that $D_{2}^+$ resp. $D_2^-$ is the extremal vectors from Theorem \ref{1dirgvmevenextremal} describing
the homomorphism dual to the Dirac operator.
\end{remark}

\begin{proof}
We will verify that the extremal vector from the theorem is annihilated by the action of positive root space generators.
This is sufficient to verify for $X_{12}, x_{23}, X_{i,i+1}$ ($3\leq i\leq n+1$), $X_{n+2-j,n+2+j}$ ($1\leq j\leq n-1$), 
$x_{2,2n+2}$ and $x_{1,2n+3}$ (draw a matrix similar to $(\ref{matrix2grad22})$ with $\lieg_0\simeq \gl(2,\C)\oplus \so(2n,\C)$ 
for general $n$). 

\hbox{
\hskip -10pt
\vbox{
\begin{center}
\begin{equation*}
{\tiny
\label{matrix2gradn2}
\left(
\begin{tabular}{cc|cccccc|cc}
$h_{1}$ & $X_{12}$ & $x_{13}$ & \ldots  & $x_{1,n+2}$ & $x_{1,n+3}$ & \ldots & $x_{1,2n+2}$ & $x_{1,2n+3}$ & $0$ \\
$Y_{21}$ & $h_{2}$ & $x_{23}$ & \ldots  & $x_{2,n+2}$ & $x_{2,n+3}$ & \ldots & $x_{2,2n+2}$ & $0$ & $-x_{1,2n+3}$\\
\hline
$y_{31}$ & $y_{32}$ & $h_{3}$ &  \ldots  & $X_{3,n+2}$ & $X_{3,n+3}$ & \ldots  & $0$ & $-x_{2,2n+2}$ & $-x_{1,2n+2}$\\
\ldots & \ldots &\ldots &\ldots & \ldots &\ldots  & \ldots & \ldots &\ldots &\ldots \\
$y_{n+2,1}$ & $y_{n+2,2}$ & $Y_{n+2,3}$ & \ldots & $h_{n+2}$ & $0$ & \ldots & $-X_{3,n+3}$ & $-x_{2,n+3}$ & $-x_{1,n+3}$\\
$y_{n+3,1}$ & $y_{n+3,2}$ & $Y_{n+3,3}$ & \ldots & $0$ & $-h_{n+2}$ &\ldots & $-X_{3,n+2}$ & $-x_{2,n+2}$ & $-x_{1,n+2}$\\ 
\ldots & \ldots &\ldots &\ldots & \ldots &\ldots  & \ldots & \ldots &\ldots &\ldots \\
$y_{2n+2,1}$ & $y_{2n+2,2}$ &  $0$ & \ldots & $-Y_{n+1,3}$ & $-Y_{n,3}$ & \ldots & $-h_3$ & $-x_{2,3}$ & $-x_{1,3}$\\
\hline 
$y_{2n+3,1}$ & $0$ & $-y_{2n+2,2}$ &  \ldots & $-y_{n+3,2}$ & $-y_{n+2,2}$ & \ldots & $-y_{32}$ & $-h_2$ & $-X_{12}$\\
$0$ & $-y_{2n+3,1}$ & $-y_{2n+2,1}$ & \ldots & $-y_{n+3,1}$ & $-y_{n+2,1}$ & \ldots & $-y_{31}$ & $-Y_{21}$ & $-h_1$
\end{tabular}\right)
}
\end{equation*}
\end{center}
}}
The weight $\lambda=\half[-\frac{2n-1}{2},-\frac{2n+1}{2}|\half,\ldots ,\half,-\half]$ implies that 
$Y_{n+3,j} v_\lambda=0$ for $3\leq j\leq n+1$ and $Y_{j,k}v_\lambda=0$ for $j\leq n+1$.

First we show that the extremal vector in question is annihilated by the positive root spaces in $\lieg_0$.

(a) The action of $X_{12}$.
\begin{eqnarray*}
&& \mmmm X_{12}D_1^+D_2^-=D_1^+ X_{12}D_2^- + [X_{12},D_1^+]D_2^-=\\
&& \mmmm =D_1^+D_2^- X_{12}+D_1^+ [X_{12}, D_2^-] +[X_{12},D_1^+]D_2^-\\
\end{eqnarray*}
The commutator $[X_{12},D_2^-]$ is zero, because $X_{12}$ commutes with each summand in the definition ($\ref{D-}$) of
$D_2^-$. Further, $[X_{12},D_1^+]=-D_2^+$ because $[X_{12},y_{n+3,1}]=-y_{n+3,2}$ and 
$$[X_{12},y_{n+2-k,1}Y_{n+3,n+2-k}]=-y_{n+2-k,2}Y_{n+3,n+2-k}$$ for $1\leq k\leq n-1.$
Therefore, 
$$X_{12}D_1^+ D_2^-=-D_2^+ D_2^-$$
Similarly, 
$$X_{12}D_2^+ D_2^- Y_{21}=D_2^+ D_2^-[X_{12},Y_{12}]=D_2^+ D_2^-,$$
because $[X_{12}, Y_{12}]v_\lambda=(h_1-h_2)v_\lambda=v_\lambda$ for 
$$\lambda=[-\frac{2n-1}{2}, -\frac{2n+1}{2}|\half,\ldots,\half,-\half].$$
So, $X_{12}(D_1^+ D_2^- + D_2^+ D_2^- Y_{21})=0$. The space $\lieg_{-2}$ is one-dimensional,
generated by $y_{2n+3,1}$. The action of $X_{12}$ is $X_{12}y_{2n+3,1}=[X_{12},y_{2n+3,1}]=0$
so we see that $X_{12}y=0$ for $y$ defined by (\ref{extrvect}).

(b) The action of $X_{i,i+1}$ for $3\leq i\leq n$.
Again, we have 
$$
X_{i,i+1}D_1^+ D_2^-=[X_{i,i+1},D_1^+]D_2^- + D_1^+ [X_{i,i+1},D_2^-]
$$
and 
$$
X_{i,i+1}D_2^+ D_2^-Y_{21}=[X_{i,i+1},D_2^+]D_2^-Y_{21} + D_1^+ [X_{i,i+1},D_2^-]Y_{21}
$$
We will show that all the commutators are zero. We have
\begin{eqnarray*}
&& \mmmm [X_{i,i+1},D_1^+]=[X_{i,i+1},y_{n+3,1}-y_{n+1,1}Y_{n+3,n+1}-\ldots -y_{31}Y_{n+3,3}].
\end{eqnarray*}
The only nonzero terms here are $[X_{i,i+1}, -y_{i+1,1}Y_{n+3,i+1}]=-y_{i,1}Y_{n+3,i+1}$ and
$[X_{i,i+1}, -y_{i,1} Y_{n+3,i}]=y_{i,1} Y_{n+3, i+1}$ and  they cancel each other.
Similarly,
\begin{eqnarray*}
&& \mmmm [X_{i,i+1},D_2^+]=[X_{i,i+1},y_{n+3,2}-y_{n+1,2}Y_{n+3,n+1}-\ldots -y_{32}Y_{n+3,3}]=\\
&& \mmmm =-y_{i,2}Y_{n+3, i+1}+y_{i,2}Y_{n+3, i+1}=0,\\
&& \mmmm \\
&& \mmmm [X_{i,i+1},D_1^-]=[X_{i,i+1},y_{n+2,1}-y_{n+1,1}Y_{n+2,n+1}-\ldots -y_{31}Y_{n+2,3}]=\\
&& \mmmm =-y_{i,1}Y_{n+2,i+1}+y_{i,1}Y_{n+2,i+1}=0 ,\\
&& \mmmm \\
&& \mmmm [X_{i,i+1},D_2^-]=[X_{i,i+1},y_{n+2,2}-y_{n+1,2}Y_{n+2,n+1}-\ldots -y_{32}Y_{n+2,3}]=\\
&& \mmmm =-y_{i,2}Y_{n+2, i+1}+y_{i,2}Y_{n+2, i+1}=0.
\end{eqnarray*}
The action of $X_{i,i+1}$ on $\lieg_{-2}$ is $X_{i,i+1} y_{2n+3}=0$ as well. 
So, we are done for $X_{i,i+1}$, $3\leq i\leq n$. 

(c) Action of $X_{n+1,n+2}$. We have
\begin{eqnarray*}
&& \mmmm [X_{n+1,n+2}, D_1^+]=[X_{n+1,n+2},y_{n+3,1}-y_{n+1,1}Y_{n+3,n+1}-\ldots -y_{31}Y_{n+3,3}]=0,
\end{eqnarray*}
because all the summands in $D_1^+$ commute with $X_{n,n+1}$. Similarly, $[X_{n+1,n+2}, D_2^+]=0$.
The other commutators are nontrivial:
\begin{eqnarray*}
&& \mmmm [X_{n+1,n+2}, D_1^-]=[X_{n+1,n+2}, y_{n+2,1}-y_{n+1,1}Y_{n+2,n+1}-\ldots -y_{3,1}Y_{n+2,3}]=\\
&& \mmmm =y_{n+1,1}-y_{n+1,1}(h_{n+1}-h_{n+2})-y_{n,1}Y_{n+1,n}-\ldots - y_{3,1}Y_{n+1,3},\\
&& \mmmm \\
&& \mmmm [X_{n+1,n+2}, D_2^-]=[X_{n+1,n+2}, y_{n+2,2}-y_{n+1,2}Y_{n+2,n+1}-\ldots -y_{3,2}Y_{n+2,3}]=\\
&& \mmmm =y_{n+1,2}-y_{n+1,2}(h_{n+1}-h_{n+2})-y_{n,2}Y_{n+1,n}-\ldots - y_{3,2}Y_{n+1,3}.\\
\end{eqnarray*}
We obtain:
\begin{eqnarray*}
&&\mmmm X_{n+1,n+2}(D_1^+ D_2^- + D_2^+ D_2^- Y_{21})=D_1^+ [X_{n+1,n+2}, D_2^-]+D_2^+ [X_{n+1,n+2},D_2^-]=\\
&&\mmmm =D_1^+ (y_{n+1,2}-y_{n+1,2}(h_{n+1}-h_{n+2})-y_{n,2}Y_{n+1,n}-\ldots - y_{3,2}Y_{n+1,3})+\\
&&\mmmm +D_2^+(y_{n+1,2}-y_{n+1,2}(h_{n+1}-h_{n+2})-y_{n,2}Y_{n+1,n}-\ldots - y_{3,2}Y_{n+1,3})Y_{21}.
\end{eqnarray*}
The first term is zero, because $Y_{n+1,n+1-j}$ has zero action on $v_\lambda$ for $1\leq j\leq n-2$ 
(because the weight $\lambda$ has $\half$ both on the $n+1$ and $(n+1-j)$'th position) and 
$y_{n+1,2}(1-(h_{n+1}-h_{n+2}))v_\lambda=y_{n+1,2}(1-(\half-(-\half)))v_\lambda=0$. 
For the second term, 
\begin{eqnarray*}
&&\mmmm (y_{n+1,2}-\ldots - y_{3,2}Y_{n+1,3})Y_{21}=Y_{21}(y_{n+1,2}-\ldots -y_{32}Y_{n+1,3})+\\
&&\mmmm +[Y_{21}, (y_{n+1,2}-\ldots -y_{32}Y_{n+1,3})]=-(y_{n+1,1}-y_{n+1,1}(h_{n+1}-h_{n+2})-\\
&&\mmmm -y_{n,1}Y_{n+1,n}-\ldots- y_{3,1}Y_{n+1,3}),
\end{eqnarray*}
and this is zero from the same reason as above.
For the $\lieg_{-2}$ term, $X_{n,n+1}y_{2n+3,1}=[X_{n,n+1}y_{2n+3}]=0$ as well.

(d) Action of $X_{n+1,n+3}$. 
The commutators are
\begin{eqnarray*}
&&\mmmm [X_{n+1,n+3}, D_1^+]=[X_{n+1, n+3}, y_{n+3,1}-y_{n+1,1}Y_{n+3,n+1}-\ldots -y_{31}Y_{n+3,3}]=\\
&&\mmmm =y_{n+1,1}-y_{n+1,1}(h_{n+1}+h_{n+2})-y_{n,1}Y_{n+1,n}-\ldots -y_{3,1}Y_{n+1,3},\\
&&\mmmm \\
&&\mmmm [X_{n+1,n+3}, D_2^+]=[X_{n+1, n+3}, y_{n+3,2}-y_{n+1,2}Y_{n+3,n+1}-\ldots -y_{32}Y_{n+3,3}]=\\
&&\mmmm =y_{n+1,2}-y_{n+1,2}(h_{n+1}+h_{n+2})-y_{n,2}Y_{n+1,n}-\ldots -y_{3,2}Y_{n+1,3},\\
&&\mmmm \\
&&\mmmm [X_{n+1,n+3},D_1^-]=[X_{n+1,n+3},D_2^-]=0.
\end{eqnarray*}

Therefore,
\begin{eqnarray*}
&&\mmmm X_{n+1,n+3}(D_1^+ D_2^- + D_2^+ D_2^- Y_{21})=\\
&&\mmmm =[X_{n+1,n+3}, D_1^+] D_2^- + [X_{n+1, n+3}, D_2^+]D_2^- Y_{21}.
\end{eqnarray*}
The first term is zero, because 
\begin{eqnarray*}
&&\mmmm =(y_{n+1,1}-y_{n+1,1}(h_{n+1}+h_{n+2})-y_{n,1}Y_{n+1,n}-\ldots -y_{3,1}Y_{n+1,3})D_2^-=\\
&&\mmmm =y_{n+1,1}D_2^- - y_{n+1,1}D_2^-(h_{n+1}+h_{n+2}) - y_{n+1,1}[(h_{n+1}+h_{n+2}),D_2^-]-\\
&&\mmmm -y_{n+1,1}[Y_{n+1,n},D_2^-]-\ldots -y_{3,1}[Y_{n+1,3},D_2^-],
\end{eqnarray*}
because the action of $Y_{n+1,j}$ on $v_\lambda$ is zero for $3\leq j\leq n$. It is easy to check
from (\ref{D-}) that the commutators 
\begin{equation*}
[Y_{n+1,j},D_2^-]=-y_{n+1,2}Y_{n+2,j}+y_{n+1,2}Y_{n+2,j}=0
\end{equation*}
and $(h_{n+1}+h_{n+2})v_\lambda=(\half-\half)v_\lambda=0$,
so we obtain
$$X_{n+1,n+3}D_1^+ D_2^-=y_{n+1,1} (D_2^- - [(h_{n+1}+h_{n+2}),D_2^-]).$$
A bit effort yields that $[h_{n+1}+ h_{n+2}, D_2^-]=D_2^-$, so 
$$X_{n+1,n+3} D_1^+ D_2^-=0.$$
Similarly, we can show that $X_{n+1,n+3}D_2^- D_2^+ Y_{21}=0$ and clearly, 
$$[X_{n+1,n+3},y_{2n+3,1}]=0,$$
so it has zero action on $\lieg_{-2}$ as well.

(e) Action of $X_{n+2-j,n+2+j}$ for $2\leq j\leq n-1$. In this case, the commutators are
\begin{eqnarray*}
&&\mmmm [X_{n+2-j, n+2+j}, D_1^+]=[X_{n+2-j, n+2+j}, y_{n+3,1}-y_{n+1,1}Y_{n+3,n+1}-\ldots -y_{31}Y_{n+3,3}]\\
&&\mmmm =-y_{n+3-j,1}[X_{n+2-j,n+2+j}, Y_{n+3,n+3-j}] -y_{n+2-j,1} [X_{n+2-j,n+2+j}, Y_{n+3,n+2-j}]=\\
&&\mmmm \hbox{(the other commutators are zero)}\\
&&\mmmm =y_{n+3-j,1}X_{n+2-j, n+2} -y_{n+2-j,1} X_{n+3-j,n+2},\\
&&\mmmm \\
&&\mmmm [X_{n+2-j, n+2+j}, D_2^+]= y_{n+3-j,2}X_{n+2-j, n+2} -y_{n+2-j,2} X_{n+3-j,n+2}\,\,\hbox{(similarly),}\\
&&\mmmm \\
&&\mmmm [X_{n+2-j, n+2+j},D_1^-]=y_{n+3-j,1}X_{n+2-j, n+3} -y_{n+2-j,1} X_{n+3-j,n+3},\\
&&\mmmm \\
&&\mmmm [X_{n+2-j, n+2+j}, D_2^-]= y_{n+3-j,2}X_{n+2-j, n+3} -y_{n+2-j,2} X_{n+3-j,n+3}.
\end{eqnarray*}
Therefore, 
\begin{eqnarray*}
&&\mmmm X_{n+2-j, n+2+j}(D_1^+ D_2^-)=[X_{n+2-j, n+2+j}, D_1^+] D_2^-=\\
&&\mmmm =y_{n+3-j,1} X_{n+2-j, n+2} D_2^- - y_{n+2-j,1} X_{n+3-j,n+2}D_2^-,
\end{eqnarray*}
To show that this is zero, observe that the weight $\lambda=\half[-(2n-1),-(2n+1)|1,\ldots,1,-1]$
differs from the weight $\lambda'=\half[-(2n-1),-(2n-1)|1,\ldots,1,-1]$ from lemma \ref{1dirgvmevenextremal}
only on the second position and we know from that lemma that the positive root spaces $X$ in $\lieg_0$ have zero
action on $D_2^- v_{\lambda'}$; a simple check shows that the second position in $\lambda$ is not used
in the computations in this case.

Similarly, we can show that $X_{n+2-j,n+2+j}D_2^+ D_2^- Y_{21}=0$ and $$X_{n+2-j, n+2+j}y_{2n+3,1}=0$$
as well.

(f) Now, we compute the action of $x_{23}$, a root space generator from $\lieg_1$. The commutators are
\begin{eqnarray*}
&&\mmmm [x_{23}, D_1^+]=[x_{23},y_{n+3,1}-y_{n+1,1}Y_{n+3,n+1}-\ldots -y_{31}Y_{n+3,3}]=\\
&&\mmmm =-Y_{21}Y_{n+3,3},\\
&&\mmmm \\
&&\mmmm [x_{23}, D_2^+]=[x_{23},y_{n+3,2}-y_{n+1,2}Y_{n+3,n+1}-\ldots -y_{32}Y_{n+3,3}]=\\
&&\mmmm =-Y_{n+3,3}+Y_{n+1,3}Y_{n+3,n+1}+\ldots + Y_{43}Y_{n+3,4}-(h_{2}-h_{3})Y_{n+3,3}=\\
&&\mmmm =-Y_{n+3,3}+(Y_{n+3,n+1}Y_{n+1,3}+[y_{n+1,3}, Y_{n+3,n+1}])+\ldots -(Y_{n+3,3}(h_2 -h_3)+\\
&&\mmmm +[(h_2 -h_3),Y_{n+3,3}])=-n Y_{n+3,3} - Y_{n+3,3}(h_2 -h_3)+Y_{n+3,n+1}Y_{n+1,3}+\ldots \\
&&\mmmm \ldots+Y_{n+3,4}Y_{43}.
\end{eqnarray*}
Similarly, we obtain
\begin{eqnarray*}
&&\mmmm [x_{23}, D_1^-]=-Y_{21}Y_{n+2,3},\\
&&\mmmm \\
&&\mmmm [x_{23}, D_2^-]=-n Y_{n+2,3} - Y_{n+2,3}(h_2 -h_3)+Y_{n+2,n+1}Y_{n+1,3}+\ldots \\
&&\mmmm \ldots+Y_{n+2,4}Y_{43}.
\end{eqnarray*}
Let us now compute 
\begin{equation}
\label{firstpart}
x_{23}D_1^+ D_2^-=[x_{23},D_1^+]D_2^- + D_1^+[x_{23}, D_2^-].
\end{equation}
The first term is
\begin{eqnarray*}
&&\mmmm [x_{23}, D_1^+] D_2^-=-Y_{21}Y_{n+3,3}D_2^-=-Y_{n+3,3}Y_{21}D_2^-=-Y_{n+3,3}D_2^- Y_{21}-\\
&&\mmmm -Y_{n+3,3}[Y_{21},D_2^-]=-Y_{n+3,3}D_2^- Y_{21} + Y_{n+3,3}D_1^- = -Y_{n+3,3}D_2^- Y_{21}+\\
&&\mmmm +[Y_{n+3,3}, D_1^-],
\end{eqnarray*}
where we used the relations $[Y_{21}, D_2^-]=-D_2^+$ and $Y_{n+3,3}v_\lambda=0$.

The second term is
\begin{eqnarray*}
&&\mmmm D_1^+[x_{23},D_2^-]=D_1^+(-n Y_{n+2,3} - Y_{n+2,3}(h_2 -h_3)+Y_{n+2,n+1}Y_{n+1,3}+\ldots\\
&&\mmmm \ldots + Y_{n+2,4}Y_{43})=D_1^+ Y_{n+2,3}(-n-(-\frac{2n+1}{2}-\half))=D_1^+ Y_{n+2,3}
\end{eqnarray*}
(the other terms are zero because $Y_{k,3}v_\lambda=0$ for $4\leq k\leq n+1$).
So, $(\ref{firstpart})$ is computed to be 
\begin{equation}
\label{firstpartresult}
-Y_{n+3,3}D_2^- Y_{21} + [Y_{n+3,3}, D_1^-]+ D_1^+ Y_{n+2,3}.
\end{equation}

Further, we compute
\begin{equation}
\label{secondpart}
x_{23}D_2^+ D_2^- Y_{21}=[x_{23},D_2^+]D_2^- Y_{21} + D_2^+ [x_{23}, D_2^-] Y_{21}
\end{equation}
We will show that the second term is zero:
\begin{eqnarray*}
&&\mmmm [x_{23}, D_2^-]Y_{21}=(-n Y_{n+2,3} - Y_{n+2,3}(h_2 -h_3)+Y_{n+2,n+1}Y_{n+1,3}+\ldots \\
&&\mmmm \ldots+Y_{n+2,4}Y_{43})Y_{21}=-n Y_{n+2,3}Y_{21}-Y_{n+2,3}Y_{21}(h_2-h_3)-\\
&&\mmmm -Y_{n+2,3}[Y_{21},h_2 -h_3]=Y_{n+2,3}Y_{21}(-n-(-\frac{2n+1}{2}-\half)-1)=0,
\end{eqnarray*}
where we used that $Y_{21}$ commutes with $Y_{k,3}$ and $Y_{k,3}v_\lambda=0$ again.
The first term in $(\ref{secondpart})$ is 
\begin{eqnarray*}
&&\mmmm [x_{23}, D_2^+]D_2^- Y_{21}=(-n Y_{n+3,3} - Y_{n+3,3}(h_2 -h_3)+Y_{n+3,n+1}Y_{n+1,3}+\ldots \\
&&\mmmm \ldots+Y_{n+3,4}Y_{43})D_2^- Y_{21}=-n Y_{n+3,3}D_2^- Y_{21} - Y_{n+3,3}D_2^- (h_2-h_3)Y_{21}-\\
&&\mmmm -Y_{n+3,3}[h_2-h_3, D_2^-]Y_{21}+Y_{n+3,n+1}Y_{n+1,3}D_2^-Y_{21}+\ldots +Y_{n+3,4}Y_{43}D_2^- Y_{21}.
\end{eqnarray*}
The commutator 
\begin{eqnarray*}
&&\mmmm [X_{k,3},D_2^-]=[X_{k,3}, y_{n+2,2}- y_{n+1,2}Y_{n+2,n+1}-\ldots -y_{32}Y_{n+2,3}]=\\
&&\mmmm [X_{k,3}, -y_{k,2}Y_{n+2,k}]-[X_{k,3},y_{32}Y_{n+2,3}]\,\,\,\hbox{(the other terms are zero)}=\\
&&\mmmm =y_{k,2}Y_{n+2,3}-y_{k,2}Y_{n+2,3}=0
\end{eqnarray*}
for $4\leq k\leq n+1$. Therefore, 
$$Y_{n+3,k}Y_{k,3}D_2^- Y_{21}=Y_{k,3}D_2^- Y_{21}Y_{n+3,k}=0$$
because $Y_{n+3,k}v_\lambda=0$. Another few calculations show that the commutator
$$[h_2-h_3,D_2^-]=-D_2^-.$$
We obtain
\begin{eqnarray*}
&&\mmmm [x_{23},D_2^+]D_2^- Y_{21}=-n Y_{n+3,3}D_2^- Y_{21} - Y_{n+3,3}D_2^- Y_{21}(h_2-h_3)- \\
&&\mmmm -Y_{n+3,3}D_2^-[h_2-h_3,Y_{21}] + Y_{n+3,3}D_2^- Y_{21}=\\
&&\mmmm =Y_{n+3,3}D_2^-Y_{21} (-n-(-\frac{2n+1}{2}-\half)-1+1)=Y_{n+3,3}D_2^-Y_{21}.
\end{eqnarray*}
So, $(\ref{secondpart})$ is equal to
\begin{equation}
\label{secondpartresult}
Y_{n+3,3}D_2^-Y_{21}.
\end{equation}
Summing up $(\ref{firstpartresult})$ and $(\ref{secondpartresult})$ we obtain
\begin{equation}
\label{partialresult}
x_{23}(D_1^+ D_2^- + D_2^+ D_2^- Y_{21})=[Y_{n+3,3}, D_1^-]+D_1^+ Y_{n+2,3}.
\end{equation}
The commutator is
\begin{eqnarray*}
&&\mmmm [Y_{n+3,3}, D_1^-]=[Y_{n+3,3}, y_{n+2,1}-y_{n+1,1}Y_{n+2,n+1}-\ldots -y_{31}Y_{n+2,3}]=\\
&&\mmmm =-y_{2n+2,1}-y_{n+2,1}Y_{n+4,3}-y_{n,1}Y_{n+5,3}-\ldots-y_{41}Y_{2n+1,3}-\\
&&\mmmm -y_{n+3,1}Y_{n+2,3}
\end{eqnarray*}
and
\begin{eqnarray*}
&&\mmmm D_1^+ Y_{n+2,3}=(y_{n+3,1}-y_{n+1,1}Y_{n+3,n+1}-\ldots -y_{31}Y_{n+3,3})Y_{n+2,3}=\\
&&\mmmm =y_{n+3,1}Y_{n+2,3}-y_{n+1,1}(Y_{n+2,3}Y_{n+3,n+1} + [Y_{n+3,n+1}, Y_{n+2,3}]) -\ldots\\
&&\mmmm \ldots - y_{41} (Y_{n+2,3}Y_{n+3,4}+[Y_{n+3,4}, Y_{n+2,3}]) - y_{31}Y_{n+2,3}Y_{n+3,3}.
\end{eqnarray*}
Because $Y_{n+3,k}v_\lambda=0$ for $3\leq k\leq n+1$ and the commutators are 
$$[Y_{n+3,k}, Y_{n+2,3}]=-Y_{2n+5-k,3}$$ for $4\leq k\leq n+1$ and $[Y_{n+3,3}, Y_{n+2,3}]=0$, we obtain
\begin{eqnarray*}
&&\mmmm D_1^+ Y_{n+2,3}=y_{n+3,1}Y_{n+2,3} + y_{n+1,1} Y_{n+4,3} + y_{n,1} Y_{n+5,3}+\ldots +y_{41}Y_{2n+1,3}.
\end{eqnarray*}
Summing this with the above equation, we obtain that $(\ref{partialresult})$ is equal to
$-y_{2n+2,1}$.

Further the action of $x_{23}$ on $y_{2n+3,1}$ is
$$
x_{23}y_{2n+3,1}=[x_{23}, y_{2n+3,1}]=-y_{2n+2,1}.
$$
We obtain the desired equality
$$
x_{23}(D_1^+ D_2^- + D_2^+ D_2^- Y_{21} - y_{2n+3,1})=0.
$$

(g) The action of $x_{2,2n+2}$. We will skip the details for now. The commutators are
\begin{eqnarray*}
&&\mmmm [x_{2,2n+2}, D_1^+]=y_{31}Y_{2,n+2},\\
&&\mmmm \\
&&\mmmm [x_{2,2n+2}, D_2^+]=(n-1) X_{3,n+2}-Y_{n+3,n+1}X_{3, n+4}-Y_{n+3,n}X_{3,n+5} -\ldots -\\
&&\mmmm -Y_{n+3,4}X_{3,2n+1}-y_{32}x_{2,n+2},\\
&&\mmmm \\
&&\mmmm [x_{23}, D_2^-]=(n-1) X_{3,n+3}-Y_{n+2,n+1}X_{3, n+4}-Y_{n+2,n}X_{3,n+5} -\ldots -\\
&&\mmmm -Y_{n+2,4}X_{3,2n+1}-y_{32}x_{2,n+3}.\\
\end{eqnarray*}
Further, we obtain
\begin{equation}
\label{hahaha}
x_{2,2n+2}D_1^- D_2^+=[x_{2,2n+2}, D_1^+]D_2^- + D_1^+ [x_{2,2n+2}, D_2^-]
\end{equation}
The second term is zero, because the commutator contains positive root spaces on the end of
each term, so $(\ref{hahaha})$ is 
\begin{eqnarray*}
&&\mmmm y_{31} x_{2,n+2}D_2^-=y_{31}[x_{2,n+2}, D_2^-]=y_{31}((h_2-h_{n+2})+(h_3-h_{n+2})-\ldots\\
&&\mmmm \ldots -(h_{n+1}-h_{n+2}))=y_{31}(-n+(n-1))=-y_{31}.
\end{eqnarray*}
For the second part,
\begin{eqnarray*}
&&\mmmm x_{2,2n+2}D_2^+ D_2^- Y_{21}=[x_{2,2n+2},D_2^+]D_2^- Y_{21}=(n-1) X_{3,n+2}-Y_{n+3,n+1}X_{3, n+4}-\\
&&\mmmm -Y_{n+3,n}X_{3,n+5} -\ldots -Y_{n+3,4}X_{3,2n+1}-y_{32}x_{2,n+2}.
\end{eqnarray*}
The action of $X_{3,n+2}$ and $X_{3,n+k}$ for $4\leq k\leq n+1$ on $D_2^- Y_{21}$ can be easily computed 
to be zero. Further, 
\begin{eqnarray*}
&&\mmmm [x_{2,n+2}, D_2^-]Y_{21}=((h_2-h_{n+2})+(h_3-h_{n+2})-\ldots -(h_{n+1}-h_{n+2}))Y_{21}=\\
&&\mmmm =Y_{21}((h_2-h_{n+2})+\ldots +(h_{n+1}-h_{n+2}))+Y_{21}=Y_{21}(-1+1)=0.
\end{eqnarray*}
Finally, 
$$x_{2,2n+2} y_{2n+3,1}=[x_{2,2n+2}, y_{2n+3,1}]=-y_{31}$$
and the action of $x_{2,2n+2}$ on $(\ref{extrvect})$ is zero.

(h) The action of $x_{1,2n+3}$. Because this is from $\lieg_2$, the action on $\ref{extrvect}$ is
\begin{equation}
\label{chacha}
[x_{1,2n+3}, D_1^+] D_2^- + [x_{1,2n+3}, D_2^+]D_2^- Y_{21} + [x_{1,2n+3}, -y_{2n+3,1}]
\end{equation}
The last term is clearly $h_1+h_2=-\frac{2n-1}{2}-\frac{2n+1}{2}=-2n$. For the first term,
\begin{eqnarray*}
&&\mmmm [x_{1,2n+3}, D_1^+]=x_{2,n+2}-x_{2,n+4}Y_{n+3,n+1}-x_{2,n+5}Y_{n+3,n}\ldots - x_{2,2n+2}Y_{n+3,3}\\
&&\mmmm \\
&&\mmmm [x_{2,n+2},D_2^-]=(h_2-h_{n+2})+(h_3-h_{n+2})+\ldots +(h_{n+1}-h_{n+2})=-1\\
&&\mmmm \\
&&\mmmm -x_{2,n+3+k}Y_{n+3,n+2-k}D_2^-=\sum_{2\leq j\leq n+1, j\neq n+2-k} (h_{n+2-k}+h_j)\,\, +\\
&&\mmmm +(h_{n+2-k}-h_{n+2})
\end{eqnarray*}
for $1\leq k\leq n-1$.
Summing this and substituting for $h_j$ the coordinates of $\lambda$, we obtain
$$x_{1,2n+3}D_1^+ D_2^-=-n.$$
For the second term of $(\ref{chacha})$ we obtain
\begin{eqnarray*}
&&\mmmm [x_{1,2n+3}, D_2^+]=-x_{1,n+2}+x_{1,n+4}Y_{n+3,n+1}+x_{1,n+5}Y_{n+3,n}\ldots + x_{1,2n+2}Y_{n+3,3},\\
&&\mmmm \\
&&\mmmm [-x_{1,n+2},D_2^-]Y_{21}=-X_{12}Y_{21}=-[X_{12},Y_{21}]=-1,\\
&&\mmmm \\
&&\mmmm x_{1,n+3+k}Y_{n+3,n+2-k}D_2^-Y_{21}=-X_{12}Y_{21}=-1
\end{eqnarray*}
for $1\leq k\leq n-1$.
Summing this we obtain
$$x_{1,2n+3}D_2^+ D_2^-Y_{21}=-n.$$
We see that the action of $x_{1,2n+3}$ on $(\ref{extrvect})$ is again $(-n-n-(-2n))=0$.
\end{proof}

We know from lemma $\ref{2dirgvmeven}$ that the homomorphism \gvmhom is not standard. 
We showed that the BGG graph and the singular Hasse graph coincide. Let us denote the
weights in this graph by $\lambda, \mu, \nu, \xi$:
\begin{center}
\begin{equation}
\label{weightnames}
\includegraphics{bgg.13}
\end{equation}
\end{center}
Similarly as in the last theorem (it is, in fact much easier, but we omit the computations), 
we could prove the following:
\begin{theorem}
\label{thirdoperatorextrvect}
The homomorphism $M_\liep(\xi)\to M_\liep(\nu)$ is given by the extremal vector
\begin{equation*}
D_1^-+D_2^-Y_{21}
\end{equation*}
where $D_1^-$ and $D_2^-$ are again defined by $(\ref{D-})$.
\end{theorem}

\begin{theorem}
\label{complex}
The BGG graph above is a complex. 
\end{theorem}

\begin{proof}
Let $\lambda,\mu, \nu, \xi$ be the weights as in the picture above.
We will show that the composition $M_\liep(\nu)\stackrel{i}{\to} M_\liep(\mu)\stackrel{j}{\to} M_\liep(\lambda)$ is zero.
This is zero exactly if $j\circ i$ sends the highest weight vector in $M_\liep(\nu)$ to zero.
Let $\alpha=1\otimes v_\lambda$ be the highest weight vector in $M_\liep(\lambda)$.
We know that the map $j$ sends the highest weight vector $1\otimes v_\mu$ in $M_\liep(\mu)$ 
to $y_{ext}\otimes_{\univ(\liep)}v_\lambda$ where $y_{ext}=D_2^+$ (see Theorem \ref{1dirgvmevenextremal}). 
We will denote it, for simplicity, $y_{ext}\alpha$.
Because $j$ is a $\lieg$-homomorphism, it sends $\tilde{y}\otimes_{\univ(\liep)}v_\mu$ to
$\tilde{y} y_{ext}\alpha$ for any $\tilde{y}\in\univ(\lieg)$. The composition $j\circ i$ sends $1\otimes v_\nu$ to $\tilde{y} y_{ext} \alpha$, where
$\tilde{y}$ is the extremal vector from Theorem \ref{2dirgvmevenextremal}. Therefore, we need to show that
$$(D_1^+ D_2^- + D_2^+ D_2^- Y_{21}-y_{2n+3,1})D_2^+$$ has zero action on $v_\lambda$.
Recall now that $\lambda=[\ldots,|\half,\ldots,\half]$ so now the action of $Y_{i,j}$ on $v_\lambda$
is zero for $i\leq n+2$. First, we compute the term $D_2^- D_2^+$  that we will denote by
$\Delta_2$. A few computations yield a very nice and symmetric expression
\begin{equation*}
\Delta_2=y_{2n+2,2}y_{32} + y_{2n+1,2}y_{42}+\ldots y_{n+3,2} y_{n+2,2}
\end{equation*}
Moreover, assuming the weight $\lambda'=[\ldots|\half,\ldots,\half,-\half]$ (corresponding to the $S^-$ representation),
we obtain the same expression $\Delta_2=D_2^+ D_2^-$ in that case.
Note that $Y_{21}v_\lambda=0$, because $\lambda(H_{\alpha_1})=\lambda(h_1-h_2)=0$.
We compute
\begin{eqnarray*}
&&\mmmm D_1^+ D_2 ^- D_2^+ + D_2^+ D_2^- Y_{21}D_2^+=D_1^+ \Delta_2 + D_2^+ D_2^- [Y_{21}, D_2^+]=\\
&&\mmmm =D_1^+ \Delta_2 - D_2^+ D_2^- D_1^+
\end{eqnarray*}
For the second term, note that $D_1^+ v_\lambda$ is a weight vector of weight 
$$[-\frac{2n+1}{2},-\frac{2n-1}{2}|\half,\ldots,\half,-\half]$$ 
It is not $\liep$-dominant but, however, similarly as for $\mu$,
the action of $Y_{n+3,k}$ is zero on such weight vectors 
and we easily compute that the action of $D_2^+ D_2^-$ on such weight vectors
is again equal to $\Delta_2$.
Therefore, the equation above is equal to 
\begin{eqnarray*}
&&\mmmm [D_1^+,\Delta_2]=\bigl[(y_{n+3,1}-y_{n+1,1}Y_{n+3,n+1}-\ldots -y_{3,1}Y_{n+3,3}),\,\,\,\, (y_{2n+2,3}y_{32}+\ldots \\
&&\mmmm \ldots +y_{n+3,2}y_{n+2,2})\bigr].
\end{eqnarray*}
The commutators are 
\begin{eqnarray*}
&&\mmmm [y_{n+3,1},\Delta_2]=[y_{n+3,1}, y_{n+3,2}y_{n+2,2}]=y_{n+3,2}y_{2n+3,1}=y_{2n+3,1}y_{n+3,2},\\
&&\mmmm \\
&&\mmmm [-y_{j,1}Y_{n+3,j}, \Delta_2]=-y_{j,1}[Y_{n+3,j}, \Delta_2]-[y_{j,1}, \Delta_2]Y_{n+3,j}=\\
&&\mmmm =-y_{j,1}[Y_{n+3,j}, y_{j,2}y_{2n+5-j,2}+y_{n+3,2}y_{n+2,2}]-[y_{j,1}, y_{2n+5-j,2}]Y_{n+3,j}=\\
&&\mmmm \hbox{(the other commutators are zero)}\\
&&\mmmm =-y_{j,1}(y_{n+3,2}y_{2n+5-j,2} - y_{n+3,2}y_{2n+5-j,2})-y_{2n+3,1}Y_{n+3,j}=\\
&&\mmmm =-y_{2n+3,1}Y_{n+3,j}
\end{eqnarray*}
for $3\leq j\leq n+1$. So, we proved that
$$[D_1^+, \Delta_2]=y_{2n+3,1}D_2^+,$$
and it follows that
$$(D_1^+ D_2^- + D_2^+ D_2^- Y_{21} - y_{2n+3,1})D_2^+=0.$$ 

We will show now that the composition $M_\liep(\xi)\to M_\liep(\nu)\to M_\liep(\mu)$ is zero.
Using Theorem \ref{thirdoperatorextrvect} we have to prove that the action of 
\begin{equation}
\label{secondzero}
(D_1^- + D_2^- Y_{21})(D_1^+ D_2^- + D_2^+ D_2^- Y_{21} - y_{2n+3,1})
\end{equation}
on $v_\mu$ is zero. Note that the action of $D_2^+D_2^-$ on $v_\mu$ is
$\Delta_2$ and similarly, 
$$D_1^+ D_1^-=y_{2n+2,1}y_{31} + y_{2n+1,1}y_{41}+\ldots y_{n+3,1} y_{n+2,1}=:\Delta_1.$$
We can easily check that if $D_1^- D_1^+$ acts on a weight vector with weight 
$[\ldots|\half,\ldots,\half]$ (the $S^+$-representation of the $\so$-part), then
$D_1^- D_1^+=\Delta_1$ as well. We adjust $(\ref{secondzero})$:
\begin{eqnarray*}
&&\mmmm (D_1^-+ D_2^- Y_{21})(D_1^+ D_2^- + D_2^+ D_2^- Y_{21})=D_1^- D_1^+D_2^- + D_1^- D_2^+ D_2^- Y_{21} + \\
&&\mmmm + D_2^- Y_{21}D_1^+ D_2^- + D_2^-Y_{21}D_2^+ D_2^- Y_{21}=\\
&&\mmmm =\Delta_1 D_2^- + D_1^- \Delta_2 Y_{21}+D_2^- D_1^+ Y_{21} D_2^- + D_2^- [Y_{21}, D_1^+]D_2^- +\\
&&\mmmm + D_2^- D_2^+ Y_{21} D_2^- Y_{21} + D_2^- [Y_{21}, D_2^+] D_2^- Y_{21}.
\end{eqnarray*}
The commutators are $[Y_{21}, D_1^+]=0$, $[Y_{21}, D_2^+]=-D_1^+$ and $[Y_{21}, D_2^-]=-D_1^-$, so we have
\begin{eqnarray*}
&&\mmmm \ldots =\Delta_1 D_2^- + D_1^- \Delta_2 Y_{21} + D_2^- D_1^+ D_2^- Y_{21} + D_2^- D_1^+ [Y_{21}, D_2^-] +\\
&&\mmmm + D_2^- D_2^+ D_2^- Y_{21} Y_{21} - D_2^- D_2^+  D_1^- Y_{21}-D_2^- D_1^+ D_2^- Y_{21}.\\
\end{eqnarray*}
Further, note that $Y_{21} Y_{21} v_\mu=0$ because, as a representation of the copy of $\sl(2,\C)$
generated by $H_{\alpha_1}, X_{12}, Y_{21}$, $V_\mu$ generates an irreducible representation. But
$H_{\alpha_1}(v_\mu)=\mu(h_1-h_2)=-\frac{2n-1}{2}-(-\frac{2n+1}{2})=1$, so this representation is
$2$-dimensional and $Y_{21}^2 v_\mu=0$. Continuing, we get
\begin{eqnarray*}
&&\mmmm \ldots =\Delta_1 D_2^- - D_2^- \Delta_1 + (D_1^- \Delta_2 - \Delta_2 D_1^-)Y_{21}=\\
&&\mmmm =[\Delta_1, D_2^-] + [D_1^-, \Delta_2]Y_{21}
\end{eqnarray*}
Some more computations yields that
\begin{eqnarray*}
&&\mmmm [\Delta_1, D_2^-]=y_{2n+3,1} D_1^-=D_1^- y_{2n+3,1}\quad \hbox{and}\quad [D_1^-,\Delta_2]Y_{21}=\\
&&\mmmm =y_{2n+3,1}D_2^-Y_{21}=D_2^- Y_{21} y_{2n+3,1}
\end{eqnarray*}
This proves that $(\ref{secondzero})$ acts trivial on $v_\mu$.

\end{proof}

\subsection{Translation of the extremal vector to an operator}
Let us consider the weights $\lambda, \mu, \nu, \xi$ defined by $(\ref{weightnames})$.
We will now revise the results of \ref{morediracs} (in case $k=2$)
and add some further comments on
that. We start with a complex Lie algebra $\so(\C,\beta)$ of matrices $2(n+2)\times 2(n+2)$
fixing the scalar product $\beta(x,y)=\sum_j x_j y_{2n+5-j}$. The homomorphism
\gvmhom is described by the extremal vector $D_2^+=y_{n+3,2}- y_{n+1,2}Y_{n+3,n+1}\ldots -y_{32}Y_{n+3,3}$
in our standard formalism.
We choose another product
$$\gamma(x,y)=x_1 x_{2n+4} + x_2 x_{2n+3} +  \sum_{j=3}^{2n+2} x_j^2,$$
of signature $(2n+2,2)$. We denote the Lie algebra of real matrices fixing this
scalar product by $\so(2(n+2), \gamma)$. Its complexification $\so(2(n+2),\gamma)^c\simeq \so(\C, \gamma)$
is isomorphic to the complex Lie algebra 
$\so(\C, \beta)$, because all complex scalar products are conjugate.
It is easy to check that the explicit isomorphism $\varphi: \so(\C, \beta)\to \so(\C, \gamma)$
is given, in matrices, by $A\mapsto C^{-1}AC$ where $\gamma=C^t\beta  C$ and
{\small
\begin{equation*}
C=
\left(
\begin{tabular}{cc|cccc|cc}
$1$ & $0$ & $0$ & $0$ & $0$ & $0$ & $0$ & $0$ \\
$0$ & $1$ & $0$ & $0$ & $0$ & $0$ & $0$ & $0$\\
\hline
$0$ & $0$ & $\frac{1}{\sqrt{2}}$ & $0$& $0$ & $\frac{-i}{\sqrt{2}}$  & $0$ & $0$ \\
$0$ & $0$ & $0$ & $\frac{1}{\sqrt{2}}$ & $\frac{-i}{\sqrt{2}}$ & $0$ & $0$ & $0$ \\
$0$ & $0$ & $0$ & $\frac{1}{\sqrt{2}}$ & $\frac{i}{\sqrt{2}}$ & $0$ & $0$ & $0$ \\
$0$ & $0$ & $\frac{1}{\sqrt{2}}$ & $0$& $0$ & $\frac{i}{\sqrt{2}}$  & $0$ & $0$ \\
\hline
$0$ & $0$ & $0$ & $0$ & $0$ & $0$ & $1$ & $0$ \\
$0$ & $0$ & $0$ & $0$ & $0$ & $0$ & $0$ & $1$ \\
\end{tabular}\right)
\end{equation*}}
for $n=2$ and analogously for larger $n$. Let $\liep$ be the real Lie subalgebra of
$\so(2(n+2), \gamma)$ such that in the associated gradation $\lieg_0=\gl(2,\R)\oplus \so(2n)$.
We showed in \ref{therealversion} that for a complex representation $\V$ of the complex Lie algebra
$\so(\C, \gamma)^c$ (being also a representation of the real Lie algebra \hbox{$\so(2(n+2),\gamma)$} by restriction)
$$\univ(\so(2(n+2), \gamma))\otimes_{\univ(\liep)} \V\simeq \univ(\so(\C, \gamma))\otimes_{\univ(\liep^c)} \V
\simeq \univ(\so(\C, \beta))\otimes_{\univ(\liep')} \V,$$
where $\liep'=\varphi(\liep^c)$ is the parabolic subalgebra of $\so(\C,\beta)$ given by $\Sigma=\{\alpha_2\}$.
The homomorphism \gvmhom maps the highest weight vector in $M_\liep(\mu)$ to 
$D_2^+\otimes_{\univ(\liep')} v_\lambda=(y_{52}-y_{32}Y_{53})\otimes_{\univ(\liep')} v_\lambda$
in the case $n=2$ (see \ref{1dirgvmevenextremal}). 
However, for the Verma modules induced by $\so(\C, \gamma)$
the extremal vector is $\tilde{D}_2^+=(\varphi(y_{52})-\varphi(y_{32})\varphi(Y_{53}))\otimes_{\univ(\liep)} v_\lambda$
where $\varphi$ is the isomorphism $\so(\C,\beta)\to \so(\C, \gamma)$. We denote, for simplicity,
$\tilde{y}_{ij}:=\varphi(y_{ij})$ and will omit the tilde in $\tilde{D_i^\pm}$.  We obtain

{\tiny
\begin{equation}
\label{transl}
\tilde{y}_{52}=
\frac{1}{\sqrt{2}}\left(
\begin{tabular}{cc|cccc|cc}
$0$ & $0$ & $0$ & $0$ & $0$ & $0$ & $0$ & $0$ \\
$0$ & $0$ & $0$ & $0$ & $0$ & $0$ & $0$ & $0$\\
\hline
$0$ & $0$ & $0$ & $0$ & $0$ & $0$ & $0$ & $0$\\
$0$ & $1$ & $0$ & $0$ & $0$ & $0$ & $0$ & $0$\\
$0$ & $-i$ & $0$ & $0$ & $0$ & $0$ & $0$ & $0$\\
$0$ & $0$ & $0$ & $0$ & $0$ & $0$ & $0$ & $0$\\
\hline
$0$ & $0$ & $0$ & $-1$ & $i$ & $0$ & $0$ & $0$\\
$0$ & $0$ & $0$ & $0$ & $0$ & $0$ & $0$ & $0$\\
\end{tabular}\right)
\tilde{y}_{32}=
\frac{1}{\sqrt{2}}\left(
\begin{tabular}{cc|cccc|cc}
$0$ & $0$ & $0$ & $0$ & $0$ & $0$ & $0$ & $0$ \\
$0$ & $0$ & $0$ & $0$ & $0$ & $0$ & $0$ & $0$\\
\hline
$0$ & $1$ & $0$ & $0$ & $0$ & $0$ & $0$ & $0$\\
$0$ & $0$ & $0$ & $0$ & $0$ & $0$ & $0$ & $0$\\
$0$ & $0$ & $0$ & $0$ & $0$ & $0$ & $0$ & $0$\\
$0$ & $i$ & $0$ & $0$ & $0$ & $0$ & $0$ & $0$\\
\hline
$0$ & $0$ & $-1$ & $0$ & $0$ & $-i$ & $0$ & $0$\\
$0$ & $0$ & $0$ & $0$ & $0$ & $0$ & $0$ & $0$\\
\end{tabular}\right)\end{equation}}
Let $\lieg=\so(2(n+2),\gamma)$ be the real Lie algebra and $\liep$ its parabolic subalgebra, inducing the gradation
of $\lieg$.
We know that $\lieg_0^{ss}\simeq \sl(2,\R)\otimes \so(2n)$ and as a $\lieg_0^{ss}$-module,
$\lieg_{-1}\simeq (\R^2)^*\otimes \R^{2n}$. We can define
natural coordinates on $\lieg_{-1}$ denoted by $y_{1,1}, \ldots, y_{1,2n}, y_{2,1}, \ldots, y_{2,2n}$ 
such that if $\{\epsilon_1, \epsilon_2\}$ is a basis of $(\R^2)^*$ dual to $e_1, e_2$ and $e_j$ is a basis of 
$\R^{2n}$ so $\epsilon_j\otimes e_k$ has coordinates $y_{jk}=1$ and $y_{mn}=0$ for $m\neq j$ or $n\neq k$.
In \ref{morediracs}, we identified sections of $\Spin(2(n+2),\gamma)\times_P S^\pm$ with spinor valued functions
on $\lieg_-$ and further restricted to functions that are only functions on $\lieg_{-1}$. 
Let $f: \R^{4n}\to S^-$ be such a function and $Df$ its image, where $D$ is the differential operator dual to the real Verma modules homomorphism.
We showed that $D=(D_1, D_2)$, where $D_i$ is the Dirac operators $\sum_j e_j \frac{\partial}{\partial y_{ij}}$. 
Let $s^+$, resp. $s^-$, be the highest weight vector 
in $S^+$ resp. $S^-$ as $\so(2n)$-modules. As a $\lieg_0^{ss}$-modules, $\V_\lambda^*\simeq S^-$ and $\V_\mu^*\simeq \C^2\otimes S^+$ for
$n$ odd and $\V_\lambda^*\simeq S^+$ and $\V_\mu^*\simeq \C^2\otimes S^-$ if $n$ is even. We will assume the first case, the other one
is similar. 
The module $\V_\mu\simeq (\C^2)^*\otimes S^-$ has a highest weight vector $\epsilon_2\otimes s^-$.
Then we know from $(\ref{duality})$ and $(\ref{transl})$ that
\begin{eqnarray}
\label{mess}
\nonumber &&\mmmm (\epsilon_2\otimes s^-)(Df(0))=\frac{1}{\sqrt{2}}\frac{\partial}{\partial y_{22}} s^+(f) |_0 +
\frac{-i}{\sqrt{2}}\frac{\partial}{\partial y_{23}} s^+(f) |_0 -\\
&&\mmmm - \frac{1}{\sqrt{2}}\frac{\partial}{\partial y_{21}} (\tilde{Y}_{53} s^+) (f) |_0 -
 \frac{i}{\sqrt{2}}\frac{\partial}{\partial y_{24}} (\tilde{Y}_{53} s^+) (f) |_0 .
\end{eqnarray}
The left hand side is equal to $s^-((D_2 f)(0))$  and the right hand side can be interpreted as the
action of $D_2^+$ on $s^+(f)$. Formally,
\begin{eqnarray*}
&&\mmmm s^-((D_2 f)(0))=D_2^+ (s^+)(f),
\end{eqnarray*}
and for any $u=Y_1\ldots Y_k$ a product of negative root spaces in $\so(2n)$
\begin{eqnarray}
\label{D2D_2}
&&\mmmm (u s^-)((D_2 f)(0))=(u D_2^+) (s^+)(f)
\end{eqnarray}
(such $u s^-$ generate all $S^-$).
Similarly,
\begin{eqnarray*}
&&\mmmm (Y_{21} \epsilon_2\otimes s^-)((D f)(0))=(Y_{21} D_2^+) (s^+)(f).
\end{eqnarray*}
The left hand side is $-\epsilon_1\otimes s^-$ (because the action of the negative root space $Y_{21}$ in $\sl(2,\C)$ 
on $\epsilon_2$ is $-\epsilon_1$) and for the right hand side, note that
$Y_{21} v_\lambda=0$, because all the $\sl(2,\R)$ acts trivially on 
$\V_\lambda=\C\otimes S^+$.
So, we obtain
\begin{eqnarray*}
&&\mmmm (-\epsilon_1 \otimes s^-)((D f)(0))=[Y_{21}, D_2^+] (s^+)(f)=-D_1^+ (s^+) (f),
\end{eqnarray*}
(because we know from the proof of Theorem $\ref{complex}$ that $[Y_{21}, D_2^+]=-D_1^+$) and 
\begin{eqnarray*}
&&\mmmm (s^-)((D_1 f)(0))=D_1^+ (s^+) (f).
\end{eqnarray*}
Both sides are equal to an equation that differs from $(\ref{mess})$ only by
differentiating $\frac{\partial}{\partial y_{1j}}$ instead of $\frac{\partial}{\partial y_{2j}}$. 
Similarly as in $(\ref{D2D_2})$, we obtain that for any product $u=Y_1\ldots Y_k$ 
of negative root spaces in $\so(2n)$,
\begin{eqnarray}
\label{D2D_1}
&&\mmmm (u s^-)((D_1 f)(0))=(u D_1^+) (s^+) (f).
\end{eqnarray}
Informally, we can say that the action of $D_2^+$ or $D_1^+$ on $f$ is the action of 
$D_1$ or $D_2$. 

Similarly, if we start with representations $V_{\lambda'}=S^-$ and 
$\V_{\mu'}=(\C^2)^*\otimes S^+$, the $D_2^-$ and $D_1^-$ would act
as $D_2$ and $D_1$.

Now consider the homomorphism $M_\liep(\nu)\to M_\liep(\mu)$. As $\lieg_0^{ss}$-modules,
$\V_\nu=\V_\mu\simeq (\R^2)^*\otimes S^-$. Let us denote by ${\cal D}$ the operator dual to the
Verma module homomorphism $M_\liep(\nu)\to M_\liep(\mu)$ and again, restrict it to functions
on $\lieg_{-1}$. 
The extremal vector is $D_1^+ D_2^- + D_2^+ D_2^- Y_{21}-\tilde{y}_{2n+3,1}$ (see \ref{2dirgvmevenextremal}).
The highest weight vector in $\V_\nu^*$ is again $\epsilon_2\otimes s^+$ and
the duality in $(\ref{duality})$ yields
\begin{eqnarray}
\label{operatorcD}
&&\mmmm (\epsilon_2\otimes s^-)(\cD g)(0)=(D_1^+ D_2^- + D_2^+ D_2^- Y_{21})(\epsilon_2\otimes s^-)(g),
\end{eqnarray}
where on the right hand side, the symbols $D_i^\pm$ acts by differentiating, similarly as in $(\ref{mess})$
and the derivations are evaluated in $0$.
Note that we omitted the term $\tilde{y}_{2n+3,1}$, because it is from $\lieg_{-2}$ 
and we assume that the function $f$ is constant
in the $\lieg_{-2}$, so differentiating in this direction is trivial.
So, $\cD$ acts between the spaces
$$\cD: C^\infty((\R^2)^*\otimes \R^{2n}, \R^2\otimes S^+)\to C^\infty((\R^2)^*\otimes \R^{2n}, \R^2\otimes S^+).$$ 
We can again define the components $\cD_1, \cD_2$ of $\cD$ so that
$$\cD: {g_1\choose g_2} \mapsto {\cD_1 (g_1, g_2)\choose \cD_2(g_1, g_2)},$$
where $g_1, g_2$ are functions on $\lieg_{-1}$. The left side of $(\ref{operatorcD})$ says that
it describes the second component $\cD_2$. The terms on the right are 
$$D_1^+ D_2^- (\epsilon_2\otimes s^-)(g)=D_1^+ D_2^- (s^- g_2),$$ 
and
$$D_2^+ D_2^-Y_{21} (\epsilon_2\otimes s^-)(g)=D_2^+ D_2^- (-\epsilon_1\otimes s^-) (g)=-D_2^+ D_2^- (s^- g_1),$$ 
so $(\ref{operatorcD})$ reads
$$s^-(\cD_2 (g)(0))=D_1^+ D_2^- (s^- g_2) - D_2^+ D_2^- (s^- g_1).$$
Let $u=Y_1\ldots Y_k$ be some products of negative root spaces in $\so(2n)$. Then
\begin{equation}
\label{cDalmostfinal}
u s^-(\cD_2 (g)(0))=u D_1^+ D_2^- (s^- g_2) - u D_2^+ D_2^- (s^- g_1).
\end{equation}
We know already that $D_2^- (s^- g_i)(0)=s^+(D_2 g_i (0))$, where $D_2$ is the Dirac operator
in variables $y_{21}, \ldots, y_{2,2n}$. Because $(\ref{duality})$ holds not only in zero but everywhere,
we can easily show that $D_2^-(s^- g_i)=s^+ (D_2 g_i)$ everywhere. The function $g_i$ is $S^+$-valued, so
$D_2 g_i$ is $S^-$-valued and we know from $(\ref{D2D_2})$ and $(\ref{D2D_1})$ that 
$$u D_j^+ D_2^- (s^- g_i)=u D_j^+ (s^+ (D_2 g_i))=u s^- (D_j D_2 g_i)\quad (j=1,2).$$
Substituting the left hand side into $(\ref{cDalmostfinal})$ we get
$$u s^- (\cD_2 (g_1, g_2))=u s^- (D_1 D_2 g_2 - D_2 D_2 g_1).$$
This holds for any $u$ and the spinors $u s^-$ generate all $S^-$, so we obtain that
the second component of the operator $\cD$ is the following combination of Dirac operators:
$$\cD_2 (g_1, g_2)=D_1 D_2 g_2 - D_2 D_2 g_1.$$

To determine the first component, we apply the action of $Y_{21}$ on $(\ref{operatorcD})$:
\begin{equation}
(-\epsilon_1\otimes s^-)(\cD g)(0)=Y_{21}(D_1^+ D_2^- + D_2^+ D_2^- Y_{21})(\epsilon_2\otimes s^-)(g).
\end{equation}
Recall the commutators $[Y_{21}, D_1^+]=0$, $[Y_{21}, D_2^+]=-D_1^+$ and $[Y_{21}, D_2^-]=-D_1^-$, so we compute
\begin{eqnarray*}
&&\mmmm Y_{21}(D_1^+ D_2^- + D_2^+ D_2^- Y_{21})(\epsilon_2\otimes s^-)=(D_1^+ D_2^- Y_{21}-D_1^+ D_1^- + D_2^+ Y_{21} D_2^- Y_{21} -\\
&&\mmmm - D_1^+ D_2^- Y_{21})(\epsilon_2\otimes s^-)=(D_1^+ D_2^- Y_{21} - D_1^+ D_1^- - D_2^+ D_1^-Y_{21} -\\
&&\mmmm - D_1^+ D_2^- Y_{21})(\epsilon_2\otimes s^-)=(-D_1^+ D_1^- - D_2^+ D_1^-Y_{21})(\epsilon_2\otimes s^-).
\end{eqnarray*}
So, applying a similar argument as for $\cD_2$, we see that 
$$\cD_1(g_1, g_2)=D_1 D_1 g_2 - D_2 D_1 g_1.$$

We see that we have found a principle, how to translate extremal vector expressions that are in terms of $D_i^\pm$
into operator. The method could be summarized by this:
\begin{itemize}
\item{The extremal vector describes the last component of the operator}
\item{The symbols $D_{j_1}^+ D_{j_2}^- D_{j_3}^+\ldots $ are translated to the composition of Dirac operators 
$D_{j_1} D_{j_2} \ldots $ acting on $g_2$}
\item{The symbols $D_{j_1}^+ D_{j_2}^- D_{j_3}^+\ldots D_{j_k}^\pm Y_{21}$ are translated to the composition of Dirac operators 
$D_{j_1} D_{j_2} \ldots$ acting on $-g_1$}
\end{itemize}

Now, we can easily determine the explicit form of the operator dual to $M_\liep(\xi)\to M_\liep(\nu)$. We know that
$\V_\xi\simeq \C\otimes S^+$, so the operator has only one component in this case. In Theorem \ref{thirdoperatorextrvect}
we found the extremal vector to be $D_1^- + D_2^- Y_{21}$ so we see that the operator is described by
$$
{h_1 \choose h_2} \mapsto D_1 g_2 - D_2 g_1.
$$

The same sequence of differential operators was derived in \cite{bluebook} by algebraic methods as 
a resolution of $2$ Dirac operators. We showed, however, that this operators are invariant with
respect to the action of $\sl(2,\R)\times \Spin(2n)$ and even more, by choosing the proper generalized conformal
weight, this operators are even $G$-invariant, where the action of $G$ includes translations
$((g_1(x,y), g_2(x,y)) \mapsto (g_1(x+u,y+v), g_2(x+u, y+v))$ and all actions of elements from $P$. 
This is an analogue of the invariance of the usual Dirac operator not only with respect to
$\Spin$ but also to all conformal transformations.

Note that we used the realization $\V_\nu\simeq (\C^2)^*\otimes S^-$ as $\lieg_0^{ss}$-module.  
But $(\C^2)^*\simeq (\C^2)$ as $\sl(2,\R)$-modules, so we could use the identification 
$\V_\nu\simeq \C^2\otimes S^-$ instead.
The difference is that for $e_1, e_2$ being a basis of $\C^2$, $e_1$ is the highest weight vector and 
$Y_{21} e_1=e_2$ (not $-e_2$). So, the operator $\cD'$ derived by this way would be $\cD_1'=\cD_2$
and $\cD_2'=-\cD_1$. The third operator would be $(h_1, h_2)\mapsto D_1 h_2 + D_2 h_1$.

Finally, let us remark that in case of the Dirac operator in $k$ variables, the first four vertices in
the BGG graph representing the weights 
\begin{eqnarray*}
&&\mmmm \lambda+\delta=\half[\ldots,5,3,1|\ldots,3,1]\\
&&\mmmm \mu+\delta=\half[\ldots,5,3,-1|\ldots,3,-1]\\
&&\mmmm \nu+\delta=\half[\ldots,5,1,-3|\ldots,3,-1]\\
&&\mmmm \xi+\delta=\half[\ldots,5,-1,-3|\ldots,3,1]\\
\end{eqnarray*}
are connected by arrows $\lambda\to\mu$, $\mu\to\nu$ and $\nu\to\xi$. The existence of a nonzero homomorphism
\gvmhom  and $M_\liep(\xi)\to M_\liep(\nu)$ was shown in Theorem $\ref{standard-nonstandard}$ and the
existence of a nonzero homomorphism $M_\liep(\nu)\to M_\liep(\mu)$ can be shown by proving
that the corresponding extremal vector is $D_{k-1}^+ D_k^- + D_k^+ D_k^- Y_{k,k-1} - y_{2n+k+1,k}$, similarly as in
Theorem \ref{2dirgvmevenextremal}, where $D_k^\pm$ and $D_{k-1}^\pm$ are defined analogously as
$(\ref{D+}), (\ref{D-})$. There is no difference in the computations. As $\lieg_0^{ss}$-modules,
$\V_\mu^*\simeq \C^k\otimes S^+$, but $\V_\nu^*$ is some more complicated representation of dimension $l\in\N$. If
we identify $\V_\nu^*\simeq \C^l\otimes S^+$, then the second order operator $\cD$ has $l$ components. 
From the above analysis, we can easily derive that if $\epsilon_l$ is
the highest weight vector in $\V_\nu\simeq (\C^l)^*$ (as $\sl(k,\R)$-module) 
then $l$-th component of the operator $\cD$ is
$$\cD_l (g_1,\ldots, g_k)=D_{k-1} D_k g_k - D_k D_k g_{k-1}$$ 
where $D_k, D_{k-1}$ are the Dirac operator in the $k$-th, resp. $(k-1)$'th, variable.

The composition $M_\liep(\nu)\to M_\liep(\mu)\to M_\liep(\lambda)$ is still zero but 
$M_{\liep}(\xi)\to M_\liep(\nu)\to M_\liep(\mu)$ is nonzero for $k>2$.

\end{document}